\documentclass{emsart}
\usepackage{makeidx}
\makeindex

\contact[candela@dm.uniba.it]{Anna Maria Candela,\\
Dipartimento di Matematica,\\ Universit\`a degli Studi di Bari,\\
Via E. Orabona 4, 70125 Bari, Italy}

\contact[sanchezm@ugr.es]{Miguel S\'anchez,\\
Departamento de Geometr\'{\i}a y Topolog\'{\i}a,
Facultad de Ciencias,\\
Universidad de Granada,\\ Avenida Fuentenueva s/n, 18071 Granada, Spain}

\numberwithin{equation}{section}

\newtheorem{theorem}{Theorem}[section]
\newtheorem{corollary}[theorem]{Corollary}
\newtheorem{lemma}[theorem]{Lemma}
\newtheorem{proposition}[theorem]{Proposition}

\theoremstyle{definition}
\newtheorem{definition}[theorem]{Definition}
\newtheorem{remark}[theorem]{Remark}
\newtheorem{example}[theorem]{Example}

\newcommand{\eps}{\varepsilon}
\newcommand{\wk}{\rightharpoonup}
\newcommand{\then}{\Longrightarrow}
\newcommand{\met}{\langle\cdot,\cdot\rangle }
\newcommand{\M}{M}
\newcommand{\mo}{M_0}
\newcommand{\cat}{\operatorname{cat}}
\newcommand{\N}{{\mathbb N}}
\newcommand{\Z}{{\mathbb Z}}
\newcommand{\R}{{\mathbb R}}
\newcommand{\LL}{{\mathbb L}}
\newcommand{\J}{{\cal J}}
\newcommand{\es}{{\cal S}}

\newcommand{\C}{{\mathbb C}}
\newcommand{\SSS}{{\mathbb S}}

\newcommand{\rf}[1]{\mbox{(\ref{#1})}}

\title[Geodesics in semi--Riemannian Manifolds]{Geodesics in semi--Riemannian Manifolds:
Geometric Properties and Variational\\ Tools}

\author[A.M. Candela, M. S\'anchez]{Anna Maria Candela\thanks{Supported by M.I.U.R. (research funds
ex 40\% and 60\%).} , Miguel S\'anchez\thanks{Partially supported by  MEC-FEDER Grant MTM2007-60731 and  J. Andal. Grant P06-FQM-01951.}}

\begin{document}
\hyphenation{re-pa-ra-me-tri-za-tion}

\begin{abstract}
Geodesics become an essential element of the geometry of a
semi--Riemannian manifold. In fact, their differences and
similarities with the (positive definite) Riemannian case,
constitute the first step to understand semi--Riemannian Geometry.
The progress in the last two decades has become impressive, being
especially relevant the systematic introduction of
(infinite--dimensional) variational methods.

Our purpose is to give an overview, from refinements of classical
results to updated variational settings. First, several properties
(and especially completeness) of geodesics in some ambient spaces
are studied. This includes heuristic constructions of compact
incomplete examples, geodesics in warped, GRW or stationary
spacetimes, properties in surfaces and spaceforms, or problems on
stability of completeness.

Then, we study  the variational framework, and focus on two
fundamental problems of this approach, which regards geodesic
connectedness. The first one deals with a variational principle
for stationary manifolds, and its recent implementation inside
Causality Theory. The second one concerns orthogonal splitting
manifolds, and  a reasonably self--contained development is
provided, collecting some steps spread in the literature.
\end{abstract}



\maketitle


\newpage
\tableofcontents

\section{Introduction}

Geodesics become an essential ingredient of semi--Riemannian
Geometry, as so they are for the (positive definite) Riemannian
one. In the indefinite case, two important difficulties arise:

{\sl (1)} there are simple questions which become elementary in the
Riemannian case but they are open in general. For some of them, it
is not clear which type of techniques will work -- say, just  a
simple brilliant idea or a full specific theory;

{\sl (2)} in the Lorentzian case, {\em Causality} is a specific theory
fully adapted to its particularities. So, other sophisticated
theories which may work in the Riemannian case, must be
implemented carefully in the framework of Causality.
Higher order indexes may be intractable, as they do not admit such a tool.
\smallskip

Accordingly, here our aim is twofold. The first one, developed in
Section \ref{s3}, is to provide a brief overview on geometric
properties of geodesics in the general semi--Riemannian case,
stressing the differences and similarities with the Riemannian
one. Even though many of them are simple items, they may involve
very different techniques, and pose open questions at a primary
level, whose strategy of solution is hard to predict. For example,
as far as we know, the following question is open: {\em must a
compact semi--Riemannian manifold which is globally conformal to a manifold
of constant curvature be complete?} Here, the reader will find related
results such as: {\sl (i)}~compact Lorentzian manifolds of constant
curvature are complete (see Theorem \ref{t-markus}), {\sl
(ii)}~lightlike completeness is conformally invariant in the
compact case (see Theorem \ref{t-comp-conf-luminosa}), {\sl
(iii)}~in locally symmetric spaces the three types of causal
completeness are equivalent (see Theorem \ref{t-lafuente}), {\sl
(iv)}~there are difficulties to find independence of causal
completeness in the compact case (see Remarks \ref{r-cr},
\ref{r-IndCom}), and {\sl (v)}~there exists a relation among: {\sl (a)} the
completeness of conformally related compact semi--Riemannian
manifolds, {\sl (b)} the causal independence of completitudes and 
{\sl (c)} the
completeness of warped products (see Remark \ref{r-dna-wac-coca}).
These five items have been studied by means of very different
techniques, with quite different levels of sophistication. 
It is not clear which techniques can solve
the open question above but, at any case, the answers to this type of
problems will yield an advance in our knowledge of
semi--Riemannian Geometry. In order to be specific, most of the
overview in Section \ref{s3} concerns geodesic completeness.
Nevertheless, the reader is introduced in the behaviour of
geodesics of different ambients: locally symmetric, homogeneous,
spaceforms, warped products, stationary... The reader is expected
to have a basic knowledge of Riemannian Geometry, and we hope
he/she will learn which Riemannian results and tools still hold or
can be adapted to the indefinite case. So, we choose a heuristic
approach in topics such as the counterexamples to completeness
(see Subsection \ref{sub3}) or the role of curvature in
singularity theorems (see Subsection \ref{sub5}). Some of the
topics on geodesics covered here, can be complemented with other
properties of semi--Riemannian manifolds, as those ones surveyed in \cite{RP}.

Our second main aim, developed in Section \ref{s4}, is to give an
overview of the infinite--dimensional variational setting  for
geodesics in Lorentzian manifolds introduced at the end of 80's by Benci,
Fortunato, Giannoni and Masiello (see, e.g., \cite{BF90,BFG91,GiaMa91}). Even though this
approach has been explained in book format by Masiello \cite{Ma},
remarkable progress has been carried out since then, even in the foundations of the theory.  First, as commented above, the interplay
between the variational theory and Causality is necessary for the
development of the approach in all its extent. This full interplay
has been achieved only very recently in \cite{CFS06}, and just
for the most simple ``fundational problem'', i.e, {\em
geodesic connectedness of static and stationary manifolds}. We
think that further results in this direction may be of big
interest for both, the variational analyst's and Lorentzian
geometer's viewpoints. On the other hand, regarding the ``second
fundational problem'', i.e., {\em geodesic connectedness in
splitting type manifolds},  we emphasize that, even though the
core of the approach is explained in  references such as the book
\cite{Ma}, the full details were filled later, and are spread in
several  papers published along the last decade. So, we explain the full approach
in a reasonably self--contained way.

The variational approach starts by recalling that, whenever
$(\M,g)$ is a semi--Riemannian manifold, its geodesics can be
found as critical points of the {\sl action
functional}\index{functional!action}
\begin{equation}\label{action}
f(\gamma) = \int_0^1 g(\gamma)[\gamma',\gamma']\ ds,
\end{equation}
where $\gamma: [0,1] \to \M$ is any curve in a suitable manifold of
functions (for more details, see later on). If $\M$ is a
Riemannian manifold, i.e., its metric $g$ has index $s_0 = 0$, the
functional $f$, now named {\sl energy
functional}\index{functional!energy}, can be directly studied by
means of classical variational tools ---as it is positive, hence,
bounded from below (see Section \ref{secriemannian}).

But, more in general, if the semi--Riemannian metric $g$ has index
$s_0 \ge 1$ ($s_0<n_0$) the corresponding functional $f$ is strongly
indefinite (i.e., unbounded both from above and from below, even
up to compact perturbations) with critical points having infinite
Morse index. Thus, different methods and suitable ``tricks'' are
needed and, at least in the Lorentzian case, we can distinguish
two different variational approaches which allow one to overcome
such a problem  (see the book \cite{Ma} or the survey
\cite{Sa01}):
\begin{itemize}
\item[{\sl (a)}] to transform the indefinite problem on a
Lorentzian manifold in a subtler (hopefully bounded from below)
problem on a Riemannian manifold; \item[{\sl (b)}] to study
directly the strongly indefinite functional $f$ but by making use
of suitable (essentially finite--dimensional) ``approximating''
methods.
\end{itemize}

Furthermore, the choice of the right manifold of curves where
functional $f$ is defined, depends on the different geometric
problem studied, and impose specific ``boundary'' properties:
\begin{itemize}
\item geodesics joining two fixed points (see Section
\ref{secriemannian}), \item closed geodesics (see, e.g.,
\cite{CS00,Ma93} and references therein), \item geodesics
connecting two given submanifolds (see, e.g., \cite{Mo96} or also
\cite{CMS00} and references therein), \item {\sl $T$--periodic
trajectories} (see, e.g., the pioneer papers \cite{BF90,Gr89}, the
survey \cite{Ca97} or \cite{BGS-annali}, and references therein).
\end{itemize}
Moreover, one can use variational tools also in order to find geodesics 
with a prescribed causal character (for the timelike ones, see, e.g., \cite{BCFS03}
and also references in \cite{Ca00}, for the lightlike ones, see the Fermat
principle in \cite{FGM95,FM95}).

The difficulty of the interplay between variational methods and
Causality makes that, in some particular cases, non--variational
approaches may yield more accurate results (see, for example,
\cite{FS-jde,FS-jmp} for geodesic connectedness, \cite{Sa06} for
closed geodesics, or \cite{Sa-nonlinearstatic,Sa-procamsperiodictraject}
for $T$--periodic trajectories). Nevertheless,
one may expect that, as we have commented for the stationary case,
the strongest results will be obtained when the variational tools
are fully implemented in Causality Theory.

Here, for simplicity, we want just to outline how to manage the
previous approaches {\sl (a)}, {\sl (b)} in a model case. More
precisely, after the introductory Riemannian case in Subsection
\ref{secriemannian}, we study geodesic connectedness:
in stationary spacetimes (method {\sl (a)}) in Subsection
\ref{secstatic}, and in orthogonal splitting spacetimes (method {\sl
(b)}) in Subsection \ref{secsplitting}.


\section{First properties}\label{s2}

Most of the material in this preliminary section is well--known.
In the Subsection \ref{sub1.1},  essentially,  notation and first
definitions are given. In Subsection \ref{sub1.2}, first we
summarize maximizing properties of timelike and lightlike
geodesics in the Lorentzian case (but we refer to \cite{MS} in
these proceedings for much more detailed results). Spacelike
geodesics, including some remarkable properties of their conjugate
points, are also considered here. The transformation of the
Levi--Civita connection and geodesics under conformal changes
(including implications for lightlike geodesics) are studied in
Subsection \ref{sub1.3}, and will be widely used throughout
Section \ref{s3}.

\subsection{Preliminaries} \label{sub1.1} Usual notation and
conventions, essentially compatible with standard books as
 \cite{BEE} or \cite{O}, will be used. We refer also to
 these references for detailed proofs of the basic
properties collected in the present subsection. As we have pointed
out, the interplay with Causality will be essential in the
Lorentzian case and, so, the contribution to these proceedings
\cite{MS} will be frequently invoked for background material on
Causality.

\begin{definition} A {\em semi--Riemannian manifold}\index{manifold!semi--Riemannian}
is a  smooth manifold $\M$,
of  dimension $n_0\geq 1$, endowed with a non--degenerate
metric $g:M \to T^{*}M\otimes T^{*}M$ of constant index $s_0$.

In the case $s_0=0$ the manifold is called
{\em Riemannian} \index{manifold!Riemannian} (or positive definite), if $s_0=n_0$
is called {\em negative definite}\index{manifold!negative definite},
if $0<s_0<n_0$ is called {\em indefinite}\index{manifold!indefinite}
and, in this case,  if $s_0=1$ ($n_0>1$) is
called {\em Lorentzian}\index{manifold!Lorentzian}.
\end{definition}

By {\em smooth}\index{manifold!smooth} we mean $C^{r_0}$; where
 $r_0=+\infty$ will be assumed
for simplicity for all the elements, except when otherwise is
stated explicitly. $M$ will be also assumed (in addition to
Hausdorff, as usual) connected, except if otherwise is specified;
thus, the constancy of $s_0$ can be deduced of the
non--degeneracy\footnote{This non--degeneracy also implies paracompactness and,
as $M$ is connected, the second axiom of numerability for the
topology (see \cite[Section 2.1]{MS} for references).} of $g$. The notion of
``causal character'' for a tangent vector $v\in TM$ in Lorentzian
Geometry, which comes from General Relativity, is
extended here and, so, $v$ is called timelike (resp., lightlike;
causal; spacelike) \index{vector!causal} depending on  if
$g(v,v)<0$ (resp. if $g(v,v)=0$ and $v \ne 0$; $v$ is either
timelike \index{vector!timelike} or lightlike,
\index{vector!lightlike} i.e., $g(v,v)\le 0$ and $v \ne 0$;
$g(v,v)>0$). Vector 0 does not lie in any of these types, even
though sometimes is useful to regard it as spacelike (see \cite{O});
according to \cite{MS}, a null vector \index{vector!null}  will be
either timelike or 0. These causal characters are naturally
extended to curves and submanifolds. We will put
$|v|=\sqrt{|g(v,v)|}$ and the {\em length} of a curve\index{curve!length}
$\gamma$ is the integral of $|\gamma'|$.

Recall that the existence and uniqueness of the Levi--Civita
connection $\nabla$ \index{Levi--Civita connection}\index{$\nabla$} depends only on
the non--degeneracy of the metric $g$ and, therefore, it can be
deduced in the semi--Riemannian case as in the Riemannian one.
Associated to $\nabla$ there is a covariant derivative $D/ds$\index{$D/ds$}
and parallel transport (which is an isometry, too).
Moreover, a smooth curve
$\gamma:I \to \M$ is a {\sl geodesic}\index{geodesic} if its velocity
is parallel or, equivalently, its acceleration vanishes:
\[
\frac{D\gamma'}{ds}(s) = 0\quad \hbox{for all $s \in I$.}
\]
If $\gamma$ is a geodesic, $g(\gamma',\gamma')$ is a constant. Thus,
the {\sl causal character of a geodesic}\index{geodesic!causal character}
is defined as timelike, lightlike etc. according to the
constant causal character of its velocity. By taking an {\em
orthonormal basis} \index{orthonormal basis} $B_p=(v_1,\dots ,
v_{n_0})$ at $T_pM$, that is, satisfying
\[
g_p(v_i,v_j)= \epsilon_i \delta_{ij}, \quad \epsilon_i =  -1, \;
\hbox{if } i=1,\dots , s_0, \quad \epsilon_i = 1, \; \hbox{if } i=
s_0+1, \dots, n_0 ,
\]
one obtains a natural isometry with $\R^{n_0}_{s_0}$\index{$\R^{n_0}_{s_0}$}, that is,
$\R^n$ endowed with the natural product of index $s_0$. In the
Lorentzian case, the semi--Riemannian manifold $\R^{n_0}_{s_0}$
will be the {\em Lorentz--Minkowski} spacetime\index{spacetime!Lorentz--Minkowski}, denoted
by $\LL^{n_0}$\index{$\LL^{n_0}$}.

As for any affine connection, the differential at $0$ of the
exponential map\index{exponential map}, $\exp$\index{$\exp$},
$$
(d\exp_p)_0: T_0(T_pM) \rightarrow T_pM
$$
is the identity up to natural identifications. Thus, any point
$p\in M$ admits a {\em starshaped neighborhood}
$U$\index{neighborhood!starshaped} ($\exp_p^{-1}$ becomes well
defined on $U$ and its image is a starshaped neighborhood of $0\in
T_pM$ in the usual sense), {\em a normal neighborhood}
\index{neighborhood!normal} (any chart $(U,\varphi \circ
\exp_p^{-1})$, where $U$ is starshaped and $\varphi:
T_pM\rightarrow \R^{n_0}_{s_0}$  a linear isometry) and a convex
neighborhood \index{neighborhood!convex} (a normal neighborhood of
all its points), see \cite[Proposition 5.7]{O}. {\em Gauss Lemma}
\index{Gauss Lemma} also makes sense, and can be proved in a
similar way as in the Riemannian case (see \cite[Theorem
10.18]{BEE}, \cite[Lemma 5.1]{O}): {\em if $p\in M$, $0\neq x\in
T_pM$ and $v_x, w_x \in T_x(T_pM)$ with $v_x$ collinear with $x$,
then:
$$
g_p((d\exp_p)_x(v_x), (d\exp_p)_x(w_x)) = \langle v_x, w_x \rangle
$$
where $\langle \cdot , \cdot \rangle $ is the scalar product in
$T_x(T_pM)$ naturally induced by $g_p$.}

The notion of {\em conjugate point} $q$ \index{conjugate point} of
a point $p$ along a geodesic $s\mapsto \gamma(s) =\exp_p(sv)$,
$v \in T_p M$ (i.e., $q= \exp_p (s_q v)$ where $s_q v$ is a
critical point of $\exp_p$) and its multiplicity
(dimension of the kernel of $d\exp_p$ at $s_qv$), is a natural extension of the usual
Riemannian definition. Nevertheless, their properties may be very
different, as we will see below. {\em Jacobi equation} (and its
solutions, the {\em Jacobi fields}) \index{Jacobi!equation}
\index{Jacobi!field} is also defined in the semi--Riemannian case
as a formal extension of the Riemannian one, that is: $DJ'/ds =
-R(J,\gamma')\gamma'$, where $R$ is the curvature tensor\index{curvature tensor}
under the convention $R(X,Y) = [\nabla_X,\nabla_Y]-\nabla_{[X,Y]}$\index{$R(X,Y)$}
for all vector fields $X,Y \in \mathfrak{X}(M)$\index{$\mathfrak{X}(M)$}.
Moreover, Jacobi fields
on $\gamma$ correspond to variational fields through longitudinal
geodesics, and  the multiplicity of a conjugate point
\index{conjugate point!multiplicity}
$q$ to $p$ along $\gamma$ is equal to the dimension of the space of Jacobi
vector fields on $\gamma$ which vanishes on $p, q$ (see \cite[Propositions
8.6, 10.10]{O}).

The {\em index form} \index{index form} of a non--lightlike
geodesic $\gamma: [a,b]\rightarrow M$ is analogously defined as\index{$I(V,W)$}:
\[
I(V,W)= \nu \int_a^b
\left(g(\frac{DV}{ds},\frac{DW}{ds})-g(R(V,\gamma')\gamma',W)\right)ds
\]
where $\nu$ is a non--null constant which depends on the speed
$|\gamma'|$ and the causal character of $\gamma$ ($\nu$ is chosen
with different values by different authors), and $V, W$ belong to
$\mathfrak{X}^\perp(\gamma)$\index{$\mathfrak{X}^\perp (\gamma)$},
i.e., the space of all the (piecewise) smooth vector fields on
$\gamma$ with vanishing endpoints pointwise orthogonal to
$\gamma'$. For lightlike geodesics, the index form $I$ can be
defined formally as above, but recall: {\sl (a)} the multiplicity
of a conjugate point along a non--lighlike geodesic $\gamma$ can
be obtained as the dimension of the radical (nullspace)  of $I$
(see \cite[Corollary 10.12]{O}), {\sl (b)} in order to extend this
result to the lightlike case, $I$ is replaced by the {\em quotient
index form} \index{index form!quotient} $\bar I$. With this aim,
one first define the quotient space $\overline{\mathfrak{X}}^\perp
(\gamma )$ \index{$\overline{\mathfrak{X}}^\perp (\gamma)$}
obtained by identifying each two $X,Y\in \mathfrak{X}^\perp
(\gamma)$ if $X-Y$ is pointwise collinear to $\gamma'$. Then,
$\bar I$\index{$\bar I(X,Y)$} is defined as the bilinear form on
$\overline{\mathfrak{X}}^\perp (\gamma)$ obtained by inducing
naturally $I$ (see \cite{MS}). Even more, in the Lorentzian case,
a Morse Index Theorem holds for both, the index form $I$ on
timelike geodesics and the quotient index form on lightlike
geodesics (see \cite[Chapter 10]{BEE}).

\subsection{Special properties of geodesics in spacetimes
depending on their causal character}\label{sub1.2} We will mean by
{\em co--spacelike} \index{geodesic!co--spacelike} any geodesic
such that the orthogonal of its velocity is a spacelike subspace
at each point, that is: all the geodesics in the Riemannian case
and timelike geodesics in the Lorentzian one. They present similar
properties of extremization \index{geodesic!extremizing
properties} (minimization or maximization) for the ``distance''
$d$ associated to the metric $g$. By such a distance we mean the
true canonical distance associated to $g$, if $g$ is Riemannian
(i.e., (\ref{distanza}) below), but the time--separation in the
Lorentzian one. Concretely (following \cite{MS}), if $g$ is
Lorentzian, we assume in the present subsection that it admits a
{\em time--orientation} \index{time--orientation} (i.e., a
continuous choice of causal cones, which will be called {\em
future cones}). $(M,g)$, with the additional choice of a
time--orientation, is a {\em spacetime}\index{spacetime}, and the
 {\em time--separation} (or Lorentzian distance)\index{time--separation
(or Lorentzian distance)}
is defined for any $p,q\in M$ as the supremum $d(p,q)$  of the
lengths of the future--directed causal curves from $p$ to $q$ 
(or 0 if no such a curve exists).
A spacetime $\M$ is {\em globally
hyperbolic}\index{spacetime!globally hyperbolic} if there exists a
(smooth) spacelike Cauchy hypersurface\index{Cauchy hypersurface}
$\es$ in in $M$ (i.e., a subset which is crossed once by any inextendible 
timelike curve). This property will have important implications 
for both, the spacetime and its geodesics (see \cite[Section 3.11]{MS}). 

As  a first definition, a causal (resp. Riemannian) geodesic
with endpoints $p, q$ will be called {\em maximizing} (resp.
minimizing) \index{geodesic!maximizing causal}
\index{geodesic!minimizing Riemannian} when its length is
equal to the maximum  between $d(p,q)$ and $d(q,p)$ (resp.
$d(p,q)$). Recall that, even though $d$ is not symmetric in a
general spacetime,  the notation will be simplified below writing
just $d(p,q)$ for the maximum between $d(p,q)$ and $d(q,p)$. The
extremizing properties of co--spacelike geodesics for $d$ can be
deduced by means of a standard study of the index form. Even more,
many of them can be extended to lightlike geodesics  as a specific
case.

Next, we sketch a very rough summary. We emphasize that, in the
Lorentzian case,  these properties depend heavily on the causal
structure of the spacetime and, so, we refer to the contribution
\cite{MS} for detailed definitions and properties. Here, our
purpose is just to stress the similitudes and differences of the
extremizing properties of timelike geodesics and Riemannian
geodesics, and how they can be extended to the lightlike case (see
\cite[Theorems 2.26, 2.27]{MS} for precise statements). So, in the next
paragraphs, the extremizing properties will be stated first for
timelike geodesics, and the corresponding Riemannian property is
also pointed out. Then, the lightlike case will be studied.
Finally, some comments on the spacelike case will be added.

\paragraph{Timelike  and co--spacelike geodesics.} It is well known that
conjugate points along a timelike (resp.
Rieman\-nian) geodesic in a Lorentzian (resp. Riemannian) manifold
{\em cannot have points of accumulation}. Even more:
\begin{itemize}
\item[(1)]Any timelike geodesic locally maximizes the time separation $d$ in a similar
way as any Riemannian geodesic locally minimizes its corresponding distance
$d$. Nevertheless, there are two important differences:
\begin{itemize}
\item Riemannian geodesics locally minimize the lengths of all the
smooth curves connecting two fixed points $p, q$, while
the timelike geodesics maximize only the lengths of the causal curves connecting
$p, q$;
\item the Riemannian minimizing property holds for both, the
restriction $d|_U$ of the distance $d$ on $M$ to any suitably small
neighborhood $U$, and the distance $d_U$ naturally associated to
the restriction $g_U$ of the metric  to $U$. 
In general, the maximizing property for timelike
geodesics  holds only for $d_U$. It holds for $d|_U$ in 
{\em strongly causal spacetimes}  (in particular, in
globally hyperbolic ones). Recall that, for example, there are
spacetimes (the {\em totally vicious} ones) with $d(p,q)= \infty$ for all $p,q$, that is, $d|_U \equiv \infty$ for all$U$.
\end{itemize}
\item[(2)] Let $\gamma:[a,b]\rightarrow M$ be any timelike (resp.
Riemannian) geodesic which connects two points $p,q\in M$
non--conjugate along $\gamma$. Then: {\em $\gamma$ will have
strictly maximum (resp. minimum) length among neighboring
curves
connecting $p,q$ (obtained by means of a variation with fixed
endpoints, and up to a reparametrization) if and only if there is
no conjugate point to $p$ along $\gamma$.}
\item[(3)] If a timelike (resp. smooth) curve
$\gamma:[a,b]\rightarrow M$ maximizes (resp. minimizes) the
time--separation (resp. distance) $d$ then, up to a
reparametrization, it is a geodesic without conjugate points
except at most $p=\gamma(a)$ and $q=\gamma(b)$.
\item[(4)] In a globally hyperbolic spacetime (resp. a complete
Riemannian manifold) each two {\em chronologically related} points
(resp. each two points) can be connected by means of a
$d$--maximizing (resp. minimizing) timelike (resp. Riemannian)
geodesic. Even more, each inextendible timelike geodesic $\gamma:
[0,b)\rightarrow M$ maximizes in some subinterval $[0,c)\subset
[0,b)$,  $c\in (0,b]$. If $c<b$ then $\gamma(c)$ is called the
{\em cut point}\index{cut point} of $p=\gamma(0)$ along $\gamma$, and cannot appear
beyond the first conjugate point.\\
Moreover, such results can be extended in the Lorentzian case to
timelike homotopy classes \index{homotopy classes!timelike} (i.e.,
classes of homotopy where the longitudinal curves are timelike).
So, fixing $p,q\in M$, if the spacetime is globally hyperbolic
(resp. the Riemannian manifold is complete) then, in each timelike
homotopy (resp. homotopy) class of curves with fixed endpoints $p,
q$,  there is one connecting timelike (resp. Riemannian) geodesic,
with maximum (resp. minimum) length among the curves in that
class.
\end{itemize}

\paragraph{Lightlike geodesics.} As in the timelike case, conjugate points
 can never accumulate on a Lorentzian lightlike geodesic. Moreover, they 
 cannot appear neither in mani\-folds of constant curvature nor in
dimension 2 (see \cite[Proposition 2.34]{MS}). The previous four points also
hold, with the following modifications:
\begin{itemize}
\item[(1')] As in (1), any lightlike geodesic $\gamma$ maximizes
$d$ locally. But, now, this means  that, for any $p=\gamma(s)$
there exists a neighborhood $U\ni p$ (say, any convex
neighborhood) such that if $q\in U$ lies on $\gamma$ then no other
causal curve contained in $U$ connects $p,q$.
\item[(2')] This property is analogous to (2) but recall that,
when $\gamma$ is timelike then curves close to $\gamma$ are causal
(and timelike), too. The analogous property does not hold for
lightlike curves. So, a lightlike geodesic
$\gamma:[a,b]\rightarrow M$ which connects two non--conjugate
points $p,q\in M$ is, up to a reparametrization, the unique causal
curve  among neighboring causal curves connecting $p,q$, if and only if
there is no conjugate point to $p$ along $\gamma$ (see \cite[Section 2]{MS}).
\item[(3')] If a causal curve $\gamma:[a,b]\rightarrow M$
maximizes  the time separation  $d$ then, up to a
reparametrization, it is a geodesic without conjugate points
except at most the endpoints, and then $\gamma$ is lightlike if
and only if $d(p,q)=0$.
\item[(4')] All the assertions in (4) hold just replacing
``chronologically related'' by ``causally related'',  and
``timelike'' geodesics or homotopy class by ``causal''
one.\index{homotopy classes!causal} Nevertheless, causal homotopy
classes have remarkable specific properties. For example (see
\cite{MS-cmp} for a detailed study):\\
{\em  Let $\gamma: [0,b)\rightarrow M $ be a lightlike geodesic
with a cut point $\gamma(\lambda_c)$, $\lambda_c \in(0,b)$. If
$\gamma(\lambda_c)$ is not a conjugate point, then:
\begin{itemize}
\item[(1)] no other lightlike geodesic which connects $\gamma(0)$ and
$\gamma(\lambda_c)$ is causally homotopic to $\gamma$;
\item[(2)] if $(M,g)$ is globally hyperbolic, there exists at least
another lightlike geodesic $\hat \gamma$ (necessarily non--causally
homotopic to $\gamma$) which connects $\gamma(0)$ and
$\gamma(\lambda_c)$. \end{itemize}} \end{itemize}
(For more properties of the cut locus see \cite[Chapter 9]{BEE}; Morse
theory for lightlike geodesics was introduced by Uhlenbeck
\cite{Uh}).

\paragraph{Spacelike geodesics.}  For spacelike geodesics in
Lorentzian manifolds (of dimension $n_0\geq 3$),  as well as for
non co--spacelike geodesics in any semi--Riemannian manifold (of
index $s_0<n_0-1$), no maximization nor minimization properties
hold, even though such geodesics are still critical points of the
action functional (\ref{action}). \\
Furthermore, as a difference with both, the Riemannian case and the case of
 causal geodesics
in a Lorentzian spacetime, conjugate points along a spacelike geodesic
may present
accumulation points and even cover a compact interval, so the Morse Index Theorems
as in \cite[Section 10]{BEE} cannot be used.
In fact, in \cite[Section 11]{He} the author constructs a spacelike geodesic
in a Lorentzian spacetime which has a continuum of conjugate points
while, later on, in \cite[Theorem 3.4]{PiTa03} the authors prove that
{\sl ``taken any compact subset $F \subset\ ]a,b]$ there exists
a 3--dimensional Lorentzian manifold $(\M,g)$ and a spacelike geodesic
$\gamma: [a,b] \to \M$ such that $\gamma(t)$ is conjugate to $\gamma(a)$
along $\gamma$ if and only if $t \in F$.''}
Anyway, suitable index theories can be developed also in this setting
in order to prove Morse--type theorems applicable to the spacelike case
(see \cite{He,PT} or also \cite{PPT} for an index theory in more general
semi--Riemannian manifolds).

\subsection{Conformal changes}\label{sub1.3}
For any indefinite semi--Riemannian manifold, the lightlike vectors
determine the conformal class of the metric (this is a consequence
of a simple algebraic result, see \cite[Proposition 2.6]{MS}). Thus, two
indefinite semi--Riemannian metrics $g, g^*$ on $M$ are {\em
pointwise conformal} \index{metric!pointwise conformal}, i.e.,
$g^*=\Omega g$ for some (non--vanishing) smooth function $\Omega
\in C^\infty(M)$, \index{$C^\infty(M)$} if and only if they have
equal lightlike vectors. In what follows, we will asume:
$$
g^*=\Omega g, \quad \quad \hbox{with} \quad \Omega>0, \quad
\Omega=e^{2u}, \quad u\in C^\infty(M).$$ A straightforward
computation from Koszul formula yields the relation between the
corresponding Levi--Civita connections $\nabla, \nabla^*$:
\[
\nabla^*_X Y = \nabla_X Y + X(u) Y + Y(u) X - g(X,Y) \nabla u ,
\]
where $\nabla u$ denotes the $g$--gradient of $u$, and the equality
holds for any vector fields $ X, Y\in \mathfrak{X}(M)$.
Consequently, one obtains a relation
between the corresponding covariant derivatives $D/dt, D^*/dt$ and
the following equality between the accelerations of any curve
$\gamma: I\subset \R \rightarrow M$:

\begin{equation} \label{e-cambioconfderivcovar}
\frac{D^*\gamma'}{dt} = \frac{D\gamma'}{dt} + 2 du (\gamma'(t))
\cdot \gamma' - g(\gamma',\gamma')\nabla u.
\end{equation}
In the particular case that $\gamma$ is a lightlike $g$--geodesic,
one has:

$$\frac{D^*\gamma'}{dt}  = f(\gamma) \gamma',
\quad \hbox{with} \quad f(\gamma) \equiv 2 \frac{d(u\circ \gamma)}{dt}
\quad  \hbox{on} \; I.
$$
That is, $\gamma$ is a  pregeodesic \index{pregeodesic} (i.e.,
geodesic up to a reparametrization) for $g^*$, and the concrete
parametrization as a geodesic can be written as $\tilde\gamma : J
\rightarrow M$ with
\begin{equation} \label{e-reparam-pregeodlum}
\tilde\gamma (s) = \gamma (t(s)), \quad
 s'(t)= s_0' e^{\int_{t_0}^{t}f(\sigma) d\sigma} \equiv C \,
 e^{2u(\gamma(t))}, \quad \hbox{for all $t\in I$.}
\end{equation}
This relation shows that the set of all the lightlike pregeodesics
is a conformal invariant. The result can be sharpened, by showing
that also their conjugate points and  multiplicities are
conformally invariant (see \cite[Theorem 2.36]{MS}). Moreover, as a simple
consequence of the possible maximal domains of definition $I, J$,
we have the following completeness result.

\begin{theorem} \label{t-comp-conf-luminosa}
Let $g$, $g^*= \Omega g, \Omega=e^{2u},$ two indefinite
semi--Riemannian metrics on the same manifold $M$. Let $\gamma:
I\subset \R \to \M$ be an inextendible lighlike geodesic for $g$, and
$\tilde\gamma: J\subset \R \to \M$ a reparametrization as an
inextendible geodesic of $g^*$.
\begin{itemize}
\item[(1)] If $ \inf(\Omega) >0$ and $\gamma$ is complete ($I=\R$), then
$\tilde\gamma$ is complete.
\item[(2)] If $ \sup(\Omega ) <0$ and $\gamma$ is incomplete, then
$\tilde\gamma$ is incomplete.
\end{itemize}
In particular, if $M$ is compact then $g$ is lightlike
complete if and only if so is $g^*$.
\end{theorem}

\begin{proof} By using the explicit reparametrization
(\ref{e-reparam-pregeodlum}), one has:\\
{\sl Case (1)} As $u$ is also lower bounded, $|ds/dt| \geq \epsilon >0$
for some $\epsilon >0$. As $I=\R$ the image  $J$ of $s(t)$ also covers all $\R$.\\
{\sl Case (2)} Analogously,  $(0<)ds/dt < N$, for some  $N>0$. Thus,
if, say   $I=(a,b)$ with $b < +\infty$, then the image of  $s(t)$
cannot reach the value $s(t_0) + N (b-t_0)$.
\end{proof}

\section{An overview on geodesics in different ambient manifolds and geodesic
completeness} \label{s3}

\subsection{First results} \label{sub1} Here, our main goals are: 
{\sl (a)} to show an example of independence of causal
completitudes, {\sl (b)} to prove dependence in the locally symmetric
case (see Theorem \ref{t-lafuente}), {\sl (c)} to give a criterion on
completeness, as an alternative way to Hopf--Rinow Theorem (see Proposition
\ref{pcomplcriterion}), and (d) to pose the problem of
independence of completeness in the compact case (see Remark
\ref{rincomp-nulldirect}).

\paragraph{Independence of completeness.}
In general, for an indefinite manifold, the term ``completeness''
means just geodesic completeness, as there is no any distance
canonically associated to the metric.
Nevertheless,\index{completeness (spacelike, lightlike and
timelike)} as there exist spacelike, lightlike and timelike
geodesics, one can speak on {spacelike, lightlike and timelike
completeness}, depending on the type of geodesics which are
complete. There are explicit examples by Kundt, Geroch and Beem
which show the full logical independence of the three types of
completeness (see \cite[pp. 203]{BEE} for detailed references).
And, in fact, it is easy to construct an example which is spacelike
incomplete and both, timelike and lightlike complete (see Theorem
\ref{twarpedcompl} for others).

\begin{example}
Concretely, consider in $\R^2$ the Lorentz metric
\[
g^*=\ e^{2u(x,y)}(dx^2-dy^2)
\]
conformal to the usual $g_0=dx^2-dy^2$,
where $u: \R^2 \rightarrow \R $ satisfies:
\begin{itemize}
\item[{\sl (i)}] $\partial_yu(x,0) = 0$ (in particular $u$ can be chosen
$x$--axis symmetric);
\item[{\sl (ii)}] it vanishes outside the horizontal
strip $S =\{(x,y)\in \R^2: |y|\le 1\}$;
\item[{\sl (iii)}]
$\displaystyle\int_{-\infty}^{+\infty} e^{u(x,0)}dx < +\infty$.
\end{itemize}
Now, taking into account formula (\ref{e-cambioconfderivcovar}),
one checks easily that the natural reparametrization of the
$x$--axis is a $g_0$--pregeodesic (by using condition {\sl (i)}), and
it is incomplete when parametrized as a (spacelike) geodesic
according to (\ref{e-reparam-pregeodlum}) (by using {\sl (iii)}).
Nevertheless, any causal geodesic $\gamma$ is complete; in fact,
if $\gamma(s_0)$ belongs to the strip $S$, $\gamma$ will leave $S$
because it cannot remain imprisoned \index{curve!imprisoned}
\index{curve!partially imprisoned} in the compact subset
$J^+(\gamma(s_0)) \cap S$ (see in this proceedings \cite[Section
3.6.2]{MS} for the notions of imprisoned and partially imprisoned
curve); this is also easy to check directly because $\gamma$ can
be reparametrized as a curve $y\mapsto (x(y),y)$ with $|dx/dy|\leq
1$) and, outside of $S$, $\gamma$ is a geodesic of $\LL^2$.
\end{example}

\begin{remark}\label{r-IndCom}
(1) {\em Completeness in any of the three causal types implies
inextendibility}. That is, if a  semi--Riemannian manifold $(M,g)$
is {\em extendible}, \index{manifold!extendible semi--Riemannian}
i.e., it is isometric to an open subset included strictly in
another (connected) semi--Riemannian manifold $(\tilde M, \tilde
g)$, then it is incomplete in the three causal senses, spacelike,
lightlike and timelike. In fact, recall that any point in the
boundary $\tilde p\in \partial M$ can be joined with a point $p$
of $M$ by means of a broken geodesic $\gamma$ of any causal
type\footnote{This property is needed only for some $p$  (thus,
one can choose $p$ in a convex neighborhood of $\tilde p$ in
$\tilde M$) but it holds for any $p\in M$.}. Looking $\gamma$ as a
curve starting at $p$, each geodesic piece will be included in
$M$, i.e., we have the contradiction $\tilde p\in M$.

(2) As far as we know,  {\em there are no explicit examples which
show some independence of any of the three causal type of
completeness for compact $M$} (see also Remarks
\ref{rincomp-nulldirect} and \ref{rcomplcriterion}(2) below).
Nevertheless, Carri\`ere and Rozoy \cite{CR} have suggested that
one such example will exist in a torus (see Remark \ref{r-cr}
below).
\end{remark}

\paragraph{Dependence in locally symmetric spaces.}  Recall that
 in a {\em symmetric}
 semi--Riemannian manifold $(M,g)$  there exists a
 {\em global symmetry} \index{manifold!globally symmetric} at each $p\in M$ (an isometry
$I: M\rightarrow M$ which fixes $p$ with $dI_p =-$Id, where Id is
the identity at $T_pM$) and in the {\em locally symmetric}
\index{manifold!locally symmetric} case this symmetry can be found
only in a neighborhood of each $p$. Locally symmetric manifolds
can be characterized as those with a parallel curvature tensor
and, thus, they include constant curvature ones.

The following result is due to  Lafuente \cite{La} (point (2)) and
 Furness and Arrowsmith \cite{FA} (point (3))\footnote{Elementary results
extending naturally the Riemannian ones,  such as {\sl ``a
semi--Riemannian manifold is complete if and only if so it is its
universal covering''}, or more refined versions such as  \cite[Corollary
7.29]{O}, will be used without further mention along the remainder of this article.}.

\begin{theorem} \label{t-lafuente} 
(1) A symmetric semi--Riemannian manifold is complete.

(2) In a locally symmetric manifold there is a full dependence
among the three types of causal completeness, i.e., completeness
in a causal sense implies completeness in the three senses.

(3) Any compact simply connected locally symmetric manifold is
complete.
\end{theorem}
\begin{proof} (1) If
$\gamma:[0,b)\rightarrow M$ were a geodesic inextendible
 to $b<+\infty$,  a contradiction is obtained by using the
global symmetry at  $p=\gamma(2b/3)$.

(2) The following extension to the semi--Riemannian case, of
classical Cartan's result\index{Theorem!Cartan} is well--known (see \cite[Theorem 8.17]{O}):
{\em if $M, \tilde M$ are complete connected locally symmetric
manifolds, $\tilde M$ is simply connected and $L: T_{\tilde
p}\tilde M\rightarrow T_pM$ is a linear isometry which preserves
the curvature, then there exists a unique semi--Riemannian covering
map $\phi: \tilde M\rightarrow M$ such that $d\phi_{\tilde p}
=L$.} In the proof of this result, one construct the covering map
by chaining local isometries along broken geodesics. One can check
that only broken geodesics of a causal type are necessary
(see \cite{La}) and, thus, the result follows if ``completeness'' is
weakened in ``completeness in one causal sense''. Thus, if $M$ is
complete in a causal sense, its universal covering $\tilde M$ will
be symmetric (the isometry $\phi: \tilde M\rightarrow \tilde M$
obtained for $L=-$Id will be a global symmetry) and, by (1),
complete.

(3) Any simply connected symmetric space $M$ can be affinely
inmersed in an affine symmetric space of the same dimension and,
moreover, when $M$ is compact, the affine inmersion becomes a
diffeomorphism (see \cite{FA}). Then, the result is a consequence of the
point (1).
\end{proof}

\paragraph{Alternative to Hopf--Rinow: a criterion
on completeness.} In the case $(M,g)$ is Riemannian,
classical Hopf--Rinow Theorem\index{Theorem!Hopf--Rinow} ensures the
equivalence between: {\sl (a)}
geodesic completeness, {\sl (b)} completeness of the associated distance
$d_g$, and {\sl (c)} the property that the closed and $d_g$--bounded
subsets of $M$ are compact. Obviously, in an indefinite manifold
{\sl (b)} does not makes sense, but one could still wonder if the
compactness of $M$ would imply completeness. It is well--known that
the answer is negative (see below). However, we reason now
where a proof (which is not based directly on the properties of
$d_g$) would fail; indeed, this will suggest some alternatives for
the indefinite case.

Assume that $M$ is compact and $g$ semi--Riemannian, take a
geodesic $\gamma: [0,b) \rightarrow M$ with $b<+\infty$, and try to
extend it beyond $b$. Recall: {\sl (i)} given any sequence $s_n \nearrow b$,
$(\gamma(s_n))_n$ converges to some $p\in M$,  up to a
subsequence, {\sl (ii)} $\gamma'$ can be seen as an integral curve in
the tangent manifold $TM$ of the geodesic vector field $G$ on
$TM$, {\sl (iii)} in the Riemannian case, the constancy of $c\equiv
|\gamma'|$ also implies that, up to a subsequence,
$(\gamma'(s_n))_n$ converges to some $p\in M$ (the bundle of
spheres of radius $c$ is compact). Now, recall that, by a
well--known result  (for example, \cite[Lemma 1.56]{O}), if an
integral curve $\rho:[0,b)\rightarrow M'$ of a vector field $X\in
\mathfrak{X}(M')$ is so that $(\rho (s_n))_n$ is convergent in
$M'$, then $\rho$ is extendible beyond $b$. Thus, the
extendibility in the Riemannian case follows putting $M'=TM, X=G,
\rho=\gamma'$.

Recall that {\sl (iii)} is the crucial step where the positive character of the metric (as well as the compactness of $M$) is used. 
Thus, previous proof yields directly the following
criterion.

\begin{proposition} \label{pcomplcriterion}
Let $(M,g)$ be a semi--Riemannian manifold and $\gamma:[0,b)
\rightarrow M$, $b<+\infty$, a geodesic. The following 
statements are equivalent:
\begin{itemize}
\item[{\sl (i)}] $\gamma$ is extendible beyond $b$;
\item[{\sl (ii)}] for some (and then for any) complete Riemannian metric $g_R$
on $M$, $|\gamma'|_R = g_R(\gamma',\gamma')^{1/2}$ is bounded;
\item[{\sl (iii)}] there exists a sequence $s_n \nearrow b$ such that
$(\gamma'(s_n))_n$  converges in $TM$.
\end{itemize}
\end{proposition}

\begin{remark} \label{rcomplcriterion}
 As a consequence, if $\gamma$ is an incomplete geodesic then
$\gamma'$ cannot be contained in a compact subset of $TM$.
Nevertheless,  even if $M$ is compact, a complete geodesic may
have $\gamma'$ not contained in a compact subset (see Remark
\ref{rwarpedindef}(2)).
\end{remark}

\paragraph{Limits of incomplete geodesics.}
 Assume in Proposition \ref{pcomplcriterion} that $\gamma$ is
incomplete and $M$ is compact and, thus, $|\gamma'(s_n)|_R
\rightarrow +\infty$. The set of all the (oriented) directions of
$TM$ can be regarded as the $g_R$--unit sphere bundle $S_RM\subset
TM$\index{$S_RM$}. If $M$ is compact, then so is $S_RM$ and, therefore, the
sequence of directions $v_n= \gamma'(s_n)/|\gamma'(s_n)|_R$
converges to some $v\in S_RM$, up to a subsequence. Necessarily,
$v$ is lightlike:
\[
g(v,v) = \lim_{n\rightarrow +\infty} g(v_n,v_n) =
\lim_{n\rightarrow +\infty} c/|\gamma'(s_n)|_R =0,
\]
where $c=g(\gamma',\gamma')\in \R$. Thus, summing up:
\begin{proposition} \label{pincomp-nulldirect}
Let $(M,g)$ be an incomplete compact semi--Riemannian manifold and
${\cal I} \subset S_RM$\index{${\cal I}$}
 the set of incomplete directions\index{set of incomplete directions}
(i.e., any geodesic with initial velocity in one such direction is
incomplete).
Then, the closure of ${\cal I}$ contains lightlike directions.
\end{proposition}

\begin{remark} \label{rincomp-nulldirect}  If ${\cal I}$ were closed then incompleteness for
compact $M$ would imply lightlike incompleteness. In the
non--compact case, it is trivial to show that, in general, ${\cal
I}$ is neither closed nor open (removing some points of $\R^2$
suffices, see \cite{RSjmp1}). In the compact one,  the question
 has been somewhat controversial, because the  closedness of ${\cal I}$
had been implicitly assumed in \cite{Yu}. Nevertheless, explicit
counterexamples which show  that ${\cal I}$ is not necessarily
closed (nor open) were constructed in \cite{RSjmp1} (see Remark
\ref{rtorosincomp1y2}(2) below). Remarkably, in these
counterexamples, the limit of some incomplete timelike and
spacelike directions is a complete lightlike direction, {\em but
there are other incomplete lightlike directions}. This lead
Romero and S\'anchez \cite{RSjmp1} to pose {\em the independence
of incompleteness in the compact case as an open question}.
\end{remark}

\subsection{Completeness under conformal symmetries}\label{sub2}
Here, we prove how our criterion of completeness (Proposition
\ref{pcomplcriterion}) can be applied to manifolds with timelike
conformal symmetries (Theorem \ref{tgeod1}), as well as conformally
homogeneous manifolds (Theorem \ref{t-marsden}), obtaining then a
generalization  of a result by Marsden (Remark \ref{r-marsden}).
Recall that a vector field $K\in \mathfrak{X}(M)$ is called
{\sl conformal Killing} \index{vector field!conformal Killing} if the
Lie derivative ${\cal L}$ \index{${\cal L}$} satisfies
$$
{\cal L}_K g= 2\sigma g
$$
(the local flows of $K$ are conformal maps) for some $\sigma \in
C^\infty(M)$. If $\gamma$ is a geodesic:
\begin{equation} \label{egeod-prod-Kill}
\frac{d}{ds}g(\gamma',K)= c \sigma \circ \gamma, \quad \hbox{for }
\; c=g(\gamma',\gamma').
\end{equation}
In the case $\sigma\equiv 0$, $K$ is {\em Killing} \index{vector field!Killing}
and $g(\gamma',K)$ is a constant.

\paragraph{General technique.} As a first application of Proposition
\ref{pcomplcriterion} for the Lorentzian case (and using the
notation there in what follows) we have the following result by
Romero and S\'anchez \cite{RSprocAMS}:

\begin{proposition} \label{p-pAMS}
Let $(M,g)$ be a compact Lorentzian manifold. If it admits a
timelike conformal Killing vector field then $g$ is complete.
\end{proposition}

\begin{proof}
Because of Lorentzian signature, the orthogonal bundle $K^\perp$
is Riemannian and, as $g(\gamma',\gamma')$ is constant,  it is
enough to check that the projection of $\gamma'$ in ${\rm span}(K)$ lies
in a compact subset. But this follows directly, because
(\ref{egeod-prod-Kill}) implies that $g(\gamma',K)$ is bounded on
$[0,b)$, if $ b<+\infty$.
\end{proof}

\begin{remark} \label{r-pAMS} The result does not hold if $K$ is allowed to be
causal (in particular, non--vanishing) at some points (see Remark
\ref{rtorosincomp1y2}).
\end{remark}

Obviously, the previous result can be extended to the case of index
$s_0$ if there are $s_0$ pointwise--independent Killing vector
fields. Even more, the compactness assumption can be dropped if
some additional conditions are imposed, yielding the following
general result (see \cite{RSgeoded}):

\begin{theorem} \label{tgeod1}
 A semi--Riemannian manifold $(M, g)$ of index $s_0$
is complete if there exist $s_0$ timelike conformal--Killing vector
fields $K_1, \dots K_{s_0}$ satisfying:
\begin{itemize}
\item[{\sl (i)}] the Gram matrix $\{g(K_i, K_j)\}$ has inverse ${g^{ij}}$, and
$\sum_{i,j} (g^{ij})^2$ is bounded;
\item[{\sl (ii)}]
functions $\sigma_i$ satisfying (\ref{egeod-prod-Kill}) for
the corresponding $K_i$, are bounded;
\item[{\sl (iii)}]
the associated Riemannian metric $g_R$ is complete, where
$g_R$ is given by
\[
g_R(X, Y) = g(X, Y),\quad g_R(A,B) = -g(A, B),\quad
g_R(X, A) = 0,
\]
for any $A, B \in$ ${\rm span}\{K_1,\dots,K_{s_0}\}$
and $X, Y$ belonging to its $g$--orthogonal complement.
\end{itemize}
\end{theorem}

\begin{remark} \label{r-ejemp-fibrado}
An example of metric where this result is applicable, is provided
by warped fiber bundles, such as the following Kaluza-Klein type. Let
$P(B,G)$, $\pi_B: P\rightarrow B$, a principal bundle on the
complete Riemannian manifold $(B,g_B)$ with structural group $G$,
and let $\omega: P \rightarrow \mathcal{G}$ be a connection 1--form
on the Lie algebra $\mathcal{G}$ of $G$. For any positive function
$f: B\rightarrow \R$ whose infimum satisfies $\inf f(B)>0$, and  any
bi--invariant metric $g_G$ on $\mathcal{G}$ (if it is semi--simple,
its Killing form suffices), define $g= \pi_B^* g_B + f^2 \omega^*
g_G$. The fundamental vectors of the bundle yield enough Killing
vector fields to ensure completeness.
\end{remark}

\paragraph{Homogeneous manifolds.} Previous technique also works for (conformally)
homogeneous manifolds.\index{manifold!homogeneous}
\index{manifold!conformally homogeneous} Recall that in a {\em
homogeneous manifold} $M$ any point $p \in M$ can be mapped to
another one $q\in M$ by means of an isometry. Remarkably then, any
tangent vector $v\in T_pM$ can be extended to a Killing vector $V$
(for a conformal metric, $V$ will be Killing conformal).

\begin{theorem}\label{t-marsden} A compact semi--Riemannian manifold which is globally conformal to
a homogeneous semi--Riemannian manifold is complete.
\end{theorem}
\begin{proof}
In order to apply Proposition \ref{pcomplcriterion}, take
$\displaystyle p=\lim_{n\to +\infty}\gamma(t_n)$, extend a basis of $T_pM$ to a set of
conformal Killing vector fields (which will be a pointwise basis
in a neighborhood of $p$) and use (\ref{egeod-prod-Kill}).
\end{proof}

\begin{remark}\label{r-marsden} (1) This is an extension of a result by Marsden
\cite{Marsden}, who considered the homogeneous case. In this case
$TM$ can be divided in compact subsets which are invariant by the
geodesic flow. Nevertheless, this property may not hold in the
conformal case (see Remarks \ref{rcomplcriterion} and
\ref{rwarpedindef}(2)).\\
(2) It is well--known that a homogeneous Riemannian manifold is
complete (use that at one point $p\in M$ the closed ball of some
radius  $r>0$ is compact; by homogeneity, this holds with the
same $r$ for any $q\in M$, and this property yields completeness
easily). This does not hold in the indefinite case (see below)
and, in fact, there are even incomplete compact locally
homogeneous manifolds (see \cite{GuLa}).
\end{remark}

\subsection{Heuristic construction of incomplete examples}
\label{sub3}
Now, our purpose is to construct examples of compact incomplete
semi--Riemannian manifolds. Even though  one such example can be
exhibited directly, we will construct it in a heuristic way, and
will find typical related items along the construction: incomplete
homogeneous manifolds, incomplete closed lightlike geodesics,
Misner's cylinder...

\paragraph{Incomplete homogeneous manifolds.} Let us start
with Lorentz--Minkowski spacetime\index{spacetime!Lorentz--Minkowski}
$\LL^2=(\R^2,g_0)$\index{$\LL^2$},
with its metric expressed in lightlike coordinates, and a
semiplane $M_0$: \begin{equation} \label{e-m0}
 g_0 = du \otimes dv + dv \otimes
du, \quad M_0=\{(u,v)\in \LL^2: u>0\}.
\end{equation}
Obviously, for any $\mu\neq 0$, the map
$$ \phi_{\mu}(u,v) = (\mu u, v/\mu) $$ is an isometry of $\LL^2$
inducible in $M_0$. As so they are the translations in the
direction of the $v$ axis too, {\em $(M_0,g_0)$ is an  incomplete
homogeneus semi--Riemannian manifold}.


\paragraph{Quotients by isometries and Misner's
cylinder.}\label{subMiner}
Now, consider the isometry group\index{group!isometry}
$G$ of both $\LL^2$ and $M_0$, generated
by $\phi\equiv \phi_2$, i.e., $ \phi(u,v) = (2u, v/2)$.
Consider first the action of $G$ on $\LL^2$. As the origen is a
fixed point, the action of $G$ is not free. Nevertheless, it
yields a free and discontinuous action\footnote{An action $\cdot$
of a group $G$ on a topological space $X$ is {\em free}
\index{action! free} if the equality $g\cdot x = x$ for some $x\in
X, g\in G$ implies that $ g$ is the the unit element $e$ of $G$.
It is {\em discontinuous} \index{action! discontinuous} if, for
any sequence $(g_m )_m$ of distinct elements of $G$ and $x\in
X$, the sequence $(g_m \cdot x)_m$ is not convergent. The
requirements for a {\em properly discontinuous} action
\index{action! properly discontinuous} are: (i) if $x, x'\in X$ do
not lie in an orbit by $G$ ($x'\neq g\cdot x$ for all $g\in G$)
there are neighborhoods $U\ni x, U'\ni x'$ such that $gU\cap
U'=\emptyset$ for all $g\in G$, (ii) the isotropy group $G_x$ at
each $x\in X$ ($G_x=\{ g\in G: g\cdot x = x \}$) is finite, and
(iii) for all $x \in X$, there exists an open neighborhood $U \ni
x $ such that: $gU\cap U =\emptyset$ for all $g \in G\backslash
G_x$.} on $\R^2_* = \R^2\backslash \{(0,0)\}$. \index{$\R^2_*$} In
the Riemannian case, these two conditions for an isometry action
are sufficient to obtain a (Hausdorff) quotient manifold, as the
action will be properly discontinuous (see \cite[I, Chapter 1]{KN},
especially Proposition 4.4). Nevertheless, this is not enough in the
indefinite case and, in fact, the quotient $\R^2_* /G$ is
non--Hausdorff.

\begin{remark}
Let ${\cal M}_n$ be the group of rigid motions of $\R^n$\index{group!rigid motions}
(semidirect product of translations and rotations) and ${\cal A}_n$ the one of
 affine maps\index{group!affine maps} (idem with translations and linear automorphisms). As
 a consequence of previous discussion, if a subgroup $G <{\cal
 M}_n$ acts freely and discontinuously on $\R^n$ the quotient is a
 (Hausdorff) manifold. Nevertheless,  the same
 assertion for $G <{\cal A}_n$ was pointed out as open
 by Charlap in \cite[pp. 4]{Char}.
\end{remark}
The action of $G$ on $M_0$ in (\ref{e-m0}) is properly
discontinuous and, in fact, the quotient  $C_0=M_0/G$ is
topologically a cylinder (see $C_0$ as $\{(u,v) \in \LL^2: 1\leq
u\leq 2\}$ with each $(1,v)$ identified to $(2,v/2)$), endowed
with a Lorentzian metric. This Lorentzian manifold (a simplified
example for the remarkable geometric properties of Taub--NUT
spacetime \cite{Misner}) will be called {\em Misner's
cylinder}\index{Misner's cylinder}. Obviously, $C_0$ is incomplete
in the three causal senses, as so is $M_0$. More strikingly, $M_0$
contains a closed incomplete lightlike geodesic, concretely the
projection $\gamma_0$ of the $M_0$--geodesic $\tilde
\gamma_0(s)=(-s,0)$ for all $s<0$.


\paragraph{Closed incomplete geodesics.} The reason for the
incompleteness of $\gamma_0$ can be explained as follows. Choose
$s_1<0$, put $s_2=s_1/2$, $T=s_2-s_1>0$ and recall that, from
$s_1$ to $s_2$ the geodesic gives a round
($\gamma_0(s_1)=\gamma_0(s_2))$, but the velocity after this round
 satisfies $\gamma'(s_1)=2\gamma(s_0)$. Thus, $\gamma_0$ will spend a time $T/2$ in
 giving a second round and, in general, $T/2^k$ for $k$ new
 rounds. This  yields directly  its incompleteness. As a
 straightforward generalization we have the following definition.

\begin{definition}
Let $(M,g)$ be a semi--Riemannian manifold, $\gamma: I=(a,b)
\rightarrow M$ a (inextensible) non--constant geodesic. Then,
$\gamma$ is {\em closed} \index{geodesic!closed} if there exists
$\lambda
>0$ and $s_1, s_2 \in I$, $s_1<s_2$, such that $\gamma(a)=\gamma(b)$
and:
\begin{equation} \label{dclosedg}
\gamma'(s_1)= \lambda \gamma'(s_2).
\end{equation}
In particular, $\gamma$ is called {\em periodic} if $\lambda = 1$. \index{geodesic!periodic}
\end{definition}

If equality (\ref{dclosedg})
holds for some non--constant geodesic $\gamma$ then:

{\sl (i)} Necessarily, $\lambda >0$. In fact, $\lambda \neq -1$
 because, in this case, the standard uniqueness of geodesics
 yields some $\epsilon>0$ such that
 $\gamma(s_1+s) = \gamma(s_2-s)$ for all $s\in [0,\epsilon]$.
 Even more, the points which satisfy this equality form an open and closed
 subset of $[s_1,s_2]$, obtaining the contradiction  $ \gamma'(s_c) =- \gamma'(s_c)$ at
 $s_c=(s_2-s_1)/2$. A similar contradiction arises if $\lambda<0$.

{\sl (ii)} If $\lambda \neq 1$ then necessarily $\gamma$ is lightlike
 (otherwise, the constancy of $g(\gamma',\gamma')$ yields a
 contradiction) and, by the reasoning above, incomplete.
 Thus, the notions {\em closed} and {\em periodic}
 are interchangeable for spacelike and timelike geodesics.

{\sl (iii)} If $\gamma$ is periodic then there exists a minimum $T>0$
(the period) \index{period} such that $\gamma(t)= \gamma(t+T)$
for all $t$. In particular, $\gamma$ is complete.

Summing up:

\begin{proposition}
Let $\gamma: I \rightarrow M$ be a closed geodesic in $(M,g)$.

(1) If $\gamma$ is not lightlike then it is periodic.

(2) $\gamma$ is complete if and only if it is periodic. In fact,
if (\ref{dclosedg}) holds with $\lambda >1$ (resp. $\lambda<1$)
then $I=(-\infty,b)$, $b<+\infty$ (resp. $I=(a,+\infty)$, $a>-\infty$).
\end{proposition}

\paragraph{Construction of incomplete compact manifolds.}
It is not difficult to realize that the metric of Misner's
cylinder $(C_0,g_0)$ in Subsection \ref{subMiner} can be modified
outside a strip $|v|\leq R$ in order to obtain a new Lorentzian
metric which can be induced in a torus by identifiying the points
with $v= 2R$. A way to make this explicitly is the following (see
\cite{RSjmp1,Sgrenoble}).

Consider the global change of variables $\Psi: \R^2 \rightarrow
M_0= \R^+\times \R$, $u(x,y)=e^y$, $v(x,y)=xe^{-y}$ and the pull--back
\begin{equation} \label{e-met-toro-incomp}
g_\tau= \Phi^*g_0 = dx \otimes dy + dy\otimes dx + \tau(x) dy^2
\end{equation}
where $\tau(x)=-2x$. The generator $\phi$ of Misner's group\index{group!Misner's} $G$
corresponds to the translation $(x,y)\mapsto (x,y+\log 2)$,
and the incomplete geodesic $\tilde \gamma$ of $M_0$ to a
reparametrization of the axis $x\equiv 0$. So, if we choose as
$\tau(x)$ a periodic function which behaves as $-2x$ around $x=0$,
an incomplete metric on $\R^2$ is obtained. This can be induced in
a quotient  torus obtained by means of two translations in the
directions of the two axis. Even more, by an explicit computation
of the geodesics of $g_\tau$ (which can be easily integrated as
$g_\tau(\gamma',\gamma')$ and $g_\tau(\gamma',\partial_y)$ are
constants) it follows:

\begin{proposition} \label{p-tau-y-toro}
If $\tau(0)=0$ but $\tau \not \equiv 0$, then metric $g_\tau$ in
(\ref{e-met-toro-incomp}) contains timelike, spacelike and lightlike
incomplete geodesics (which asymptote  the axis $x\equiv 0$).

Thus, if $\tau$ is additionally 1--periodic, the quotient torus
$T^2=\R^2/\Z^2$, with the induced metric from $g_\tau$, is
incomplete in the three causal senses.
\end{proposition}

\begin{remark} \label{rtorosincomp1y2} Under the stated hypotheses
for $(T^2,g_\tau)$:

(1) If  we choose $\tau\leq 0$ then the vector field induced in
$T^2$ from $K=\partial_y$ becomes Killing and causal, showing that
Proposition \ref{p-pAMS} cannot be extended to this case.

(2) If $ \tau'(0) =0$ then the axis $x\equiv 0$ can be
reparametrized as a complete geodesic, showing that the incomplete
directions ${\cal I}$ (see Proposition \ref{pincomp-nulldirect}) are not
closed.

(3) On the other hand, if one chooses either $\tau>0$ or $\tau<0$,
the torus will be complete (this can be either computed directly or
deduced from Proposition \ref{p-pAMS}, see also Theorem \ref{t-torus-conf-flat} below).
\end{remark}

\begin{remark} Starting by one such incomplete torus, it is not difficult to
find another interesting incomplete semi--Riemannian manifolds. For
example (see  \cite{GuLa}), one can construct an incomplete solid
torus and, taking into account that the 3--sphere can be obtained
by gluing two solid tori, a {\em  simply connected compact
incomplete Lorentzian manifold can be constructed}.
\end{remark}

\subsection{Surfaces}\label{sub4}

We have already seen quite a few bidimensional examples, which
become interesting for both, simplicity and  to test possible
results (in the spirit of \cite{Sm}). Here, we collect the properties  of any
such a Lorentzian surface $(S,g)$, focusing in the cases $S=\R^2,
S=T^2$ (as the Euler characteristic of $S$ would be 0, see
\cite[Theorem 2.4]{MS}, the former is the unique simply connected
case and the latter the unique compact orientable case). Our main
goals are to study: {\sl (a)} for $S=\R^2$, the stability of
completeness \index{stability of completeness} in the $C^r$
topologies (see Theorem \ref{t-r2-estab}), {\sl (b)}~con\-formally
flat tori, including stability of completeness (see Theorem
\ref{t-torus-conf-flat}), {\sl (c)} structure of tori with a
Killing field, including the dependence of causal completitudes
(see Theorem \ref{t-trans}), and {\sl (d)} general tori, including
the problem of dependence of completitudes, studied by means of
its two lightlike foliations\index{lighlike foliations}. Recall
that any Lorentzian surface admits lightlike coordinates and,
thus, it is locally conformally flat (see, for example,
\cite{We}). So, {\em conformally flat} \index{manifold!conformally
flat} means {\em globally}
 conformal to a flat manifold in the remainder.

\paragraph{Case $S=\R^2$.} A  remarkable difference with the
Riemannian case appears for the  conformal class of $(\R^2, g)$.
If $g$ is Riemannian, Theorem of Uniformization\index{Theorem!Uniformization}
implies that $g$ is conformally equivalent either to the unit disc or to $\R^2$
with their natural metrics. Nevertheless, if $g$ is Lorentzian
there exist infinitely many conformal classes (see Weinstein
\cite[pp. 45]{We} for a detailed study or \cite[footnote 10]{MS} in these
proceedings).

Any  metric $g$ on $\R^2$ is stably causal (see \cite[Theorem 3.55]{MS})
and, as $-g$ is also Lorentzian, the roles of spacelike and
timelike geodesics are interchangeable. In particular, no
inextendible geodesic of any causal type is imprisoned
in a compact subset in the forward or backward direction. Thus,
the results of {\em stability of geodesic incompleteness}
in \cite[Theorem 7.30]{BEE}, or Theorem \ref{t-estC1} below (we also refer to
\cite{BEE} or \cite{MS} for definitions) apply, yielding:

\begin{theorem} \label{t-r2-estab}
If the Lorentzian surface $(\R^2,g)$ is incomplete in some causal
sense, then there is a $C^1$--neighborhood $U(g)$ of $g$ such that
each metric $g_1$ in $U(g)$ is incomplete in that causal sense.
\end{theorem}

\begin{remark} The $C^1$ stability of the
completeness of surfaces as $\LL^2$
is also known (see Theorem \ref{t-est-compl-C1} below).
\end{remark}

\paragraph{Case $S=T^2$\index{torus}\index{$T^2$}.} Recall first the following result.

\begin{lemma} \label{ltecn-torus} (1) Any flat Lorentzian torus is complete and, in fact, the
 quotient of $\LL^2$ by the action generated by two independent
translations.

 (2) Infinitely many classes of
non--conformally related flat Lorentzian torus exist.

(3) The tori which admit a periodic lightlike geodesic $\gamma_0$
are dense in the set of all the flat Lorentzian tori, for any
$C^r$--topology.
\end{lemma}

\begin{proof}
(1) Completeness (and, therefore, the result) can be seen as
a particular case of  Theorem \ref{t-markus} below, but direct
proofs are possible (see \cite{FGH,FF}).

(2) Choose as generators of the translations $v=(1,0)$ and
$w_k=(0,k)$ for $ k$ a (prime) natural number. Each  quotient torus
$T_k$ admits a periodic lightlike pregeodesic (which is a
conformal invariant) in a different free homotopy class.

(3) This property follows from a commensurability argument between
the generators of the translations, using the density of the
rationals in $\R$.
\end{proof}

Now, we can see that, even though (conformally) flat Lorentzian
tori are complete, this is an unstable property (see \cite{RSjmp2}):

\begin{theorem} \label{t-torus-conf-flat}
(1) A Lorentzian torus $(T^2,g)$ is conformally flat if and only
if it admits a timelike (or spacelike) conformal vector field $K$.
Thus, it is complete.

(2) Conformally flat Lorentzian metrics on $T^2$ lie in the
closure of the set of lightlike incomplete Lorentzian metrics.
Moreover, the flat ones lie in the closure of timelike, lightlike
and spacelike incomplete Lorentzian metrics.
\end{theorem}

\begin{proof}
(1) ($\Rightarrow$) From Lemma \ref{ltecn-torus}(1), any flat
torus admits a parallel vector field of any causal character,
which will be conformal Killing for any conformal metric.

($\Leftarrow$) Easily, $K$ is Killing for the conformal metric
$g^*=g/|K|$ and, in  dimension 2, it is also parallel and $g^\star$ flat
(the last assertion follows from Proposition \ref{p-pAMS}).

(2) By Lemma \ref{ltecn-torus}(3), and the conformal invariance of
the lightlike incompleteness (see Theorem \ref{t-comp-conf-luminosa}), it
is enough to prove the property for the flat tori which admit a
periodic lightlike geodesic $\gamma_0$. Recall that this geodesic
can be lifted to the universal covering $\R^2$ as an affine
parametrization of $x\equiv 0$. Now, choose $\tau$ as in
Proposition \ref{p-tau-y-toro} arbitrarily close to 0 in the $C^r$
topology.
\end{proof}

The case of non--conformally flat Lorentzian tori admitting a
Killing vector field $K$ can be also characterized in a precise
way (see \cite{Sa-trans}); surprisingly, all of them are incomplete. In
order to describe not only this result but also the structure of
these tori, consider first the following generalization of the
metric (\ref{e-met-toro-incomp}):

\begin{equation}
\label{3.1} g[(x,y)]=E(x)dx^2+F(x)\left(dx\otimes dy+dy\otimes
dx\right) -G(x)dy^2,
\end{equation}
where $E, F, G\in C^2(\R )$ satisfy the following conditions: (i)
$EG+F^2>0$, that is, $g$ is Lorentzian, and (ii) $E, F$ and $G$
are periodic with period 1, so, the metric is naturally
inducible on a torus $T^2=\R^2/\Z^2$.  As the vector field
$K=\partial/\partial y$ is Killing, if $|G|$ were greater than 0
at every  point, then $K$ would be either timelike or spacelike,
and thus, $g$ is conformally flat. Moreover, if $G\equiv0$ then
$K$ is lightlike and $g$ is flat. Now, denote by ${\cal G}$ (resp.
${\cal G}^c$)  the set of metrics given by (\ref{3.1}) and
satisfying (i), (ii) and also: (iii) the sign of $G$ is not
constant (resp. (iii)$^c$ the function $G$ is not constant and
$|G|>0$). Recall that the geodesic equations for
$\gamma(t)=(x(t),y(t))$ can be explicitly integrated (see
the comments above Proposition \ref{p-tau-y-toro}).
This, combined with some technicalities  (for example, a
non--trivial Killing vector field cannot vanish on a torus),
 yields (see \cite{Sa-trans}):
\begin{theorem} \label{t-trans}
If a Lorentzian torus $(T^2,g)$ admits a  Killing vector field
$K\not\equiv 0$, then $K$ does not vanish at any point and:
\begin{itemize}
\item[{\sl (1)}] The metric $g$ is flat if and only if $g(K,K)$ is constant.
\item[{\sl (2)}] The metric $g$ is conformally  flat  if  and  only  if
$g(K,K)$  has a definite sign (strictly positive, strictly
negative or identically zero), and if and only if $g$ is
geodesically complete (in the three  causal senses).\\
Moreover, $g$ is conformally flat but non--flat if and only if it
is isometric to one of the ${\cal G}^c$--tori
constructed above, up to a covering.
\item[{\sl (3)}] The metric $g$ is non--conformally flat if and only if $g$ is
geodesically incomplete in one (and then in the three) causal
senses, and if and only  if  it is isometric to
one of the ${\cal G}$--tori constructed above, up to a covering.
\end{itemize}
\end{theorem}

\begin{remark} (1) For tori admitting a Killing vector field there
is equivalence between the three types of causal completeness.

(2) For a torus conformal to a previous one,  all the assertions
in Theorem \ref{t-trans} hold with obvious modifications, except
that, in principle,  one can ensure only the equivalence with
lightlike incompleteness in the case (3).

(3) One of these incomplete tori is the celebrated Clifton--Pohl
one\index{Clifton--Pohl torus}, usually defined as the quotient of
 ($\R^2_*$, $g=(u^2+v^2)^{-1}(du\otimes dv+dv\otimes du)$) by the
 isometry group generated by the homothety
$(u,v)\mapsto (2u,2v)$ for all $(u,v)\in\R^{2}_*$. Necessarily,  it
is a ${\cal G}$--torus ($K = u \partial/\partial u + v \partial/\partial v$
is Killing), but this can be checked directly (see \cite[Remark 5.2]{Sa-trans}).
\end{remark}

When no (conformal) Killing vector field on $T^2$ exists, the study of the
geodesics becomes more difficult. The completeness of geodesics in
simultaneously all the elements of a conformal class of time-orientable Lorentzian
tori has been studied from the next point of view by Carri\`ere
and Rozoy in \cite{CR}. Consider the two (conformally invariant)
foliations $\cal F$, $\cal H$ yielded by the lightlike geodesics
of $T^{2}$. Then one has (we refer to \cite{CR} and references
therein for the explanation of concepts relative to foliations):

\begin{theorem}
\label{teo-cr} {\rm (1)} For arbitrary $\cal H$:

{\rm (1a)} if $(T^{2},g)$ is lightlike complete then $\cal F$ and
$\cal H$ are $C^{0}$--linearizable;

{\rm (1b)} if $\cal F$ and $\cal H$ are $C^{1}$--linearizable then
$g$ is lightlike complete.\\
{\rm (2)} In the particular case that $\cal H$ is a foliation by
circles, then $\cal F$ is obtained from the suspension of a
diffeomorphism $\phi$ of the circle, and then:

{\rm (2a)} $g$ is lightlike complete if and only if $\phi$ is
$C^{0}$--conjugate to a rotation,

{\rm (2b)} if $\phi$ is $C^{1}$--conjugate to a rotation then $g$
is complete.
\end{theorem}

\begin{remark} \label{r-cr}  As a
consequence of Theorem~\ref{teo-cr} one has that, at least when
${\cal H}$ is a foliation by circles, ``generically'' any incomplete
Lorentzian tori must be lightlike incomplete (for almost all real
number,  $C^0$--conjugate to a rotation implies $C^1$).
Nevertheless, there are residual cases which  are compatible with
$g$ incomplete but lightlike complete; this suggests that such a
possibility will hold.
\end{remark}

\subsection{Influence of curvature}\label{sub5}

In semi--Riemannian Geometry, as in the Riemannian one,  the
curvature determines the metric (for example, in the sense of
Cartan's result, see the proof of Theorem \ref{t-lafuente}(2)
or \cite{Re}).
Nevertheless, very subtle questions about a global
property such as completeness, appear in the indefinite case. Here, we
show first that the curvature does not characterize completeness,
even in the compact case (Theorem \ref{t-putadacurv}). This stresses
Markus' conjecture as well as the completeness of compact
manifolds of constant curvature (Theorem \ref{t-markus}). The nice
behaviour of geodesics in the  complete case (i.e., spaceforms),
is described around Theorem \ref{t-calabi}. Finally, even though
relativistic {\em Singularity Theorems} become a substantial topic
in its own right (see \cite{Se}), we give some flavour of them in
comparison with Riemannian Geometry (Theorems \ref{thaw},
\ref{t-mishuevos}).

\paragraph{Complete and incomplete tori with the same
curvature.} A first strong sense of independence between curvature
and completeness is the following (see \cite{RSjmp1}):

\begin{theorem} \label{t-putadacurv}
There are two metrics $g_1$, $g_2$ on the same torus $T^2$, with
the same curvature at each point and such that the first one $g_1$
is incomplete and the second one $g_2$ is complete.
\end{theorem}

\begin{proof}
Consider the metric (\ref{e-met-toro-incomp}), which has Gauss
curvature  $\tau''/2$. Now, choose as $\tau$ for $g_1$ any
1--periodic function $\tau_1$ with $\tau_1(0)=0 \neq  \tau_1'(0)$
(it is incomplete by Proposition \ref{p-tau-y-toro}). For $g_2$ put
$\tau_2=\tau_1 + N$, where $N$ is chosen such that $\tau_2>0$, and
notice that it is complete by Remark \ref{rtorosincomp1y2}(3).
\end{proof}

\begin{remark}
This result can be generalized to obtain complete and incomplete
Lorentzian torus with prescribed curvature $k=k(x)$, satisfying
the condition of compatibility with Gauss--Bonnet
Theorem\index{Theorem!Gauss--Bonnet} (i.e., $\int_0^1k(x)dx=0$).
Applications to solutions of D'Alembert equation can be seen in
\cite{Sa-trans}.
\end{remark}

\paragraph{Markus conjecture.} In
spite of previous result, one can wonder if a strong restriction
on the curvature, as flatness, will imply completeness. Such a
question is a particular case of Markus conjecture\index{Markus conjecture} on affine
manifolds.

An {\em affine manifold} \index{manifold!affine} is a manifold $M$
locally modelled on open subsets of $\R^n$ such that the changes
of coordinates are elements of the group of affine transformations
${\rm Aff}(\R^n)$\index{${\rm Aff}(\R^n)$}. Fixing $p\in M$ one has the natural {\em
representation of holonomy}\index{representation of holonomy}
for the fundamental group $\pi_1(M)$,
$ h: \pi_1(M) \mapsto h(\pi_1(M)) =\Gamma \subset
\hbox{Aff}(\R^n)$, which is independent of $p$ up to conjugacy.
Let $L(\Gamma)\subset $Gl$(\R^n)$ be the linear part of $\Gamma$
obtained from the natural projection.
\begin{quote}
{\em Conjecture} (Markus): An unimodular (i.e., $L(\Gamma) \subset
SL(\R^n)$) compact affine manifold is complete.
\end{quote}
In the flat Lorentz case one has naturally an affine manifold with
$L(\Gamma)$ included (up to an orientable covering) in the special
Lorentz group SO$_1(n)$. Carri\`ere \cite{Ca} solved this case by:
(i) introducing an invariant for any subgroup $G\subset
$Gl$(\R^n)$, the discompactness disc$(G)$, which measures at what
extent $G$ fails to be compact (in particular, disc(SO$(n)$)$=1$),
(ii) proving Markus' conjecture for disc(SO$(n)$)$\leq 1$. Even
more, Klingler \cite{Kl} extended this result to include the
manifolds of constant curvature. Summing up:
\begin{theorem} \label{t-markus}
Any  compact Lorentzian manifold  $(M,g)$ of constant curvature is
complete.
\end{theorem}
\begin{remark}
(1) As far as we know, the question remains open for $g$
indefinite with higher index.

(2) Obviously, here the techniques are very different to those ones for
Proposition \ref{pcomplcriterion} and, in particular, they are not
conformal invariant. Thus, one can wonder if completeness will
hold for any conformal metric $g^* = \Omega g$. This will hold  if
$M$ is a $n$--torus or a nilmanifold, because in this case $(M,g)$
will admit a timelike Killing vector field and Proposition
\ref{p-pAMS} applies.

(3) Notice that:
\[
\hbox{constant curvature $\Rightarrow$ locally
symmetric $\Rightarrow$ locally homogeneus.}
\]
As we have seen, at
least in the compact Lorentzian case the first condition implies
completeness,  but the last one does not (see Remark \ref{r-marsden}).
As far as we know, also the intermediate compact locally symmetric case
remains open  (see also Theorem \ref{t-lafuente}).
\end{remark}

\paragraph{Spaceforms.} A semi--Riemannian $n$--dimensional manifold $M$ of
index $s$  is a {\em spaceform} \index{manifold!spaceform} if it is
complete with constant curvature. This curvature will be regarded
as normalized to $\epsilon= -1, 0, 1$. Simply connected spaceforms
are called {\em  model spaces} \index{manifold!model space} and are
characterized by $n, \epsilon, s$. Essentially, the model spaces
are either $\R^n_s$ or the pseudosphere \index{pseudosphere}
$\SSS^n_s$ \index{$\SSS^n_s$} (spacelike vectors of norm 1 in
$\R^{n+1}_{s}$); but for $n=2$ the pseudosphere is topologically
$\SSS^1\times \R$ and its universal covering must be taken. Any
spaceform can be constructed as the quotient of its model space by
a group $\Lambda$. As we have seen above, if the assumption of
completeness is dropped then the universal covering may be a
proper open subset of the model space, which makes the study more
difficult.

Geodesics on $\SSS^n_s$ can be constructed by intersecting
$\SSS^n_s$ with a plane in $\R^{n+1}_{s}$ which crosses the
origin. Such geodesics are well known (see \cite[Proposition 5.38]{O}) and,
in particular, {\it no indefinite pseudosphere $\SSS^n_s$, $0<s <
n$, is geodesically connected.} Nevertheless, even non--flat
spaceforms may be geodesically connected. We recall the following
result by Calabi and Markus \cite{CM} (which, in particular, solves
completely the geodesic connectedness of Lorentzian spaceforms of
positive curvature with $n\geq 3$) and refer also to \cite{O,Sa01}
for further information.

\begin{theorem} \label{t-calabi} (1) Two points $p, q \in \SSS^n_1$ ($n\geq 2$)
are connectable by a geodesic
if and only if $\langle p,q\rangle_1 >-1$,  where
$\langle\cdot,\cdot \rangle_1$ is the  inner product of $\LL^{n+1}
\equiv \R^{n+1}_1$.\\
(2) Any  spaceform $M = \SSS^n_{1} /\Lambda, M \neq \SSS^n_{1}$ is
starshaped from some point $p\in M$.\\
(3) A spaceform $M = \SSS^n_{1}/\Lambda$ is geodesically connected
if and only if it is not time--orientable.
\end{theorem}

The proof of (1) follows from a direct computation of the geodesics.
For the remainder, the essential idea is that (whenever $2s \leq
n$), the group  $\Lambda$ is finite. Then, up to conjugacy,
$\Lambda \subset O(1)\times O(n) \subset O_1(n+1)$, and the proof
follows by studying the barycenter of the orbits, which must lie
in the timelike axis of $\R_1^{n+1}$.

\paragraph{Singularity theorems.} In Riemannian Geometry, it is
well--known that negative (sectional) curvature implies divergence
of geodesics, while positive curvature, or even just positive Ricci
curvature, implies convergence and focalization. This result, in
combination with properties of the distance function, implies
bounds on the diameter (and, then, compactness) for complete
Riemannian metrics with Ricci curvature greater than some positive
constant (see Myers' Theorem). In Lorentzian manifolds the properties
of convergence and focalization hold analogously for timelike
geodesics (and, with some particularities, for lightlike ones)
under positive Ricci curvature on timelike vectors (or negative
sectional curvature on timelike planes)\footnote{The imposition of
inequalities only on timelike vectors (or on timelike planes)
becomes essential for the mathematical non--triviality of the
problem, as well as for physical interpretations. So, inequality
${\rm Ric}(v,v)\geq 0$ for any tangent vector $v$ imply  (in dimension
$\geq 3$) that the manifold is Einstein. But this inequality
imposed only for timelike $v$ is mathematically natural, and
admits the physical interpretation that {\em gravity, on average,
attracts} (timelike convergence condition).}. Nevertheless, the
inexistence of a true distance in a Lorentzian manifold and the
particular properties of their global structure, make so that the
natural conclusion cannot be the finiteness of the diameter but
the causal incompleteness of the manifold.

In fact, in the framework of General Relativity, singularity
theorems are incompleteness theorems for causal geodesics.
 Essentially, there are two
types which prove: (1)  the existence of  incomplete timelike
geodesics in a cosmological setting,  and (2)  the existence of an
incomplete lightlike geodesic in the context of gravitational
collapse and black holes. In general, they use some refined
properties of Causality which lie out of the scope of the present
paper. But we can give some ideas about them by studying a ``Big
Bang'' singularity theorem which has a clear correspondence with
the Riemannian ideas commented above.

Recall the following Hawking's Singularity
Theorem\index{Theorem!Hawking's singularity} (the necessary
concepts on Causality can be seen in these proceedings \cite{MS};
see \cite{HE}, \cite{O} or \cite{Se} for a detailed exposition of
this  singularity theorem and others):

\begin{theorem} \label{thaw} Let $(M,g)$ be a spacetime such that:
\begin{itemize}
\item[1.] it is globally hyperbolic;
\item[2.] some spacelike Cauchy hypersurface $S$ is strictly expanding $H
\geq C >0$ ($H$ is the mean curvature with respect to the future
direction);
\item[3.]  ${\rm Ric}(v,v)\geq 0$ for timelike $v$.
\end{itemize}
Then, any past--directed timelike geodesic  $\gamma$ is
incomplete.
\end{theorem}

\begin{proof} The last two hypotheses imply that  any past--directed geodesic
$\rho$ normal to $S$ contains a focal point if it has length
$L'\geq \frac{1}{C}$. Thus, once $S$ is crossed, no  $\gamma$ can
have a point $p$ at length $L>\frac{1}{C}$  (otherwise,  a
length--maximizing timelike geodesic from $p$ to $S$ with length
$L'\geq L$ would exist by global hyperbolicity, a contradiction).
\end{proof}

Now, as an exercise, the reader can prove the Riemannian result
below and stress the isomorphic roles of:
\begin{enumerate}
\item Global hyperbolicity $\longleftrightarrow$ Riemannian
completeness
\item ${\rm Ric}(v,v)\geq 0$ for timelike $v$ $\longleftrightarrow$
${\rm Ric}(v,v)\geq 0$ for all tangent  $v$
\item Timelike incompleteness (bounded lengths of timelike
geodesics starting at $S$) $\longleftrightarrow$ Finite Riemannian
distance to $S$.
\end{enumerate}

\begin{theorem} \label{t-mishuevos}
 Let $(M,g)$ be a Riemannian manifold such that:
\begin{itemize}
\item[1.] it is complete;
\item[2.] some closed (as a subset) hypersurface $S$ separates $M$ as
the disjoint union $M=M_-\cup S \cup M_+$,  and $S$ is strictly
expanding towards $M_+$: $H \leq -C <0$ (with appropiate sign
convention for $H$);
\item[3.] ${\rm Ric}(v,v)\geq 0$ for all $v$.
\end{itemize}
Then, ${\rm dist}(p,S)\leq 1/C$ for all $p\in M_-$.
\end{theorem}

\subsection{Warped products} \label{sub6} Here, we study the
behaviour of geodesics in a general warped product. For geodesic
completeness, the details of the fiber become irrelevant (Theorem
\ref{t-warpedcomplete}), and leads to the notion of {\em warped
completeness}. This can be characterized very accurately in the
case of definite base (Theorem \ref{t-warpedcomplete}), and is related
to some open questions in the case of indefinite base (Remark
\ref{r-dna-wac-coca}).

A {\em warped product} \index{manifold!warped product} $B\times_fF$
\index{$B\times_f F$} of the semi--Riemannian manifolds $(B, g_B)$
(base) and $(F,g_F)$ (fiber) with warping function $f:
B\rightarrow \R, f>0,$ is the product manifold $B\times F$ endowed
with the warped metric
\begin{equation} \label{e-mwarped}
g= \pi_B^* g_B + f^2 \pi_F^* g_F, \end{equation} where $\pi_B,
\pi_F $ are the natural projections of the product $B\times F$.
The geodesic equations and elements of curvature of a warped
product (expressed in the general semi--Riemannian setting) has
been systematically studied by O'Neill \cite[Chapter 7]{O}. It is
straightforward to check that if the base and fibers are complete
and Riemannian then the warped product is complete, too (use
Proposition \ref{pcomplcriterion}), But this does not hold in the 
indefinite case, as stressed by the
simple Beem and Buseman counterexample\index{counterexample!Beem
and Buseman} $\R\times_f\R, g=dx^2-e^x dy^2$. A careful study of
geodesic completeness was carried out in \cite{RSgeoded}, which
will be our main reference. Extensions to other type of multiply
warped manifolds can be found in \cite{Unal1,Unal2}.

\paragraph{Warped completeness.}
A curve  $\gamma =(\gamma_B,\gamma_F)$ in $B\times F$ is a
geodesic if and only if it satisfies:
\begin{equation} \label{egeodwarped}
\left\{\begin{array}{l}\displaystyle
\frac{D  \gamma'_B}{dt} = \frac{C}{(f\circ \gamma_B)^3}\ \nabla^Bf
\\
\displaystyle\frac{D  \gamma'_F}{dt} = -\frac{2}{f\circ \gamma_B}
\frac{d(f\circ \gamma_B)}{dt}\ \gamma'_F,
\end{array}\right.
\end{equation}
where $C= (f\circ \gamma_B)^4 g_F(\gamma'_F,\gamma'_F)$ is
necessarily a constant, $\nabla^B$\index{$\nabla^B$}
denotes the $g_B$--gradient, and
$D/dt$\index{$D/dt$} denotes the covariant derivative in the corresponding
manifold.

\begin{remark}  The equation for $\gamma_F$ implies that it is
a pregeodesic for $g_F$ (recall (\ref{e-reparam-pregeodlum})). As
the equation for $\gamma_B$ is independent of $\gamma_F$, up to
the constant $C$, one must study only an equation on $B$ and a
reparameterization.
\end{remark}
Now, taking into account the previous remark, the proof of Proposition
\ref{pcomplcriterion},  and some well--known facts (say, if a
geodesic is continuously extendible to a point then it is also
extendible as a geodesic (see \cite[Lemma 5.8]{O}) one easily has:

\begin{lemma} \label{lwarped} If the fiber is complete, then
for a geodesic $\gamma:[0,b)\rightarrow B\times F$, $b<+\infty$,
the following properties are equivalent:

 (1) $\gamma$  is extendible as
a geodesic beyond $b$;

(2) $\gamma_B$ is continuously extendible beyond $b$;

(3) $\gamma'_B$  lies in a compact subset of $TB$.

\noindent Moreover, if $g_B$ is  Riemannian,  previous conditions
are also equivalent to:

(4) $\gamma_B$  lies in a compact subset of $B$.
\end{lemma}

 From this result one can see that the role of the chosen fiber is
irrelevant for the completeness of $\gamma$, except for the fact
that, if it is Riemannian, the value of the constant $C$ in
\rf{egeodwarped} is always positive. More precisely (see
\cite{RSgeoded}):

\begin{theorem} \label{t-warpedcomplete} If the fiber of a warped product $B\times_fF$ is
incomplete then the warped product is incomplete (and in the three
causal senses, if the base is indefinite).

If $(F,g_F)$ is complete and indefinite,  the following assertions
are equivalent:
\begin{itemize}
\item[(i)] $B\times_fF$ is (resp. timelike, lightlike or spacelike)
complete;
\item[(ii)] for any other complete and indefinite fiber $(F',g_{F'})$ the
warped product $B\times_fF'$ is (resp. timelike, lightlike or
spacelike) complete.
\end{itemize}
\end{theorem}

\begin{definition}
A triple $(B,g_B,f)$ is called {\em   (resp. timelike, lightlike
or spacelike) warped complete} \index{warped complete triple} if for any
complete fiber $(F,g_F)$ the warped product $B\times F$ is (resp.
timelike, lightlike or spacelike) complete.
\end{definition}

More technically, any solution $\gamma_B$ of the differential
equation \rf{egeodwarped} for $(B,g_B,f)$ will be called a {\em
(warped) geodesic projection} \index{warped geodesic projection}
and it will be called timelike, lightlike or spacelike depending
on if the (necessarily constant) value of $D=
g_B(\gamma'_B,\gamma'_B) + C/(f\circ \gamma_B)^2$ is negative, 0
or positive. Theorem \ref{t-warpedcomplete} can be paraphrased by
saying that: {\em a warped product with a complete indefinite
fiber is (resp. timelike, lightlike or spacelike) complete if and
only if all its (resp. timelike, lightlike or spacelike) geodesic
projections are complete}.

\begin{remark}
 All the study will hold also if $B\times F$ is replaced by a
fiber bundle $E(B,F)$ with base $B$ and fiber $F$, endowed with a
metric as \rf{e-mwarped} in each trivializing neighborhood,  where
$\pi_B, \pi_F $ are now the natural projections of the fiber
bundle. An example would be the one in Remark
\ref{r-ejemp-fibrado}, when the connection $\omega$ is flat.
\end{remark}

\paragraph{Definite basis.} Now, consider that the metric $g_B$ is
Riemannian and $d$ is its distance. As the base of the warped
product is totally geodesic, if $g_B$ is incomplete then
$(B,g_B,f)$ is at least  spacelike warped incomplete.
Nevertheless, completeness in other causal senses may still hold.
Now, we will focus in the case $g_B$ complete, and will consider
an incomplete base in Subsection \ref{subGRW}.

In the case that $B$ is compact, Lemma \ref{lwarped}(4) yields
trivially the completeness. Otherwise, the behavior of $f$ at
infinity belongs crucial. Fix a point $x_0\in B$ and define the
continuous function $f_{inf}:[0,+\infty) \rightarrow \R$ so that
\[
f_{inf}(r) = \, \min\{f(x): d(x,x_0)=r\}.
\]
By using (\ref{egeodwarped}), the minimum increasing of the
parameter $s$ of $\gamma_B$ when $d(\gamma_B(s),x_0)=r$ can be
bounded from below and, if $f$ is radial from $x_0$ (i.e.,
$f_{inf}=f$), also from above. So, one obtains (see
\cite{RSgeoded}):

\begin{theorem}
\label{twarpedcompl} Let $(B,g_B)$ be Riemannian and complete. If
$B$ is compact, then the triple $(B,g_B,f)$ is warped complete.
Otherwise, $(B,g_B,f)$ is:
\begin{itemize}
\item[(a)]  warped complete if
\begin{equation}
\label{1-comp} \int_{0}^{+\infty}
\frac{f_{inf}}{\sqrt{1+f_{inf}^{2}}}\ dr = +\infty ;
\end{equation}
\item[(b)] timelike and lightlike complete  if
\begin{equation}
\label{2-comp} \int_{0}^{+\infty} f_{inf} dr = +\infty;
\end{equation}
\item[(c)] timelike complete if  either $f_{inf}$ satisfies {\rm
(\ref{2-comp})} or $f_{inf}$ is unbounded.
\end{itemize}
Moreover, if $f$ is radial, then the converses to (a), (b), (c)
hold, too.
\end{theorem}

\begin{remark}
 (1) The sufficient conditions in the items above satifies
 \rf{1-comp} $\Rightarrow$ \rf{2-comp} (and the latter implies the sufficient condition
 (c)) but the converses do not
hold. Obviously, if $\inf f(B) > 0$ then the  triple is warped
complete.

(2) If we consider a ``twisted product'' \index{twisted product},
i.e., instead of a warping  function $f$ in the definition of $g$,
a twisting one $h: B\times F \rightarrow \R$, $h>0$,  then none
of previous results hold. In fact, take $B=F=S^1$ (the standard
unit circumpherence in $\C$) and put
$h(e^{i\theta_1},e^{i\theta_2})= e^{\sin(\theta_1-\theta_2)}$. A
simple computation shows that the twisted product is incomplete
(in fact, isometric to one of the ${\cal G}$ torus  in Theorem
\ref{t-trans}).
\end{remark}

\paragraph{Indefinite basis.} If $(B,g_B)$ is indefinite then
no such accurate results hold. In fact, there are even some open
questions which are related with other studied problems. As a
 first result obtained by combining the ideas for Theorems \ref{tgeod1},
 \ref{twarpedcompl}, we have:

 \begin{proposition} \label{pwarpedindef}
Let $(B,g_B)$ be an indefinite manifold with $s$ conformal Killing
vectors fields $K_1, \dots , K_{s}$ satisfying the
hypotheses in Theorem \ref{tgeod1}. For any smooth function $f$ on
$B$ such that $\inf f(B)>0$ and all $K_i(f)$ are bounded, the triple
$(B,g_B,f)$ is warped complete.
 \end{proposition}

\begin{remark} \label{rwarpedindef}
(1) The boundedness of $K_i(f)$ cannot be dropped, otherwise there
are explicit counterexamples in \cite[Counterexample
3.17]{RSgeoded}. This suggests that no accurate and general result
on completeness such as Theorem \ref{twarpedcompl},  can be
obtained for indefinite bases.

(2) Choose as base the natural Lorentzian torus $\LL^2/\Z^2$ and
$f(x,y) = e^{\sin (x-y)}$. The corresponding triple is warped complete
as an application of Proposition \ref{pwarpedindef}. Nevertheless,
any non--constant geodesic projection $\gamma_B$ with $x(t)\equiv
y(t)$ is so that $\gamma'_B$ is not contained in a compact
subset of $TB$. Thus, choosing a compact fiber, an example of {\em
complete compact semi--Riemannian manifold with geodesics not
contained in any compact subset of the tangent bundle} is obtained
(compare with Remarks \ref{rcomplcriterion} and
\ref{r-marsden}(1)).
\end{remark}
Nevertheles, in spite of the first remark, the case of lightlike
geodesics depends only on the bounds of $f$ and is related to
the completeness of conformal metrics. Concretely, one can show:

\begin{theorem} \label{tdna-wac-coca} Let $(B, g_B)$ be a semi--Riemannian
manifold, $f\in C^\infty(M), f>0$.

(1) $\sup f(B) < +\infty$ and $(B,g_B,f)$ lightlike warped complete,
imply $g_B/f^2$ complete.

(2) $\inf f(B) >0$ and $(B,g_B,f)$  not lightlike warped complete,
imply $g_B/f^2$  not complete.
\end{theorem}
\begin{remark}\label{r-dna-wac-coca} Consider the following open questions:
\begin{itemize}
\item[$(Q1)$] A lightlike complete  compact indefinite manifold is
complete.
\item[$(Q2)$] Any triple $(B,g_B,f)$ with $(B,g_B)$
complete compact indefinite is warped complete.
\item[$(Q3)$] A
complete compact indefinite manifold which is globally conformal
to a complete one (in particular, to a Lorentzian spaceform) is
complete.
\end{itemize}
As a consequence of Theorem \ref{tdna-wac-coca} and the conformal
invariance of lightlike completeness in the compact case, Theorem
\ref{t-comp-conf-luminosa}, we have:

\begin{equation} \label{edna-wac-coca}
(Q1)\quad \Longrightarrow \quad (Q2)\quad
\Longrightarrow\quad (Q3).
\end{equation}
As pointed out in Remark \ref{r-cr}, there are tracks  suggesting
that $(Q1)$ is not true. But the other questions (and, specially, $(Q3)$
for Lorentzian spaceforms) remain  open.
\end{remark}
\subsection{GRW spacetimes}\label{subGRW} Here, we apply previous
results to a simple and important class of warped products, the
{\em GRW spacetimes}. Those results allow one to prove the stability
of completeness for the $C^0$ topology in the class of GRW
spacetimes (see Theorem \ref{t-estab-GRW}). We also study stability in the
class of all the spacetimes, applying general results which
generically yield  $C^1$ stability.

Following \cite{ARSgrg}, a {\em Generalized Robertson--Walker} (GRW)
\index{spacetime!Generalized Robertson--Walker or GRW} spacetime is a
warped product $I \times_f F$ \index{$I\times_f F$} with base
$(I,-dt ^{2} )$, $I\subset \R$ interval, fiber a Riemannian
manifold $(F,g_{F} )$ and any (positive) warping function $f$;
with natural identifications $g=-dt ^{2} + f^{2}(t) g_{F}$. They
generalize classical Friedman--Lemaitre--Robertson--Walker spacetimes
because the fiber is not necessarily a model space. For
relativistic motivations and several properties, see \cite{ARSgrg,Sa-grg1,Sa-jgp},
and for some generalizations to more fibers,
\cite{Sa-grg2,Unal1,Unal2}.

\paragraph{$C^0$ stability  in the GRW class.} Even though  the base of a GRW spacetime may be
non--complete, some simple considerations in addition to Theorem
\ref{twarpedcompl} allow one to characterize its completeness as
follows.

\begin{lemma}
\label{l-GRW}  Let $I\times_f F$ be a GRW spacetime and fix $c\in
I=(a,b)$. If $F$ is incomplete then the GRW is incomplete in the
three causal senses. Otherwise:

\begin{enumerate}
\item The GRW  is timelike complete if and only if
\begin{equation}
\label{1-compGRW} \int_{a}^{c} \frac{f}{\sqrt{1+f^{2}}}\ dt =
\int_{c}^{b} \frac{f}{\sqrt{1+f^{2}}}\ dt = +\infty .
\end{equation}
Otherwise, all the timelike geodesics (not tangent to $I$) are
incomplete.
\item The GRW  is lightlike complete if and only if
\begin{equation}
\label{2-compGRW} \int_{a}^{c} f dt = \int_{c}^{b} f dt = +\infty
\end{equation}
Otherwise, all lightlike geodesics are incomplete.
\item The GRW is spacelike complete if and only if either $f$
satisfies {\rm (\ref{2-compGRW})} or if, whenever
\[\int_{a}^{c} f dt< +\infty \qquad \hbox{(resp.}\;
\int_{c}^{b} f dt < +\infty \hbox{)}
\]
holds, the function $f$ is unbounded in $(a,c)$ (resp. $(c,b)$).
\end{enumerate}
\end{lemma}

\begin{remark} (1) Clearly, the first integrals in
\rf{1-compGRW} and \rf{2-compGRW} are related to the completeness of
the causal geodesic $(t(s), x(s))$ towards the past direction
$t(s) \searrow a$, and the second ones towards the future
$t(s)\nearrow b$.

(2) Timelike $\Rightarrow$ lightlike $\Rightarrow$ spacelike
completeness. There are counterexamples to the  converses, but
these converses do hold if $f$ is bounded.
\end{remark}

Now, recall that: (i) both the complete and the incomplete
Riemannian metrics on $F$ are  $C^0$ open subsets of the set of
all the metrics on $F$, and (ii) the conditions on $f$ which
characterize causal completeness are both, open and closed  in the
set of all the positive warping functions on $I$. Therefore, one
can prove that each one of the subsets of timelike, lightlike and
spacelike geodesically complete GRW metrics on $I\times F$, are
open and closed in the set ${\cal GRW}(I\times F)$\index{${\cal GRW}(I\times F)$}
of all the GRW metrics on $I\times F$. Summing up (see \cite{Sa-grg1}):

\begin{theorem}
\label{t-estab-GRW} Both, completeness and incompleteness of each
causal type, are $C^r$ stable for all $r\geq 0$ in the set ${\cal
GRW}(I\times F)$.
\end{theorem}

Finally, we remark  that Lemma \ref{l-GRW} not only characterizes
incompleteness but also it says exactly which causal geodesics are
incomplete: if $I\neq \R$ all timelike geodesics are incomplete
either to the future or to the past, if $I=\R$ the geodesics
tangent to the base are complete, but all the other timelike ones
are either complete or incomplete, etc. So, if we call $Cau(M)$\index{$Cau(M)$}
the subset of the set of directions $S_RM$ (see Remark
\ref{pincomp-nulldirect}) which are causal, we have:

\begin{proposition}
\label{prop-est-direcGRW} In any GRW spacetime, the set of
complete causal directions  is a closed subset of ${\rm Cau}(M)$
(and $S_RM$). Thus, the incompleteness of causal geodesics is a
stable property in the set of causal directions for a given GRW
metric.
\end{proposition}

\begin{remark} Nevertheless, in general
the set of complete directions is not open in Cau$(M)$ or $S_R(M)$
(assume timelike incompleteness with $I=\R$) nor closed in $S_RM$
(take any spacelike complete non--lightlike complete GRW).
\end{remark}

\paragraph{Stability in Lor(M).} Recall that Theorem \ref{t-estab-GRW} only ensures
the stability of completeness and incompleteness in the class of
GRW spacetimes. We can wonder if the stability for a given GRW
spacetime will hold in the class Lor$(M)$ of all the Lorentzian
metrics on $M=I\times F$. This study (which can be carried out
mainly for causal geodesics), is  explained in detail in the book
\cite[Chapter 7]{BEE}, see also \cite{Be, BeEh-Cambridge}. So,
here we summarize the results very briefly. Generically, such
results yield $C^1$ stability, as geodesic equations (Christoffel
symbols) depend on the first derivatives of the metric tensor.
Nevertheless, under a simplifying technical hypotheses on the
fiber, the following $C^0$ result can be obtained.

\begin{theorem}
Let $I\times F$ be a GRW spacetime with homogeneous fiber
$(F,g_F)$. Then  timelike  incompleteness is $C^0$ stable in
Lor$(M)$.
\end{theorem}

\begin{proof} (Sketch).
Firstly recall that all GRW spacetimes with complete fiber are
globally hyperbolic, being the slices at constant $t$ Cauchy
hypersurfaces\footnote{In fact, any GRW spacetime is stably
causal, as it admits the time function $t$. If $g_F$ is complete,
then it is globally hyperbolic, because $t$ becomes a Cauchy time
function.} (see \cite[Section 3.11]{MS}). Then, one takes an adapted
auxiliary Riemannian norm, and estimates upper bounds for the
increasing of the  length of any timelike curve when each two such
slices are crossed (so, one proves that
some timelike geodesics are of finite length).
\end{proof}

For lightlike geodesics more accurate estimates are necessary and,
so, $C^1$ stability would be needed. Nevertheless, it is possible
to prove $C^1$ stability in the more general setting of
non--partially imprisoned geodesics.

\begin{theorem}\label{t-estC1}
Let $(M,g)$ be a semi--Riemannian manifold such that an inextendible
geodesic $\gamma: (a,b)\rightarrow M$ exists which is: (i)
forward--incomplete ($b<+\infty$), and (ii) not partially
imprisoned in a compact subset as $t\rightarrow b$.\\
Then there is a $C^1$ neighborhood $U(g)$ of $g$ such that all
$g_1\in U(g)$ admits an incomplete geodesic of the same causal
type of $\gamma$.
\end{theorem}

Now, notice that in all strongly causal spacetimes (in particular,
any GRW), timelike and lightlike geodesics are never partially
imprisoned\footnote{Nevertheless, all spacelike geodesics may be
imprisoned when $F$ is compact and $f$ is unbounded. In fact, this
happens in the de Sitter spacetime\index{spacetime!de Sitter}
$\SSS^n_1$, which can be written as
a warped product with fiber the round sphere and warping function
equal to $\cosh$.}.

\begin{corollary}
In any GRW spacetime, both lightlike and timelike incompleteness
are $C^1$ stable in Lor$(M)$.
\end{corollary}

Finally, a natural assumption for the stability of causal
completeness is the existence of a {\em pseudoconvex nonspacelike
geodesic system}. \index{pseudoconvex nonspacelike geodesic
system} This means that for each compact subset $K$ there is a
second compact subset $H$ such that each causal geodesic with both
endpoints in $K$ lies in $H$. Any globally hyperbolic spacetime
contains such a system and, thus:

\begin{theorem}\label{t-est-compl-C1}
Causal completeness is $C^1$ stable for any globally hyperbolic
spacetime. In particular, this holds in any GRW spacetime with
complete fiber $F$.
\end{theorem}

\subsection{Stationary spacetimes}\label{subStat}
A connected Lorentzian manifold $(M,g)$ is called a {\em
stationary spacetime} \index{spacetime!stationary} if it admits a
timelike Killing vector field $K$. Implicitly, we assume that
such a {\em stationary vector field} $K$ has been chosen, and
consider that it time--orientates the manifold. When the orthogonal
distribution $K^\perp$ to $K$ is integrable, the spacetime is
called {\em static} \index{spacetime!static} and admits a
natural local warped structure associated to $K$. The completeness
of stationary spacetimes has been already studied (see Proposition
\ref{p-pAMS}, Theorem \ref{tgeod1}). In Subsection
\ref{secstatic}, its geodesic connectedness will be studied in
depth. Now, we will make some remarks on its structure and will deduce
explicitly its geodesic equations, extending O'Neill's formulas for warped products.

 \paragraph{Standard stationary spacetimes.} Let
$(M_{0},\langle\cdot , \cdot\rangle_{R}) $ be a Riemannian
manifold,
and $\delta$ and $\beta$ a vector field and a positive smooth
function on $M_{0}$, respectively. A {\it standard stationary
spacetime} \index{spacetime!standard stationary} is the product
manifold $M=\R \times M_{0}$ endowed with the Lorentz metric,
under natural identifications:
\begin{equation}
\label{e-standardstationary} \langle\cdot ,
\cdot\rangle_L=-\beta(x)dt^{2}+\langle\cdot ,
\cdot\rangle_{R}+2\langle\delta(x),\cdot\rangle_{R}dt.
\end{equation}
When the cross term vanishes ($\delta \equiv 0$) the spacetime is
called standard static\index{spacetime!standard static}. This is a
warped product with Riemannian base and negative definite fiber
$M_0 \times_{\sqrt{\beta}} \R$.

Every stationary spacetime is locally a standard stationary one
with $K=\partial_t$. Even more, any globally hyperbolic spacetime
admits a spacelike Cauchy hypersurface $S$ (\cite{BeSa03}, see
\cite{MS}) and, if the spacetime is stationary with a complete
Killing vector field $K$\index{vector field!complete}
(i.e., $K$ has integral curves defined on the whole real line),
 $S$ can be moved by
means of the flow of $K$, so that the splitting \eqref{e-standardstationary} 
is obtained (see
\cite[Theorem 2.3]{CFS06} for details):

\begin{theorem} \label{t0}
A globally hyperbolic stationary spacetime is a standard
stationary one, if and only if one of its timelike Killing vector fields $K$
is complete.
\end{theorem}

Nevertheless, if $K$ is static the spacetime may be non--standard
static, as no Cauchy hypersurface orthogonal to $K$ may exist.

\paragraph{Levi--Civita connection and geodesic equations.}
 From now on, $\nabla$ and $\nabla^{R}$\index{$\nabla^R$} will denote the Levi--Civita
connection of $\langle\cdot , \cdot\rangle_L$ and $\langle\cdot ,
\cdot\rangle_{R}$, respectively. For each vector field $V$ on
$M_{0}$, $V\in \mathfrak{X} (M_{0})$, its  lifting to $M$ will be
denoted $\overline{V}$: that is, $\overline{V}_{(t,x)}=V_{x}$
for all $(t,x)\in M$ (analogously, if necessary, for a vector
field on $\R$). Put
\[
\Lambda(x) =\ -\ \frac{1}{\beta(x) + \langle\delta(x),\delta(x)\rangle_{R}}\qquad
\hbox{for all $x \in \mo$.}
\]
An explicit computation from Koszul formula  yields (see detailed
computations for this section in \cite{FS-IJTP}):
$$
\begin{array}{rcl}\displaystyle
\nabla_{\partial_{t}}\partial_{t}&= &\displaystyle
-\frac{1}{2}\Lambda\langle\delta,\nabla^{R}\beta\rangle_{R}\partial_{t}
+\frac{1}{2}\Lambda\langle\delta,\nabla^{R}\beta\rangle_{R}\delta
+\frac{1}{2}\nabla^{R}\beta,
\\\displaystyle
\nabla_{\overline{V}}{\overline{W}}&= &\displaystyle \frac{1}{2}\Lambda(\langle
W,\nabla_{V}^{R}\delta\rangle_{R}+\langle
V,\nabla_{W}^{R}\delta\rangle_{R})\partial_{t}
\\ &&\displaystyle
-\frac{1}{2}\Lambda(\langle
W,\nabla_{V}^{R}\delta\rangle_{R}+\langle
V,\nabla_{W}^{R}\delta\rangle_{R})\delta+\nabla_{V}^{R}W,
\\\displaystyle
2\nabla_{\overline{V}}\partial_{t}&= &\displaystyle
2\nabla_{\partial_{t}}\overline{V}=
 -\Lambda(V(\beta)+\langle
\delta,\nabla_{V}^{R}\delta\rangle_{R}-\langle
\nabla_{\delta}^{R}\delta,V\rangle_{R})\partial_{t}
\\ &&\displaystyle
+\Lambda(V(\beta)+\langle
\delta,\nabla_{V}^{R}\delta\rangle_{R}-\langle
\nabla_{\delta}^{R}\delta,V\rangle_{R})\delta
\displaystyle +\nabla_{V}^{R}\delta
-\langle \nabla^{R}_{(\cdot)}\delta,V\rangle_{R}^{\natural}
\end{array}
$$
for any $V,W\in \mathfrak{X} (M_{0})$,
where $^{\natural}$ denotes the vector field on $M_{0}$ metrically
associated to the corresponding 1--form (that is,
 $\langle Y,\langle \nabla^{R}_{(\cdot)}\delta,V\rangle_{R}^{\natural}\rangle_{R}= \langle \nabla_{Y}^{R}\delta,V\rangle_{R}$
 for any $Y\in \mathfrak{X}(M_{0})$).

For the geodesic equations, recall that $q=\langle \gamma ',\gamma'\rangle_L $
is constant for any geodesic $ \gamma(s)=(t(s),x(s))$,
$s\in I$ and,  as $\partial_{t}$ is a Killing vector field,  $C_\gamma
=\langle \gamma ',\partial_{t}\rangle_L$ is constant, too. That is,
one has the relations:
\begin{equation}\label{e8} \left\{ \begin{array}{l}
\langle (t',x'),(t',x')\rangle_L=-\beta(x) (t')^{2}+
2\langle \delta(x),x'\rangle_{R}t'
+\langle x',x'\rangle_{R}=q, \\
\langle \partial_{t},(t',x')\rangle_L=-\beta(x) t'
+\langle \delta(x),x'\rangle_{R}=C_\gamma.
\end{array} \right. \end{equation}
Let $\mathfrak{X}_{r,s}(M_0)$ \index{$\mathfrak{X}_{r,s}(M_0)$} be
the space of $r$--contravariant, $s$--covariant tensor fields on
$M_0$ ($\mathfrak{X} (M_0)\equiv \mathfrak{X}_{1,0} (M_0)$;
$\mathfrak{X}_{1,1}(M_{0})$ is identifiable to the space of
endomorphism fields). From the formulas of the connection,
for the $x$ part one has:
\begin{equation}\label{e9}
\nabla_{x'}^{R}x'=(t')^{2}R_{0}(x)+t'R_{1}(x,x')+R_{2}(x,(x',x'))
\end{equation}
with $R_0 \in \mathfrak{X}(M_{0})$, $R_1 \in
\mathfrak{X}_{1,1}(M_{0})$, $R_2 \in \mathfrak{X}_{1,2}(M_{0})$,
so that:
\begin{equation} \label{eee} \left\{
\begin{array}{l}R_{0}(x)=-\frac{1}{2}\left( \Lambda \langle \delta,
\nabla^{R}\beta\rangle_{R}\delta+\nabla^{R}\beta \right)(x) \\
R_{1}(x,x')= -\Lambda \left(\langle \nabla^{R}\beta,x'\rangle_{R}+
{\rm rot}\delta (x',\delta)\right) \delta(x)+
{\rm rot}\delta (\cdot,x')^{\natural}  \\
R_{2}(x,(x',\overline{x}'))= \Lambda {\rm
Sym}\nabla^{R}\delta(x',\overline{x}')\delta(x).
\end{array} \right.
\end{equation}

\begin{remark}
$R_0$ vanishes if $\beta$ is constant, $R_1$ reduces to $-\Lambda
\langle \nabla^{R}\beta,x'\rangle_{R} \delta(x)$
 if $\delta$ is irrotational, and $R_2$ vanishes if $\delta$ is Killing.
\end{remark}

By substituting in (\ref{e9}) the value of $t'$ from the second
equation in (\ref{e8}), a second order equation  for the spacelike
component $x(s)$  is obtained. Then, the first relation (\ref{e8})
can be regarded as a first integral of this equation.
 Moreover, a geodesic can be reconstructed for any solution of this differential equation.
Summing up:

\begin{theorem}\label{t2} Consider a curve $\gamma (s)=(t(s),x(s)), s\in I\subset \R$, 
in a stationary spacetime $(\R\times M_{0},\langle \cdot , \cdot \rangle_L)$. The curve $\gamma$ is a geodesic if and only if
 $C_\gamma=\langle \gamma',\partial_{t}\rangle_L$ is a constant and, 
 for such a $C_\gamma$, $x(s)$ satisfies
\begin{equation}\label{e10} \nabla_{x'}^{R}x'=\overline{R}_{0}(x)+\overline{R}_{1}(x,x')+\overline{R}_{2}(x,x'\otimes x')
\end{equation}
where
\[
\begin{split}
&\overline{R}_{0}(x)\ =\  \frac{C_\gamma^2}{\beta^2(x)} R_{0}(x),  \\
& \overline{R}_{1}(x,x')\ =\  -\ \frac{C_\gamma}{\beta(x)}
 \left( 2\frac{\langle
\delta(x),x'\rangle_R}{\beta(x)}R_{0}(x)+
R_{1}(x,x')\right) \\
&\overline{R}_{2}(x,x'\otimes x')\ =\ \frac{(\langle
\delta(x),x'\rangle_{R})^2}{\beta^{2}(x)}R_{0}(x)+ \frac{\langle
\delta(x),x'\rangle_R}{\beta(x)}R_{1}(x,x')+R_{2}(x,(x',x')),
\end{split}
\]
being $R_0, R_1, R_2$ as in (\ref{eee}).
\end{theorem}

\begin{remark} \label{rem3}

Notice that for any solution $x(s)$ of (\ref{e10}), the second
equation (\ref{e8}) yields the value of $t(s)$. Moreover, all the
geodesics can be reparametrized in such a way that $C_\gamma = 0,1$.
Thus, we have the following two cases:

(a) Case $C_\gamma=0$. The geodesic $\gamma$ is always spacelike and
orthogonal to $\partial_t$, and we can also assume that the value
of $q$ is fixed equal to 1. Recall that in this case
$\overline{R}_{0}=\overline{R}_{1}=0$. In the static case,
equation (\ref{e10}) is just the equation of the geodesics of
$(M_0, \langle \cdot,\cdot\rangle_R)$; thus, these geodesics can
be regarded as trivial. Nevertheless, when $\delta $ is not null
the differential equation becomes:
$$\nabla_{x'}^{R}x'=\overline{R}_{2}(x,x'\otimes x').$$

(b) Case $C_\gamma=1$. Equation (\ref{e10}) becomes rather complicated in
general. In the static case, $x(s)$ satisfies the equation of a
classical Riemannian particle under the potential $V=-1/2\beta $,
and the constant $q/2$ is the classical energy (kinetic plus
potential) of this particle.  In
the general case this interpretation does not hold,
even though $q$ is a constant of the
motion.
\end{remark}

As we have seen, the direct study of the geodesic equations become
extremely complicated, except if some simplifying assumptions are
carried out. The simplest one is to consider the static case
$\delta \equiv 0$, where previous equations reduce to a
Riemannian dynamical system for $x$ and a reparametrization for
$t$. Nevertheless, even in this case the interplay between $x$ and the
reparametrization of $t$
may yield a non--trivial problem. In fact, in the topic of
connecting by a geodesic two fixed
 endpoints $(t_1,x_1), (t_2,x_2)$, the existence of a solution
 $x(s)$ of (\ref{e10}) (a simple matter in the static
 case) which connects $x_1$ and $x_2$ yields a $t(s)$ which must be
 controlled in order to connect $t_1, t_2$. This makes necessary
 the variational approach in Subsection \ref{secstatic} below.

\section{Variational approaches  in Lorentzian manifolds}\label{s4}

Due to the increasing knowledge in variational methods applicable
to the study of critical points of unbounded functionals, at the
end of 80's Benci, Fortunato, Giannoni and Masiello improved the
investigations on the existence of geodesics in a Lorentzian
manifold by means of variational tools. Next, we will describe the
approach first for geodesics connecting two points in the
Riemannian case, then in the two most characteristic
Lorentzian cases: stationary and orthogonal splitting.

\subsection{The Riemannian framework for geodesics connecting two points}
\label{secriemannian}

In order to apply global variational tools  to the action
functional $f$ as defined in
(\ref{action})\index{functional!action}, some ``good'' abstract
arguments and a suitable variational framework are required. We
introduce both here in the Riemannian setting not only for
completeness, but also for pedagogical reasons.

Let us remark that, in general, we are interested in the study of
critical points of a $C^1$ functional $J$ defined on a Riemannian manifold $\Omega$,
so a way to check the existence of critical levels is
condition $(C)$\index{condition!$(C)$} introduced by Palais in \cite{Pa63}.
Roughly speaking, it says that {\sl ``each subset of $\Omega$ in which $J$ is bounded
and its differential can tend to zero, has to have a critical point in its closure''}.
But, starting from this definition, the true problem is to test the existence of such a
critical point and a way to help in this ``research'' comes from the most used
Palais--Smale condition\index{condition!Palais--Smale or $(PS)$}.

\begin{definition}
Let $\Omega$ be a (possibly infinite--dimensional) Riemannian manifold
and $J \in C^1(\Omega,\R)$.
Functional $J$ satisfies the {\sl Palais--Smale condition} on $\Omega$,
briefly $(PS)$, if any sequence $(x_k)_k \subset \Omega$ such that
\[
(J(x_k))_k\; \mbox{is bounded}\quad\mbox{and}\quad
\lim_{k \to +\infty}J'(x_k) = 0
\]
has a subsequence converging in $\Omega$.
\end{definition}

Thus, once that condition $(PS)$ is satisfied, a bounded--from--below functional
has at least a critical point, its minimum. More precisely, a minimum
theorem can be stated as follows (see, e.g, \cite{Ra}).

\begin{theorem}
\label{min}
Let $\Omega$ be a complete Riemannian manifold
and $J$ a $C^1$ functional on $\Omega$ which satisfies condition $(PS)$.
If $J$ is bounded from below then it attains its infimum.
\end{theorem}

The power of condition $(PS)$ is that not only it allows one to
prove that an infimum is a minimum but also that more than one critical
point can be found if the topology of $\Omega$ is rich enough.
To this aim, different topological tools can be used but, here, for simplicity,
we just introduce the Ljusternik--Schnirelman Theory\index{Theory!Ljusternik--Schnirelman}
and its main arguments
(for more details, see, e.g., \cite{Amb,Ma,Pa66}).

\begin{definition}
Let $\Omega$ be a topological space. Given $A \subset \Omega$, the
{\sl Ljusternik--Schnirelmann
category}\index{category!Ljusternik--Schnirelmann} of $A$ in
$\Omega$, briefly $\cat_{\Omega} (A)$\index{$\cat_{\Omega} (A)$},
is the least number of closed and contractible subsets of $\Omega$
covering $A$. If it is not possible to cover $A$ with a finite
number of such sets, then $\cat_{\Omega} (A) = +\infty$.
\end{definition}
For simplicity, we denote $\cat (\Omega)$ = $\cat_{\Omega} (\Omega)$.

\begin{theorem}
\label{LjuSchni} Let $\Omega$ be a Riemannian manifold
and $J$ a $C^1$ functional on $\Omega$ which satisfies
condition $(PS)$. Given any $k \in \N$, $k > 0$,
let us define
$$
c_k = \inf_{A \in \Gamma_k} \sup_{x \in A} J(x) \quad
\mbox{with}\;\; \Gamma_k = \{ A \subset \Omega : \cat_\Omega (A)
\ge k\} .
$$
If $\Omega$ is complete, then $c_k$ is a critical value of $J$ for each
$k$ such that $\Gamma_k \ne \emptyset$ and $c_k \in \R$. Moreover,
if $J$ is bounded from below but not from above and $\cat(\Omega)=+\infty$,
then a sequence $(x_k)_k$ of critical points exists such that
$\displaystyle{\lim_{k \to +\infty} J(x_k) = + \infty}$.
\end{theorem}

\begin{remark}
The same results of Theorems \ref{min} and \ref{LjuSchni} still hold if
the completeness of $\Omega$ is replaced by the completeness of each sublevel
of $J$ or if condition $(PS)$ is replaced by weaker assumptions (see, e.g.,
\cite{Ce78} or also \cite[pp. 80]{St}).
\end{remark}

Now, in order to introduce the correct functional framework (for more details,
see, e.g., \cite{Pa63}),
if $(\M,g)$ is a smooth $n$--dimensional connected semi--Riemannian
mani\-fold, let us define $H^1(I,\M)$\index{$H^1(I,\M)$} the set of curves
$z : I \to \M$ such that for any local chart $(U,\varphi)$ of $\M$, with
$U \cap z(I) \ne \emptyset$, the curve
$\varphi \circ z$ belongs to the Sobolev space
$H^1(z^{-1}(U),\R^n)$. Taking $I=[0,1]$ (without loss of generality,
as the set of geodesics is invariant by affine
reparametrizations), it is well known that $H^1(I,\M)$ is equipped with a
structure of infinite dimensional manifold modelled on the Hilbert space
$H^1(I,\R^n)$, as, if $z \in H^1(I,\M)$, the tangent space\index{$T_zH^1(I,\M)$}
to $H^1(I,\M)$ at $z$ can be written as
\[
T_zH^1(I,\M) \equiv
\{\zeta \in H^1(I,T\M) : \pi \circ \zeta = z\},
\]
being $T\M$ the tangent bundle of $\M$ and
$\pi : T\M \to \M$ the corresponding bundle projection. In other words,
$T_zH^1(I,\M)$ is the set of the vector fields along $z$ whose components
with respect to a local chart are functions of class $H^1$.\\
It is easy to check that the action functional\index{functional!action} $f$ in (\ref{action})
is well defined on $H^1(I,\M)$; moreover, it is
at least of class $C^1$ and there results
\begin{equation}\label{derivative}
f'(z)[\zeta] =
2\ \int_0^1 g(z(s))[z'(s),\frac{D\zeta}{ds}]\ ds 
\end{equation}
if $z \in H^1(I,\M)$,
 $\zeta \in T_z H^1(I,\M)$.

Here in the following, we focus
only on one model problem, for example the existence of geodesics joining two
fixed points (as already remarked, different submanifolds have to be defined if
other ``boundary'' conditions are required).

In this case, the ``natural'' setting
of this variational problem is a suitable infinite dimensional submanifold of
$H^1(I,\M)$. More precisely, fixed $p$, $q \in\M$, we can consider\index{$\Omega^1(p,q;\M)$}
\begin{equation}
\label{omega}
\Omega^1(p,q;\M) =
\{z \in H^1(I,\M):\; z(0) = p,\, z(1) = q\}
\end{equation}
such that, if $z \in \Omega^1(p,q;\M)$, the tangent space\index{$T_z\Omega^1(p,q;\M)$} to
$\Omega^1(p,q;\M)$ at $z$ is given by
\[
T_z\Omega^1(p,q;\M) =
\{\zeta \in T_zH^1(I,\M) :\; \zeta(0) = 0 = \zeta(1)\}\ .
\]

Classical boot--strap arguments allow one to prove
the following variational principle\index{variational principle!classical}:

\begin{theorem}
\label{classic}
The curve $z:I \to\M$ is a geodesic joining $p$ to $q$ in $\M$ if and only if $z \in \Omega^1(p,q;\M)$
is a critical point of $f$ on $\Omega^1(p,q;\M)$.
\end{theorem}

Two particular subcases are of special interest::
\begin{itemize}
\item[{\sl (i)}] $\M \equiv \R$ is the 1--dimensional Euclidean
space; \item[{\sl (ii)}] $(\M,g) \equiv (\mo,\met_R)$ is a
Riemannian manifold.
\end{itemize}

In case {\sl (i)}, fixed any $t_p$, $t_q \in \R$ the set defined in (\ref{omega})\index{$W^1(t_p,t_q)$}
becomes
\[
W^1(t_p,t_q) =\{t\in H^1(I,\R): t(0)=t_p,\; t(1)=t_q\} = H_0^1(I,\R) + j^*,
\]
with $H_0^1(I,\R) =\{\tau\in H^1(I,\R): \tau(0) = 0 =\tau(1)\}$\index{$H_0^1(I,\R)$}
vector space and\index{$j^*$}
\[
j^* : s\in I \mapsto t_p + s \Delta_t \in\R,\qquad \Delta_t = t_q - t_p.
\]
Whence, $W^1(t_p,t_q)$ is a closed affine submanifold of $H^1(I,\R)$
with tangent space $T_tW^1(t_p,t_q) \equiv H_0^1(I,\R)$ for all $t \in W^1(t_p,t_q)$.

On the other hand, in case {\sl (ii)}, if $\mo$ is a smooth enough
Riemannian manifold (at least of class $C^3$), by means of the
Nash Embedding Theorem\index{Theorem!Nash} (see \cite{Na56}) we
can assume that $\mo$ is a submanifold isometrically embedded in
some Euclidean space
$\R^N$ and the embedding is closed in compact regions of $\mo$ 
(as will be the case in the proofs of the
results below). Thus, $\met_R$ is the restriction to $\mo$ of the standard
Euclidean metric of $\R^N$ with $d(\cdot,\cdot)$\index{distance}
being the corresponding distance, i.e.,
\begin{equation}
\label{distanza} d(x_1,x_2) = \inf \left\{  \int_a^b
\sqrt{\langle\gamma',\gamma'\rangle_R} ds : \; \gamma\in
A_{x_1,x_2}\right\}
\end{equation}
where $x_1$, $x_2 \in\mo$ and $\gamma\in A_{x_1,x_2}$, 
$\gamma: [a,b]\rightarrow\mo$ being any piecewise smooth curve
in $\mo$ joining $x_1$ to
$x_2$.
Hence, it can be proved that manifold $H^1(I,\mo)$ is a sub\-mani\-fold
of $H^1(I,\R^N)$ and can be
identified with the set of the absolutely continuous curves
$x: I\to \R^N$, $x(I) \subset \mo$, with square summable derivative.
Thus, fixed any $x_p$, $x_q \in \mo$, the set defined in (\ref{omega})\index{$\Omega^1(x_p,x_q;\mo)$}
can be identified as follows:
\[
\begin{split}
\Omega^1(x_p,x_q;\mo) = \{x:I \to\mo:\, & \hbox{$x$ is absolutely continuous and such that} \\
&\displaystyle x(0) = x_p,\, x(1) = x_q,\, \int_0^1\langle
x', x'\rangle_R ds < +\infty\}.
\end{split}
\]
Taking any $x \in \Omega^1(x_p,x_q;\mo)$, let us point out that
it is well defined
\[
\| x'\|^2 = \int_0^1 \langle  x', x'\rangle_R ds
\]
and the tangent space in $x$\index{$T_x\Omega^1(x_p,x_q;\mo)$} can be
identified with
\[\begin{split}
T_x\Omega^1(x_p,x_q;\mo) = &\{\xi:I \to T\mo:\
\xi(s) \in T_{x(s)}\mo\, \hbox{for all $s\in I$, $\xi$ is}\\
&\qquad \hbox{absolutely continuous and}\, \xi(0) = 0 =
\xi(1),\, \|\xi\|_* < +\infty\},
\end{split}
\]
being its Hilbert norm\index{$\|\cdot\|_*$}
\[
\|\xi\|_*^2\ =\
\int_0^1\langle\frac{D^R\xi}{ds},\frac{D^R\xi}{ds}\rangle_R ds
\]
(here, $D^R/ds$\index{$D^R/ds$} denotes the covariant derivative along $x$
relative to the Riemannian metric tensor $\met_R$).\\
Furthermore, let us point out that if $\mo$ is a complete
Riemannian manifold with respect to $\met_R$, also $H^1(I,\mo)$
and $\Omega^1(x_p,x_q;\mo)$
are complete Riemannian manifolds equipped with their scalar product.

The possibility to work in a submanifold of an Euclidean space $\R^N$
allowed Benci and Fortunato to prove a ``splitting'' lemma
useful in the proof of condition $(PS)$ as
it allows one to manage bounded sequences in $H^1(I,\mo)$ so to pass from
a ``weak'' limit to a ``strong'' one (see \cite[Lemma 2.1]{BF94}).

\begin{lemma}\label{splitting}
Let $\mo$ be a submanifold of ${\mathbb R}^N$ and
 $(x_k)_k \subset \Omega^1(x_p,x_q;\mo)$ a sequence so that
\begin{equation}\label{bounded}
(\| x'_k\|)_k\quad\hbox{is bounded.}
\end{equation}
Then, there exists $x \in H^1(I,\R^N)$ such that, up to subsequences, it is
\begin{equation}\label{limit}
x_k \rightharpoonup x\;\mbox{weakly in $H^1(I,\R^N)$,}\quad
x_k \rightarrow x \; \mbox{uniformly in $I$.}
\end{equation}
If $\mo$ is complete\footnote{We emphasize the following technical step. 
Nash's embedding is not closed if $\mo$ is not compact.
Nevertheless, if $\mo$ is complete, the bounded sequence $(x_k)_k$ 
lies in a compact subset of $\mo$, and this suffices.}, 
then $x \in \Omega^1(x_p,x_q;\mo)$;
furthermore, there
exist two sequences $(\xi_k)_k$, $(\nu_k)_k \subset H^1(I,\R^N)$
such that
\begin{equation}\label{sequences} \begin{split}
&\xi_k \in T_{x_k}\Omega^1(x_p,x_q;\mo), \quad
x_k - x = \xi_k + \nu_k \qquad \mbox{for all $k\in\N$,}\\
&\xi_k \rightharpoonup 0\;\; \mbox{weakly}\quad\mbox{and}\quad \nu_k \rightarrow 0
\;\;\mbox{strongly in $H^1(I,\R^N)$.}
\end{split}
\end{equation}
\end{lemma}

\noindent
At last, in this setting, we are able to estimate the Ljusternik--Schnirelman
category of the space of curves $\Omega^1(x_p,x_q;\mo)$
(cf. \cite{FH91}).

\begin{proposition}
\label{FadellHusseini}
If the Riemannian manifold $\mo$ is not contractible in itself
then $\Omega^1(x_p,x_q;\mo)$ has infinite category and possesses
compact subsets of arbitrarily high cate\-gory.
\end{proposition}

Now, just in order to apply all the above arguments in a very
simple case, let us consider $(\mo,\met_R)$, a smooth Riemannian
manifold so that Nash Embedding Theorem\index{Theorem!Nash} holds.
The well known Hopf--Rinow Theorem\index{Theorem!Hopf--Rinow}
states that completeness implies geodesic connectedness (see
\cite{HR31} or also \cite[pp. 4]{BEE}), and it would be also a
consequence of Weirstrass theorem for reflexive spaces (see, e.g.,
\cite[Theorem 1.2]{St}). Here, we give an alternative proof by
means of the variational setting and the topological arguments
introduced in this section. Furthermore, such tools allow one to
obtain also a multiplicity result under a minimum topological
assumption  (but which implies that the $\Omega^1(x_p,x_q;\mo)$
does have a ``rich'' topology), the  non--contractibility for
$\mo$. This result cannot be deduced by more classical geometric
methods (as those explained in Subsection \ref{sub1.2}), which
allow to obtain just one (minimizing) geodesic in each homotopy
class.

\begin{theorem}\label{hopf}
If $\mo$ is complete as metric space equipped
with $d(\cdot,\cdot)$ in (\ref{distanza}), then it
is also geodesically connected\index{geodesic
connectedness!Riemannian manifold}, i.e., each couple of its
points can be joined by a geodesic. Furthermore, if $\mo$  is not
contractible in itself, then each couple of its points can be
joined by infinitely many geodesics.
\end{theorem}

\begin{proof}
Our aim is to prove that, fixed $x_p$, $x_q \in \mo$, one or
more geodesics exist joining $x_p$ to $x_q$ in $\mo$. Or better,
by means of Theorem \ref{classic}, it is enough to investigate the
existence of critical points of the energy
functional\index{functional!energy} in (\ref{action}) with $g =
\met_R$, i.e.,
\begin{equation}\label{energy}
f(x) = \int_0^1 \langle x', x'\rangle_R ds\quad\hbox{on
$\Omega^1(x_p,x_q;\mo)$.}
\end{equation}
In order to apply Theorem \ref{min}, we just need to prove that
$f$ satisfies condition $(PS)$. So, let
$(x_k)_k \subset \Omega^1(x_p,x_q;\mo)$ be a sequence so that
$(f(x_k))_k$ is bounded and $f'(x_k) \to 0$ as $k \to +\infty$.
Clearly, by (\ref{bounded}), (\ref{energy}) and the completeness of $\mo$,
Lemma \ref{splitting} applies and, up to subsequences, $x \in \Omega^1(x_p,x_q;\mo)$
exists so that (\ref{limit}) holds. Moreover, $x_k - x \in H^1(I,\R^N)$
splits so that (\ref{sequences}) is satisfied. Now, we have only
to prove that $x_k \to x$ strongly in $H^1(I,\R^N)$; whence,
to prove that $\xi_k \to 0$ strongly in $H^1(I,\R^N)$. But we know that
$\xi_k \in T_{x_k}\Omega^1(x_p,x_q;\mo)$, so, being $x_k = x + \xi_k +\nu_k$,
from (\ref{derivative}) and (\ref{sequences}) if follows that $\varepsilon_k \searrow 0$
exists so that
\[
\varepsilon_k = \int_0^1 \langle x'_k, \xi'_k\rangle_R ds =
\int_0^1 \langle x', \xi'_k\rangle_R ds + \int_0^1 \langle
\xi'_k, \xi'_k\rangle_R ds + \int_0^1 \langle
\nu'_k, \xi'_k\rangle_R ds,
\]
where it is
\[\begin{split}
&\int_0^1 \langle x', \xi'_k\rangle_R ds \to 0\quad \hbox{for the weak convergence $\xi_k \wk 0$,}\\
&\int_0^1 \langle  \nu'_k, \xi'_k\rangle_R ds \to 0\quad
\hbox{for the strong convergence $\nu_k \to 0$.}
\end{split}
\]
Thus, it has to be also
\[
\int_0^1 \langle \xi'_k, \xi'_k\rangle_R ds \to 0;
\]
hence, $\xi_k \to 0$ strongly in $H^1(I,\R^N)$. So, $x_p$ and $x_q$
must be joined by a geodesic which minimizes the (positive) energy functional $f$.\\
On the other hand, if $\mo$ is not contractible in itself, then
the existence of infinitely many critical points of $f$ in $\Omega^1(x_p,x_q;\mo)$
is a direct consequence of Proposition \ref{FadellHusseini} and
Theorem \ref{LjuSchni}.
\end{proof}

\begin{remark} (1) Let us point out that the completeness is not 
a necessary condition for
the geodesic connectedness. In fact, if $\mo$ is not complete but
some convexity assumptions on its boundary are satisfied, then the
geodesic connectedness still holds (see \cite{BGS02}).

(2) A similar approach works for the problem of the existence of a
closed geodesic. Nevertheless, here the space of curves is formed
by loops which are not attached to any point and, so, one needs
either the compactness of the manifold or some assumption at
infinity, in order to ensure Palais Smale condition  (e.g., see
\cite{BG92,CS94} for two different sets of hypotheses which
overcome the lack of compactness, or the careful discussion in
\cite{BGS-advan}). Moreover, a topological hypothesis (as
compactness) is also needed even for the problem of the existence
of one closed geodesic: notice that the minimum critical points of
the energy functional are always constant geodesics which cannot
be accepted as solutions of the problem. As far as we know, the
general multiplicity problem remains open  (see \cite{GrMe69}).
The difficulty appears because infinite critical points of the
functional for closed geodesics may mean a single closed geodesic
which is run many times.
In fact,  topological assumptions (which implied the existence of
critical points with arbitrary lengths of $f$) yield
multiplicity of solutions in the case of geodesic connectedness,
but in the case of closed geodesics only ensure the existence of a
non-trivial one.
\end{remark}

\subsection{Variational principles for static
and stationary spacetimes:  extrinsic and intrinsic
approaches}\label{secstatic}

Differently from the Riemannian case, for Lorentzian manifolds
something like Theorem \ref{hopf} does not hold.
In fact, a counterexample is
given by the anti--de Sitter spacetime\index{spacetime!anti--de Sitter}
$\M =\ ]-\frac{\pi}{2},\frac{\pi}{2}[ \times \R$ equipped with the Lorentzian
metric
\[
\langle\cdot,\cdot\rangle_L =\frac{1}{\cos^{2}x}\ \left(-\ dt^{2} + dx^2\right)
\]
which is geodesically complete (and static), but not geodesically
connected (cf. \cite{Pe72}). Moreover, in a Lorentzian manifold
the minimizing arguments in Theorem \ref{min} cannot be used
directly as the action functional is strongly unbounded both from
below and from above.

Anyway, in some model problems, several alternatives are possible. And, in
fact, two different techniques can be used  for stationary
spacetimes:
\begin{itemize}
\item[{\sl (i)}] the extrinsic approach, when the Lorentzian
manifold ``splits'' in a Riemannian manifold and an Euclidean
space and the metric coefficients are ``time independent'' (see,
e.g., \cite{BFG91,CFS03,CaSa00}); \item[{\sl (ii)}] the intrinsic
one, when such a splitting is not given ``a priori'' but a
timelike Killing vector field exists anyway (see
\cite{CFS06,CMP03, GP99}).
\end{itemize}

Even though the latter has revealed as the most powerful
\cite{CFS06}, in this section, we want to outline both these
techniques, showing the main used tools and their limits.
\smallskip

\paragraph{
Case {\sl (i)} (Extrinsic).} The main idea in this case is
transforming the indefinite problem in a subtler problem on a
Riemannian manifold. In order to give an idea of this way of work,
let us consider the simplest model of Lorentzian manifolds of this
kind. Recall from Subsection \ref{subStat}
that, explicitly, a standard static spacetime
\index{spacetime!standard static} is a connected Lorentzian
mani\-fold $(\M ,\met_L)$ with $\M =
 \R \times \mo$ and
\begin{equation}
\label{metrica} \met_L = - \beta(x)\ d t^2 + \met_R,
\end{equation}
for a  Riemannian manifold $(\mo,\langle \cdot,\cdot\rangle_R)$,
$t$ natural coordinate of $\R$, $x \in \mo$ and a smooth positive
function $\beta: \mo \to \R$.

\begin{remark}\label{rema4-11}
Recall that, essentially, the study of geodesic connectedness in
the standard case is equivalent to that one on any static
spacetime by using the universal covering (see \cite[Theorem
2.1]{Sa06} and \cite[Section 2]{BCFS03}). This does not hold by
any means in the stationary case. In fact, as far as we know, it
is not known if a compact stationary spacetime must be
geodesically connected. And this problem cannot be reduced to the
standard one by passing through the universal covering (the
3--sphere admits a stationary metric, see for example
\cite{RSgeoded}).
\end{remark}

 For fixed $p=(t_p,x_p)$, $q=(t_q,x_q) \in \M$, thanks to this
splitting we have that the manifold of curves joining $p$ to $q$
splits, also, so it is
\[
\Omega^1(p,q;\M) \equiv W^1(t_p,t_q)\times \Omega^1(x_p,x_q;\mo)
\]
and, consequently, also its tangent space in each $z=(t,x) \in
\Omega^1(p,q;\M)$ becomes
\[
T_z\Omega^1(p,q;\M) \equiv H_0^1(I,\R) \times
T_x\Omega^1(x_p,x_q;\mo).
\]
Hence, the action functional $f$ in (\ref{action}) becomes
\begin{equation}
\label{actionstatic} f(z) =\ \int_0^1\langle z',
z'\rangle_L\ ds\ =\ \int_0^1 \left(-\beta(x) (t')^2 + \langle
x', x'\rangle_R\right) ds
\end{equation}
if $z=(t,x) \in \Omega^1(p,q;\M)$, with Fr\'echet differential
\[
 f'(z)[\zeta]  =  -\ \int_0^1
\left(\beta'(x)[\xi]\ (t')^2 + 2 \beta(x) t' \tau'\right)
ds + 2 \int_0^1 \langle x',\frac{D^R\xi}{ds}\rangle_R ds
\]
for all $\zeta=(\tau,\xi)\in T_z\Omega^1(p,q;\M)$, where
$\beta'(x)[\xi]$ denotes the differential of $\beta$ on $\xi$ at
$x\in \mo$.

As Benci, Fortunato and Giannoni showed in
their pioneer work \cite[Theorem 2.1]{BFG91},
the following new variational principle\index{variational principle!in the static case}
can be stated.

\begin{proposition}
\label{newazione} Let $z^* = (t^*, x^*) \in \Omega^1(p,q;\M)$. The
following statements are equi\-valent:
\begin{itemize}
\item[{\sl (a)}] $z^*$ is a critical point of the action functional $f$ defined
as in (\ref{actionstatic});
\item[{\sl (b)}] $x^*$ is a critical point of the 
functional $\J: \Omega^1(x_p,x_q;\mo) \to \R$
defined as
\[
 \J(x) =\ \int_0^1 \langle x',x'\rangle_R\ ds\  -\ {\Delta_t^2}\ \left(\int_0^1
\frac{1}{\beta(x)}\ ds\right)^{-1}
\]
and $t^* = \Psi(x^*)$, with $\Delta_t^2 = (t_q-t_p)^2$ and
$\Psi:\Omega^1(x_p,x_q;\mo) \to W^1(t_p,t_q)$
such that
\[
\Psi(x)(s) = t_p +  \Delta_t\ \left(\int_0^1 \frac{1}{\beta(x)}\ ds\right)^{-1}
\ \int_0^s \frac{1}{\beta(x(\sigma))}\ d\sigma .
\]
\end{itemize}
Moreover, $f(z^*) = \J(x^*)$.
\end{proposition}

\begin{proof}
The main idea of this proof is to introduce the partial derivatives
\[
\frac{\partial f}{\partial t}(z)[\tau] = f'(z)[(\tau,0)],\quad
\frac{\partial f}{\partial x}(z)[\xi] = f'(z)[(0,\xi)],
\]
for any $z=(t,x) \in \Omega^1(p,q;\M)$, $\tau \in H_0^1(I,\R)$,
$\xi \in T_x\Omega^1(x_p,x_q;\mo)$, and consider the kernel of
$\frac{\partial f}{\partial t}$, i.e.,
\[
{\cal N} = \{z \in \Omega^1(p,q;\M):\; \frac{\partial f}{\partial t}(z) \equiv 0\}.
\]
Once proved that $z=(t,x) \in {\cal N}$ if and only if $t =
\Psi(x)$ (${\cal N}$ is the graph of the $C^1$ map $\Psi$), then it is enough to
define $\J = f|_{\cal N}$, i.e., $\J(x) = f(\Psi(x),x)$ for all $x
\in \Omega^1(x_p,x_q;\mo)$, and to remark that
\[
\J'(x)[\xi] = f'(\Psi(x),x)[(0,\xi)]\quad  \mbox{for all $\xi \in
T_x\Omega^1(x_p,x_q;\mo)$.}
\]
\end{proof}

Now, we can investigate the existence of critical points of $\J$ on $\Omega^1(x_p,x_q;\mo)$.

It is quite easy to prove that this new functional is bounded from below
if coefficient $\beta$ is bounded from above (see \cite{BFG91}),
but by means of a reduction to a $1$--dimensional problem
(following some ideas introduced in \cite{CFS003}), very careful
estimates on the diverging sequences of $\Omega^1(x_p,x_q;\mo)$
allow one to prove the following result (for all the details,
see \cite[Proposition 4.1]{BCFS03}).

\begin{proposition}
\label{inf} If coefficient $\beta $ grows at most quadratically at
infinity, i.e., there exist $\lambda \ge 0$, $\mu, k \in \R$,
$\alpha \in [0,2)$ and a point $\bar x\in\mo$ such
that\footnote{Obviously, one can put $\mu=0$ without loss of
generality here (as well as in formula (\ref{bounddelta}) below
one can put $\mu_1=0$). Nevertheless, in this fashion, one can
compare better this result with the previously obtained ones.}
\begin{equation}
\label{nonlimitato}
0 < \beta(x) \le \lambda d^2(x,\bar x) + \mu d^\alpha(x,\bar x) + k
\quad \mbox{for all $x \in \mo$,}
\end{equation}
then functional $\J$ is
\begin{itemize}
\item bounded from below;
\item coercive\index{coercivity}, i.e.,
\[
 \J(x) \to +\infty \qquad \mbox{if}\quad \| x'\| \to +\infty.
\]
\end{itemize}
\end{proposition}

Thus, if (\ref{nonlimitato}) holds, not only $\J$ has an infimum on
$\Omega^1(x_p,x_q;\mo)$ but any sequence
$(x_k)_k \subset \Omega^1(x_p,x_q;\mo)$, such that $(\J(x_k))_k$ is bounded
and $\J'(x_k) \to 0$, has to be bounded and Lemma \ref{splitting} can apply.
Then, by suitable computations (for more details, see, e.g., \cite[Proposition 4.3]{BCFS03}),
we are able to prove that:

\begin{proposition}
\label{palais}
In the hypotheses of Proposition \ref{inf}
functional $\J$ satisfies condition $(PS)$ in $\Omega^1(x_p,x_q;\mo)$.
\end{proposition}

At last, from Propositions \ref{inf} and \ref{palais} and Theorem \ref{min},
respectively \ref{LjuSchni}, it follows the following result (see \cite[Theorem 1.1]{BCFS03}).

\begin{theorem}
\label{maintheoremstatic} Let $(\M ,\met_L)$ be a standard static
Lorentzian manifold with $\M = \R\times\mo$ and $\met_L$ as in
(\ref{metrica}). If the smooth Riemannian manifold $(\mo,\met_R)$
is complete and coefficient $\beta$ in (\ref{metrica}) satisfies
the at most quadratic growth condition (\ref{nonlimitato}), then,
$\M$ is geodesically connected\index{geodesic connectedness!static
spacetime}.

 Furthermore, if $\mo$ is non--contractible in itself,
then any two points can be joined by a sequence of (spacelike)
geodesics $(z_k)_k$ with diverging lengths.
\end{theorem}

Let us point out that this existence result can be extended also to a non complete
static manifold but under suitable hypotheses on its boundary
(see \cite[Section 6]{BCFS03} and references therein) while the
growth assumption on $\beta$ cannot be further improved
as there are counterexamples to geodesic connectedness if $\beta$
grows more than quadratically at infinity.

\begin{example}\label{counter1}
For fixed $\eps>0$, consider the static spacetime
\[
(\R^{2},\langle\cdot,\cdot\rangle_{L}), \quad \met_L =
-\beta_{\eps}(x)dt^{2} + dx^{2},
\]
where $\beta_{\eps}(x)$ is a  smooth function such that
\[
\left\{
\begin{array}{lll}
\beta_{\eps}(x)=1+|x|^{2+\eps} & \hbox{if} & x\in\R\setminus (-1,1)
\\ \beta_{\eps}([-1,1])\subset [1,2] .
\end{array}
\right.
\]
It can be proved that a couple of points $p$, $q \in \R^2$ exists
which cannot be joined by means of a geodesic (for all the details, see \cite[Section 7]{BCFS03}).
\end{example}

\begin{remark}
Careful estimates allow one to investigate also the number of timelike geodesics
joining $p$ to $q$ according to $|\Delta_t| \nearrow +\infty$ (see, e.g., \cite[Section 5]{BCFS03}
and the related references).
\end{remark}

Similar arguments to those ones outlined for the standard static
spacetimes, can be also used for the study of geodesic
connectedness for the  larger class of standard stationary
spacetimes according to the definition in
(\ref{e-standardstationary}).

In fact, fixed $p=(t_p,x_p)$, $q=(t_q,x_q) \in \M$, also in this
setting a suitable new variational principle can be proved (see
\cite[Theorem 2.2]{GiaMa91} or also \cite[Theorem 3.3.2]{Ma}) so
to define a new ``Riemannian'' functional
\[
\begin{split}
\J_1(x) &= \int_0^1\langle x', x'\rangle_R\
ds\ +\ \int_0^1\frac{\langle \delta(x), x'\rangle_R^2}{\beta(x)}\ ds\\
&\;  -\  \left(\int_0^1 \frac{\langle \delta(x),
x'\rangle_R}{\beta(x)}\ ds - \Delta_t\right)^2 \ \left(\int_0^1
\frac1{\beta(x)}\ ds\right)^{-1}
\end{split}
\]
on $\Omega^1(x_p,x_q;\mo)$ which is bounded from below and coercive
under ``good'' growth hypotheses on coefficients
$\beta$ and $\delta$ in (\ref{e-standardstationary}). Thus, with this extrinsic
approach, also in this case the following result can be stated
(for all the details, see \cite[Theorem 1.2]{BCF06}).

\begin{theorem}\label{mts}
Let $\M =\R \times\mo$  be a standard stationary spacetime as in
(\ref{e-standardstationary}). Suppose that the smooth Riemannian manifold
$(\mo,\met_R)$ is complete and there exists a point $\bar x\in\mo$
such that coefficient $\beta$ satisfies the at most quadratic
growth condition (\ref{nonlimitato}) and $\delta$ the at most
linear
\begin{equation}\label{bounddelta}
\sqrt{\langle\delta (x),\delta(x)\rangle_R}\ \le\ \lambda_1\
d(x,\bar x) \ +\mu_1 d^{\alpha_1}(x,\bar x) +k_1\quad \hbox{for
all $x \in\mo$}
\end{equation}
for suitable $\lambda_1 \ge 0$, $\mu_1, k_{1}\in \R$, $\alpha_1
\in [0,1)$. Then, $\M$ is geodesically connected\index{geodesic
connectedness!stationary spacetime}.

Furthermore, if $\mo$ is non--contractible in itself, then any two
points can be joined by a sequence of (spacelike) geodesics
$(z_k)_k$ with diverging lengths.
\end{theorem}

Clearly, also in this case the growth hypothesis on $\beta$ cannot
be improved (see the previous Example \ref{counter1}). Morevoer,
if $\delta$ has a more than linear growth, we are able  to find an
example (see \cite[Example 2.7]{BCF06}) whose functional $\J_1$ is
unbounded from below;  so, in general, one does not expect
geodesic connectedness for this case.

\begin{remark}
Also in stationary spacetimes it is possible to give an estimate
of the number of timelike geodesics connecting two fixed points
(see, e.g., \cite[Subsection 3.5]{Ma}); moreover, the completeness assumption
can be weakened (for a quite complete discussion on this subject, see
\cite{Ba00}).
\end{remark}

In Remark \ref{rema4-11}, we pointed 
out that geodesic connectedness is essentially equivalent
in both, the general static and the standard static case; but this was not true in the
stationary case. Moreover, the standard stationary case hides an important drawback:
the same spacetime can split as
(\ref{e-standardstationary}) in very different ways (with very different
$\beta, \delta$) because just one such splitting is not intrinsic
to the spacetime. As a simple and extreme example, Minkowski
spacetime can be written as (\ref{e-standardstationary}) either with an
arbitrary growth of $\delta$ or with an incomplete $\met_R$ (see
\cite[Section 6.2]{CFS06}). More deeply, the bounds
(\ref{nonlimitato}) for $\beta$ and (\ref{bounddelta}) for
$\delta$ do not have a geometric meaning on $\M$, except as
sufficient (but neither necessary nor intrinsic) conditions for
global hyperbolicity (see \cite{Sa97}).
\smallskip

\paragraph{Case {\sl (ii)} (Intrinsic).}  Giannoni and Piccione \cite{GP99} introduced
a different intrinsic approach to the study of geodesic
connectedness in a stationary spacetime whose definition does not
make use of any global splitting.

Recall from (\ref{egeod-prod-Kill}) that if $z:I \to \M$ is a
$C^1$ curve and $K$ is a Killing vector field on $\M$, then it is
\[
\langle z',\frac{D}{ds}K(z)\rangle_L \equiv 0\quad \hbox{on $I$}
\]
and, if $z$ is only absolutely continuous, this holds almost
everywhere in $I$. In particular, if $z$ is a geodesic this
property implies the existence of a constant $C_z \in \R$ such
that
\begin{equation}
\label{constraint} \langle z',K(z)\rangle_L \equiv C_z\quad
\hbox{for all $s \in I$.}
\end{equation}

Thus, if $(\M,\met_L)$ is a stationary spacetime,
(\ref{constraint}) gives a natural constraint to the
action functional $f$ defined in (\ref{action}).
More precisely, for each $p,q\in \M$, we can define
\[
\begin{split}
C^1_K(p,q)\index{$C^1_K(p,q)$} = \{z\in C^1(I,\M) :& z(0) = p,\, z(1) = q,\, \hbox{and}\\
&\exists C_z \in \R\ \hbox{such that $\langle z',K(z)\rangle_L \equiv C_z$}\}  \end{split}
\]
and the following variational principle can be stated
(see \cite[pp. 2]{GP99} joint to some arguments developed in
\cite[Proof of Theorem 3.3]{GP99}).

\begin{theorem} \label{tcpq}
If $z \in C^1_K(p,q)$ is a critical point of $f$ restricted to
$C^1_K(p,q)$, then $z$ is a geodesic connecting $p$ to $q$.
\end{theorem}

Even if functional $f$ is defined in $C^1_K(p,q)$, it
cannot be managed only in this space, as this space is ``too
small'' for problems of convergence. So, the ``natural'' setting
of this variational problem is $\Omega^1(p,q;\M)$ and the corresponding
constraint is
\[
\begin{split}
\Omega^1_K(p,q)\index{$\Omega^1_K(p,q)$} =
\{z\in \Omega^1(p,q;\M) : & \exists C_z \in \R\ \hbox{such that}\\
& \langle z',K(z)\rangle_L= C_z\ \hbox{a.e. on
$[0,1]$}\}.\end{split}
\]
It is quite simple to check that the closure of $C^1_K(p,q)$
with respect to the $H^1$ norm is a subset of $\Omega^1_K(p,q)$; furthermore,
Theorem \ref{tcpq}\index{variational principle!natural constraint}
can be reformulated as follows (see \cite[Theorem 3.3]{GP99}).

\begin{theorem} \label{tcpq1}
If $z \in \Omega^1_K(p,q)$ is a critical point of $f$ restricted to
$\Omega^1_K(p,q)$, then $z$ is a geodesic connecting $p$ and $q$.
\end{theorem}

Thus, in order to find sufficient conditions for the existence of
critical points of $f$ on $\Omega^1_K(p,q)$,
let us introduce the following definitions.

\begin{definition}
Fixed $c\in \R$ the set $\Omega^1_K(p,q)$ is {\em $c$--precompact}\index{precompactness}
if every sequence $(z_k)_k \subset \Omega^1_K(p,q)$ with $f(z_k) \le c$ has a
subsequence which converges weakly in $\Omega^1(p,q;\M)$ (hence, uniformly in $\M$).
Furthermore, the restriction of $f$ to $\Omega^1_K(p,q)$ is
pseudo--coercive\index{pseudo--coerciveness} if $\Omega^1_K(p,q)$ is $c$--precompact for all
$c \ge \inf f(\Omega^1_K(p,q))$.
\end{definition}

Then, the geodesic connectivity between each $p$ and $q$ will be a
consequence of the following theorem (see \cite[Theorem 1.2]{GP99}).

\begin{theorem}\label{intrinsictheo}
If $\Omega^1_K(p,q)$ is not empty and there exists
$c > \inf f(\Omega^1_K(p,q))$ such that
$\Omega^1_K(p,q)$ is $c$--precompact, then
there exists at least one geodesic joining $p$ to $q$ in
$\M$\index{geodesic connectedness!stationary spacetime}.
\end{theorem}

Moreover, if in the stationary spacetime $\M$ there exists a Killing vector field $K$
which is complete, a multiplicity result can be stated (see \cite[Theorem 1.3]{GP99}).

\begin{theorem}\label{intrinsictheo2}
Let $\Omega^1_K(p,q)$ be not empty and assume that $f$ is
pseudo--coercive in $\Omega^1_K(p,q)$. Then, if $K$ is complete
and $\M$ is non--contractible in itself, any two points can be
joined by a sequence of (spacelike) geodesics $(z_k)_k$ with
diverging critical levels $f(z_k)\nearrow +\infty$.
\end{theorem}

Thus, reduced the geodesic connectedness problem to the research
of ``good'' hypotheses which guarantee the application of Theorem
\ref{intrinsictheo}, the explicit example pointed out in
\cite[Appendix A]{GP99} is a standard stationary spacetime $\M =
\R \times\mo$ equipped with metric (\ref{e-standardstationary}) whose
coefficient $\beta$ has to be bounded from above and far away from
zero while $\delta$ must have a sublinear growth, i.e.,
(\ref{bounddelta}) holds with $\lambda_1 = 0$.

Let us point out that the main limitation of Giannoni
and Piccione's results is
that pseudo--coercivity condition is analytical and very
technical. In fact, it can be regarded as a tidy and neat version
of Palais--Smale condition for the stationary ambient.
Furthermore, in general,
the assumption $\Omega^1_K(p,q)$ non--empty must be imposed.
Nevertheless, the possibility of $\Omega^1_K(p,q)=\emptyset$ can be
ruled out if $K$ is complete (compare with \cite[Lemma 5.7]{GP99} and
\cite[Proposition 3.6]{CFS06}).

\begin{proposition}
\label{noempty}
If the timelike Killing vector field $K$
is complete, then we have
 $\Omega^1_K(p,q) \ne\emptyset$ for each $p$, $q \in \M$.
\end{proposition}

Now, we want to translate the technical condition
of pseudo--coercivity in terms of the (Lorentzian) geometry of
the manifold (for a complete proof, see \cite[Section 5]{CFS06}).

\begin{proposition}
\label{coerciveness}
Let $(\M,\met_{L})$ be a stationary spacetime with
a complete timelike Killing vector field $K$ and a complete
smooth spacelike Cauchy hypersurface
$\es$. Then, the restriction of the action function $f$ to $\Omega^1_K(p,q)$ is
pseudo--coercive for any $p$, $q \in \M$.
\end{proposition}

Thus, by means of Theorems \ref{intrinsictheo}, \ref{intrinsictheo2}
and Propositions \ref{noempty}, \ref{coerciveness},
the following results can be stated.

\begin{theorem}\label{globhyp1}
Let $(\M,\met_{L})$ be a stationary spacetime with
a complete timelike Killing vector field $K$. If $\M$ is globally
hyperbolic with a complete (smooth, spacelike) Cauchy hypersurface
$\es$, then it is geodesically connected\index{geodesic connectedness!stationary spacetime}.\\
Furthermore, if $\M$ is non--contractible in itself, then any two
points can be joined by a sequence of (spacelike) geodesics
$(z_k)_k$ such that $f(z_k)\nearrow +\infty$.
\end{theorem}

Let us point out that, if $\M$ has a complete timelike Killing
vector field $K$ and is globally hyperbolic with a complete
spacelike Cau\-chy hypersurface $\es$, from Theorem \ref{t0} it
follows that the spacetime is the product $\R \times\es$, and its
metric can be written as in (\ref{e-standardstationary}) for a certain vector field $\delta$ on
$\es$ and the identifications
\[
K(z) \equiv (0,1) \in \R \times T_x\es \quad \hbox{for all $z =
(t,x) \in \M$ ($t \in \R$, $x \in \es$),}
\]
\[
\langle K(z),K(z)\rangle_L = - \beta(x) \quad \hbox{for all}\quad
z = (t,x)\in \M.
\]
Nevertheless, neither $K$ nor $\es$ are unique and this global splitting
is not canonically associated to a spacetime even under the hypotheses of Theorem \ref{globhyp1}.
Anyway, the results will be independent of the chosen $K, \es$
and no growth hypothesis on $\beta$, $\delta$ is required.

Vice versa, if $\M = \R \times\mo$ is a standard stationary
spacetime with metric (\ref{e-standardstationary}) and the Riemannian metric
$\met_R$ on $\mo$ is complete while $\beta$, respectively
$\delta$, satisfies (\ref{nonlimitato}), respectively
(\ref{bounddelta}), then $\M$ is globally hyperbolic with
$\{0\}\times\mo$ a complete Cauchy hypersurface. Moreover, the
standard timelike Killing vector field $K = \partial_t$ is
complete. Thus, Theorem \ref{mts} follows from Theorem
\ref{globhyp1}.

Let us remark that both the completeness of the Killing vector field
$K$ and the completeness of the Cauchy hypersurface $\es$ are
needed in the proof of Theorem \ref{globhyp1}, as counterexamples
to geodesic connectedness can be found if one of these conditions
is dropped (see \cite[Section 6.3]{CFS06}).
\smallskip

At last, let us summary the results of geodesic connectedness
in a stationary spacetime $(\M,\langle\cdot,\cdot\rangle_L)$:
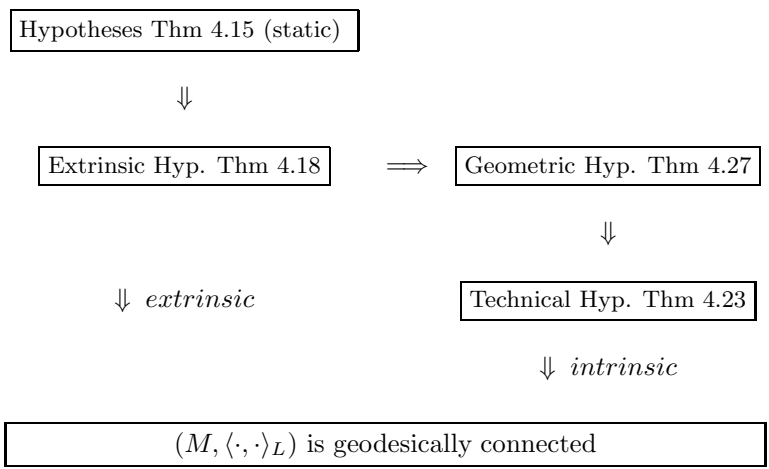
\[
\begin{array}{ccc}
\framebox{\small{Hypotheses Thm \ref{maintheoremstatic} (static)
}}&&\\
\\
 \Downarrow &&\\
\\
\framebox{\small{Extrinsic Hyp. Thm \ref{mts}}}&\Longrightarrow
&\framebox{\small{Geometric Hyp. Thm \ref{globhyp1}}}\\
\\
&&\Downarrow\\
\\
\Downarrow\; extrinsic&&\framebox{\small{Technical Hyp. Thm \ref{intrinsictheo}}}\\
\\
&&\Downarrow\;intrinsic
\end{array}
\]
\[
\framebox{\qquad\qquad\qquad$(\M,\langle\cdot,\cdot\rangle_L)$
is geodesically connected\qquad\qquad\qquad}
\]
And, as widely discussed in \cite{CFS06}, the intrinsic geometric
hypotheses on the stationary spacetime  in Theorem \ref{globhyp1} are
essentially equivalent to the technical analitycal hypotheses on
the space of curves  in Thm. \ref{intrinsictheo}.

\subsection{Time dependent metrics and saddle critical points}
\label{secsplitting}

Now, we want to take care of models of Lorentzian manifolds
whose metric is time--dependent so it is necessary
to study directly the corresponding strongly indefinite
action functional\index{functional!action} $f$ defined in (\ref{action}).
Such an approach is based essentially on a Galerkin
finite--dimensional approximation of the ``real part'' $W^1(t_p,t_q)$
and yields estimates for the connectedness of very general manifolds, as the
orthogonal splitting ones. In this setting,
the existence of at least one connecting geodesic
is obtained by means of Rabinowitz's Saddle Point
Theorem while multiplicity results follow from the Relative Category
Theory.

The basic ideas of this method were introduced in \cite{BF94}, in
order to study the geodesic connectedness in a standard stationary
spacetime, then they were applied to the splitting case in
\cite{BFM94,GiaMa95,Ma,PiSa96}. Anyway, some details were
completely analyzed later on (see \cite{CGM99,CMS00}) and, even if
such results follow from \cite[Theorems 1.5 and 1.10]{CMS00},
here, they are explained for the first time in a self--contained
presentation.

First of all, let us introduce the main tools of the abstract
theorems we need for dealing with unbounded functionals.

In order to obtain the existence of at least one solution,
let us state the following slight variant of the classical
Saddle Point Theorem\index{Theorem!Saddle Point}
(cf. \cite{BFM94,Ra} or also \cite[Theorem 3.2]{CMS00}).

\begin{theorem}
\label{saddlepoint} Let $\Omega$ be a complete Riemannian manifold
and $H$ a separable Hilbert space. Fix $H_0$ as a linear subspace
of $H$ and $j \in H$. Moreover, let $(a_l)_{l\in\N}$ be an
orthonormal basis of $H_0$. Set $W = H_0 + j$ and $Z = W
\times\Omega$.
Let $f : Z \to \R$ be a $C^1$ functional.\\
For any $m \in \N$,  let $W_m = span\{a_1,a_2,\dots,a_m\} + j$,
$Z_m = W_m \times\Omega$ and $f_m = f|_{Z_m}$.\\
Fixed $y \in \Omega$, for any $R > 0$
consider the sets:
\[
\begin{split}
S &=\; \{(j,x) \in Z : x \in \Omega\} =
\{j\} \times\Omega\ ,\\
Q(R) &=\; \{(t,y) \in Z : \|t - j \|_H \le R\}\ ,
\end{split}
\]
where $\|\cdot\|_H$ is the norm of the Hilbert space $H$.\\
Assume that $f_m$ satisfies condition $(PS)$ for any $m \in \N$
and there exists $R > 0$ such that
\[
\sup f(Q(R)) < +\infty\ ,\qquad
\sup f(\partial Q(R)) < \inf f(S)\ .
\]
Then, for any $m \in \N$, $f_m$ has a critical level $c_m \in
[\inf f(S),\sup f(Q(R))]$, where
\[
c_m = \inf_{h \in \Gamma_m} \sup_{z \in Q_m(R)} f_m (h (z))\ ,
\]
\[
\Gamma_m = \{ h \in C(Z_m,Z_m) : h(z) = z\; \mbox{for all}\;
            z \in \partial Q_m(R)\}
\]
and
\[
Q_m(R) =\ \{(t,y) \in Z_m : \|t - j\|_H \le R\}\ .
\]
\end{theorem}

Now, we need introducing the Relative Category Theory\index{Theory!Relative Category},
which generalizes the Ljusternik--Schni\-rel\-man one.
To this aim, let us give the notion of relative category and some of its main properties
(see, e.g., \cite{F,FW90,Sz90}).

\begin{definition}
\label{categoriarelativa} Let $Y$ and $A$ be closed subsets of a
topological space $Z$. The {\em category of $A$ in $Z$ relative to
$Y$}\index{category!relative}, briefly
$\cat_{Z,Y}(A)$\index{$\cat_{Z,Y}(A)$}, is the least integer $k$
such that there exist $k + 1$ closed subsets of $Z$, $A_0,A_1,
\dots, A_k$, $A = A_0 \cup A_1 \cup \dots \cup A_k$, and $k + 1$
functions, $h_l \in C([0,1] \times A_l,Z)$, $l \in \{0,1,\dots,
k\}$, such that
\begin{enumerate}
\item[{\sl (a)}] $h_l(0,z) = z$ for $z \in A_l$, $0 \le l \le k$;
\item[{\sl (b)}] $h_0(1,z) \in Y$ for $z \in A_0$, and
$h_0(\sigma,y) \in Y$ for all $y \in A_0 \cap Y$, $\sigma \in
[0,1]$; \item[{\sl (c)}] $h_l(1,z) = z_l$ for $z \in A_l$ and some
$z_l \in Z$, $1 \le l \le k$.
\end{enumerate}
If a finite number of such sets does not exist, we set
$\cat_{Z,Y}(A) = +\infty$.
\end{definition}

We have that $\cat_Z(A) = \cat_{Z,\emptyset}(A)$
is the classical Ljusternik--Schnirelman category of $A$ in $Z$.

\begin{proposition}
\label{proprieta}
Let $A$, $B$, $Y$ be closed subsets of a
topological space $Z$.
\begin{enumerate}
\item[{\sl (i)}] If $A \subset B$ then $\cat_{Z,Y}(A) \le \cat_{Z,Y}(B)$;
\item[{\sl (ii)}] $\cat_{Z,Y}(A \cup B) \le \cat_{Z,Y}(A) + \cat_Z(B)$;
\item[{\sl (iii)}] if there exists
$h \in C([0,1] \times A,Z)$ such that $h(\sigma,y) = y$ for $y \in A \cap Y$
and $\sigma \in [0,1]$, then
$\cat_{Z,Y}(A) \le \cat_{Z,Y}(B)$ where $B = \overline{h(1,A)}$.
\end{enumerate}
\end{proposition}

\begin{remark}
Let $Z$ be a topological space and $Y$ a closed subset of $Z$.
Then Proposition~\ref{proprieta}{\sl (ii)} \ implies that
the relative category and the classical Ljusternik--Schnirelman
category are connected by the inequality
\[
\cat_{Z,Y} (A) \le \cat_Z (A) \quad \mbox{for any closed set}\ A \subset Z\ .
\]
\end{remark}

It is easy to see that Definition~\ref{categoriarelativa}\
implies the following proposition.

\begin{proposition}
\label{deformazione}
Let $Z$ be a topological space and $C$, $\Lambda$ be two subsets
of $Z$ such that $C$ is a closed strong deformation retract of
$Z \backslash \Lambda$, i.e., there exists a continuous map
${\cal R} : [0,1] \times (Z \backslash \Lambda) \to Z$ such that
\[
\left\{
\begin{array}{ll}
{\cal R}(0,z) = z &\mbox{for all $z \in Z \backslash \Lambda$,}\\
{\cal R}(1,z) \in C &\mbox{for all $z \in Z \backslash \Lambda$,} \\
{\cal R}(\sigma,z) = z &\mbox{for all $z \in C$, $\sigma \in [0,1]$.}
\end{array}
\right.
\]
Then $\cat_{Z,C} (Z \backslash \Lambda) = 0$.
\end{proposition}

A further property of the relative category can be stated
(for the proof, see \cite[Proposition 2.2]{CGM99}).

\begin{proposition}
\label{retratto}
Let $Y$, $Z'$, $Y'$ be closed subsets of a
topological space $Z$ such that $Y' \subset Z'$. Suppose that there exist
a retraction $r : Z \to Z'$, i.e., a
continuous map such that
$r(z) = z$ for all $z \in Z'$,
and a homeomorphism $\Phi : Z \to Z$ such that
$\Phi(Y') \subset Y$ and $r \circ \Phi^{- 1}(Y) \subset Y'$.
Then, if $A'$ is a closed subset of $Z'$, it results that
\[
\cat_{Z,Y}(\Phi(A')) \ge \cat_{Z',Y'}(A')\ .
\]
\end{proposition}

Now, we can state a multiplicity
result for critical levels of a strongly indefinite func\-tional
(for the proof, see \cite[Theorem 1.4]{CGM99}).

\begin{theorem}
\label{teomolteplicita}
Let $Z$ be a $C^2$ complete Riemannian manifold
modelled on a Hilbert space
and $f : Z \to \R$ a $C^1$ functional.
Let us assume that there exist two subsets $\Lambda$ and
$C$ of $Z$ such that
$C$ is a closed strong deformation retract of $Z \backslash \Lambda$,
\[
\inf_{z \in \Lambda} f(z) > \sup_{z \in C} f(z)
\quad\hbox{and}\quad
\cat_{Z,C} (Z) > 0\ .
\]
If $f$ satisfies condition $(PS)$, then it has at least
$\cat_{Z,C}(Z)$ critical points in $Z$ whose critical levels are
greater or equal than $\inf f(\Lambda)$. Moreover, if
$\cat_{Z,C}(Z) = +\infty$, there exists a sequence
$(z_k)_{k}$ of critical points of $f$ such that
\[
\lim_{k \to +\infty} f(z_k) = \sup_{z \in Z} f(z)\ .
\]
\end{theorem}

\begin{remark}
\label{remanoveastratto} In Theorem~\ref{teomolteplicita}\ the
critical levels $c_k$ are characterized as follows:
\[
c_k = \inf_{B \in F_k} \sup_{z \in B} f(z)\quad \mbox{for any $k
\in \N$, $1 \le k \le \cat_{Z,C}(Z)$,}
\]
where
\[
F_k = \{ B \subset Z :\ \mbox{$B$ closed,}\; \cat_{Z,C}(B) \ge
k\}\ .
\]
\end{remark}

The aim of this section is to investigate the existence of geodesics joining two fixed
points in a Lorentzian manifold which splits in a ``good'' product
but it has a metric depending on both its coordinates.
The simplest model of this type can be defined as follows.

\begin{definition}
\label{varietasplitting} A Lorentzian manifold $(\M,\met_L)$ is an
{\em orthogonal splitting spacetime}\index{spacetime!orthogonal
splitting} if it is $\M = \R\times\mo$ and
\begin{equation}
\langle\zeta,\zeta'\rangle_L = - \beta(z) \tau \tau'+
\langle\alpha(z)\xi,\xi'\rangle_R \ , \label{metrica3}
\end{equation}
for any $z = (t,x) \in \M$ ($t \in \R$, $x\in \mo$), and $\zeta =
(\tau,\xi)$, $\zeta' = (\tau',\xi') \in T_z\M \equiv
 \R \times T_x\mo$,
where $(\mo,\met_R)$ is a finite dimensional connected Riemannian
manifold, $\alpha(z)$ is a smooth symmetric linear strictly
positive operator from $T_x\mo$ into itself and $\beta : \M \to
\R$ is a smooth and strictly positive scalar field.
\end{definition}
 Clearly, if both $\alpha$ and $\beta$ do not
depend on $t$, Definition \ref{varietasplitting} reduces to the
static metric (\ref{metrica}) and the geodesic connectedness
result is Theorem \ref{maintheoremstatic}, already proved with a
quadratic growth assumption on $\beta$ (once that
$\langle\alpha(\cdot) \cdot, \cdot\rangle_R$ is complete). 
But, in general, when
coefficients depend on time component $t$, stronger assumptions,
in particular a ``good'' control on the partial derivatives
$\alpha_t$, $\beta_t$, is needed. More precisely, the following result
will be proved below.

\begin{theorem}
\label{maintheoremsplitting} Let $\M = \R \times \mo$ be an
orthogonal splitting manifold such that $(\mo,\met_R)$ is a
complete Riemannian manifold. Assume that there exist some
constants $\lambda,\ \nu,\ N,\ K > 0$ such that the coefficients
$\alpha$, $\beta$ of its metric (\ref{metrica3}) satisfy the
following hypotheses:
\begin{equation}
\label{alfa} \lambda\ \langle\xi,\xi\rangle_R\ \le\
\langle\alpha(z)\xi,\xi\rangle_R\ ,
\end{equation}
\begin{equation}
\label{beta} \nu \le \beta(z) \le N\ ,
\end{equation}
\begin{equation}
\label{derivate} |\beta_t(z)| \le K\ ,\quad
|\langle\alpha_t(z)\xi,\xi\rangle_R| \le K
\langle\xi,\xi\rangle_R\ ,
\end{equation}
for all $z =(t,x) \in \M$, $\xi \in T_x\mo$. Furthermore, assume that
\begin{eqnarray}
\label{limsup} \limsup_{t \to +\infty}\
(\sup\{\langle\alpha_t(t,x)\xi,\xi\rangle_R :\ x \in \mo,\; \xi
\in
T_x\mo,\ \langle\xi,\xi\rangle_R = 1\}) \ \le\ 0,\\
\label{liminf} \liminf_{t \to -\infty} \
(\inf\{\langle\alpha_t(t,x)\xi,\xi\rangle_R :\ x \in \mo,\; \xi
\in T_x\mo,\ \langle\xi,\xi\rangle_R = 1\}) \ \ge\ 0.
\end{eqnarray}
Then, fixed $p=(t_p,x_p)$, $q=(t_q,x_q) \in \M$, there exists at least one
geodesic joining $p$ to $q$.\index{geodesic connectedness!ortogonal splitting spacetime}\\
Moreover, if $\mo$ is not contractible in itself and there exists
$\lambda_1 > 0$ such that
\begin{equation}
\label{alfa2} \langle\alpha(0,x)\xi,\xi\rangle_R\ \le\ \lambda_1\
\langle\xi,\xi\rangle_R\quad \hbox{for all $x \in \mo$, $\xi \in
T_x\mo$,}
\end{equation}
infinitely many such (spacelike) geodesics $(z_k)_k$ exist so that
$f(z_k)\nearrow +\infty$.
\end{theorem}

The first problem in applying directly the abstract theorems to
action functional $f$ in (\ref{action}), is that $f$ may not
satisfy condition $(PS)$ on the whole manifold $\Omega^1(p,q;\M)$.
Furthermore, a finite dimensional decomposition of $W^1(t_p,t_q)$
needs in order to apply not only Theorem \ref{saddlepoint},
obvious for its statement, but also Theorem \ref{teomolteplicita},
as we have to be sure that we deal with manifolds with
non--trivial relative category.

A way to overcome the lack of condition $(PS)$ is to use a
penalization argument on $f$ as follows.

For any $\eps > 0$, let $\psi_\eps : \R_+ \to \R_+$ be the
``cut--function'' defined as
\begin{equation}
\label{penalita} \psi_\eps (s) = \left\{
\begin{array}{ll}
0  &\mbox{if $\displaystyle{0 \le s \le \frac{1}{\eps}}$} \\
                      &\\
\displaystyle\sum_{n=3}^{+\infty}\frac{\left(s -
\frac{1}{\eps}\right)^n}{n!}
   &\mbox{if $\displaystyle{s > \frac{1}{\eps}}$.}
\end{array}
\right.
\end{equation}
Obviously, $\psi_\eps$ is a $C^2$ map on $\R_+$, it is increasing
and there exist two positive constants, $a$, $b$, such that
\begin{equation}\label{primaprop}
\psi'_\eps(s) \ge \psi_\eps(s) \ge a s - b\ , \;\ s\ \psi'_\eps(s)
\ge \psi_\eps(s) \quad \mbox{for all} \ s \in \R_+ .
\end{equation}
Thus, for any $\eps > 0$ we consider the perturbed functional
\begin{equation}
\label{funzionalepenalizzato}
\begin{array}{rl}
f_\eps (z) &= f(z) - \psi_\eps(\| t'\|^2)\\
&\displaystyle =\ \int_0^1 (\langle\alpha(t,x) x',
x'\rangle_R -\ \beta(t,x)\ (t')^2)\ ds - \psi_\eps(\|t'\|^2)
\end{array}
\end{equation}
for any $z = (t,x) \in \Omega^1(p,q;\M)$, where $\| t'\|^2 =
\int_0^1| t'|^2 ds$. Clearly, it is
\[
f_\eps(z) \le f(z)\quad \mbox{for all $z \in \Omega^1(p,q;\M)$,}
\]
furthermore, the Fr\'echet differential of $f_\eps$ at $z =
(t,x)\in \Omega^1(p,q;\M)$ is given by
\[
f'_\eps (z)[\zeta]\ =\ 2\ \int_0^1 \langle z',\frac{D
\zeta}{d s}\rangle_L\ ds -\ 2\ \psi'_\eps(\| t'\|^2)\
\int_0^1  t'\ \tau'\ ds
\]
for all $\zeta = (\tau,\xi) \in T_z \Omega^1(p,q;\M)$.

Taking any $z_\eps=(t_\eps,x_\eps)$ critical point of $f_\eps$,
some usual bootstrap arguments allow one to prove that it is
smooth and satisfies the equation
\[
\frac{D z'}{d s} - \psi'_\eps(\| t'\|^2)\ (t'',0) = 0;
\]
hence, a constant $E_\eps(z_\eps) \in \R$ exists so that
\begin{equation}
\label{energiapenalizzata} \langle z'_\eps(s),
z'_\eps(s)\rangle_L -\ \psi'_\eps(\| t'_\eps\|^2)\
(t'_\eps)^2(s) \equiv E_\eps(z_\eps)\quad \hbox{for all $s\in I$.}
\end{equation}
Thus, by integrating (\ref{energiapenalizzata}), according to
(\ref{primaprop}) it follows that
\begin{equation}
\label{maggioroenergia} E_\eps(z_\eps) =\ f(z_\eps) -\
\psi'_\eps(\|t'_\eps\|^2)\ \|t'_\eps\|^2 \le
f_\eps(z_\eps)\ .
\end{equation}

Obviously, if $L^* > 0$ is a fixed constant, by (\ref{penalita})
and (\ref{funzionalepenalizzato}) it follows that a critical point
$z_\eps = (t_\eps,x_\eps)$ of the penalized functional $f_\eps$ is
also a critical point of the action functional $f$ if it is such
that
\begin{equation}\label{supp}
\| t'_\eps\|^2 \le L^*\qquad\hbox{and}\qquad \eps \leq \frac
1{L^*}.
\end{equation}

\begin{remark}
\label{limitazioni} A direct consequence of the given hypotheses
is that on sublevels of $f$ the norm of the $t$ component
``controls'' the $x$ component. In fact, taken $z = (t,x)$ such
that $f_\eps(z) \le L$, $L \in \R$, then (\ref{alfa}),
(\ref{beta}) and (\ref{funzionalepenalizzato}) imply that
\[
\lambda\ \int_0^1\langle x', x'\rangle_R\ ds \le
\int_0^1\langle\alpha(z) x', x'\rangle_R\ ds \le L + N
\| t'\|^2 + \psi_\eps(\| t'\|^2)\ .
\]
\end{remark}

\begin{lemma}\label{Palais0}
If $\mo$ is complete and (\ref{alfa}), (\ref{beta}),
(\ref{derivate}) hold, then taken any $\eps > 0$ the corresponding
functional $f_\eps$ satisfies condition $(PS)$ on
$\Omega^1(p,q;\M)$.
\end{lemma}

\begin{proof} Up to small changes, the main ideas of this proof can be found
in \cite[Proposition 4.4]{CMS00}. Here, for completeness, we
outline the main steps.\\
Let $L \in \R$ and $(z_k)_k \subset \Omega^1(p,q;\M)$ be such that
\[
f_\eps(z_k) \le L \;\hbox{for all $k \in \N$,}\qquad f'_\eps(z_k)
\to 0\;\hbox{if $k \to +\infty$.}
\]
Firstly, suitable computations making use of hypotheses
(\ref{alfa}), (\ref{beta}), (\ref{derivate}) and properties on
$\psi_\eps$ in (\ref{primaprop}), allow one to prove that $(\|t'_k\|)_k$
is bounded; then, by Remark \ref{limitazioni}, also
$(\| x'_k\|)_k$ is bounded.
Hence, $(z_k)_k$ is a bounded sequence in $\Omega^1(p,q;\M)$.\\
Thus, since $\mo$ is complete, by means of Lemma \ref{splitting}
and careful limit estimates, there exists $z \in \Omega^1(p,q;\M)$
such that $z_k \to z$ strongly in $\Omega^1(p,q;\M)$ (up to
subsequences).
\end{proof}

Now, once we are able to find critical points of some $f_\eps$,
the problem is to ``reduce'' them to critical points of the
``original'' functional $f$ by means of some a priori estimates as
(\ref{supp}).

\begin{lemma}\label{noeps}
If the hypotheses (\ref{alfa}) -- (\ref{liminf}) hold, then for
each fixed $L > 0$ there exists $\eps_0 > 0$ such that, if $z_\eps
= (t_\eps,x_\eps)\in \Omega^1(p,q;\M)$ is such that
\begin{equation}\label{po}
f'_\eps(z_\eps) = 0\qquad\hbox{and}\qquad f_\eps(z_\eps) \le L,
\end{equation}
with $\eps \le\eps_0$, then it is $\|t'_\eps\|^2 \le \frac
1{\eps_0}$; hence, $\psi_\eps(\|t'_\eps\|^2) = 0$.
\end{lemma}

\begin{proof}
Also for this proof, the main ideas are in \cite[Propositions 5.2
and 5.3]{CMS00} but, anyway, here we
outline the main steps.\\
Firstly, we want to prove that there exists $L^* > 0$, independent
of $\eps$, such that if (\ref{po}) holds then $|t_\eps|_\infty \le
L^*$. To this aim, taken any $\mu > 0$ and by means of hypotheses
(\ref{limsup}), (\ref{liminf}), there exists $t_\mu >
\max\{|t_p|,|t_q|\}$ such that, by using (\ref{alfa}), for all $x
\in \mo$, $\xi \in T_x\mo$ it is
\begin{eqnarray*}
&\langle\alpha_t(t,x)\xi,\xi\rangle_R \le \lambda \mu
\langle\xi,\xi\rangle_R \le \mu \langle\alpha(t,x)\xi,\xi\rangle_R
\quad &\hbox{if $t \ge t_\mu$,} \\
&\langle\alpha_t(t,x)\xi,\xi\rangle_R \ge -\lambda \mu
\langle\xi,\xi\rangle_R \ge -\mu
\langle\alpha(t,x)\xi,\xi\rangle_R \quad &\hbox{if $t \le
-t_\mu$.}
\end{eqnarray*}
Assume that $z_\eps = (t_\eps,x_\eps)$ exists such that (\ref{po})
holds and $\|t_\eps\|_\infty > t_\mu$ (the latter, for example, as 
the essential supremum satisfies
$\displaystyle{\rm ess}\sup_{s\in I}t_\eps(s) > t_\mu$). Hence, an
interval $[a,b] \subset ]0,1[$ exists such that $t_\eps(s) >
t_\mu$ for all $s \in ]a,b[$, $t_\eps(a) = t_\eps(b) = t_\mu$.
Fixed $\omega > 0$, let us define function $\tau : I \to \R$ such
that
\[
\tau(s) = \left\{
\begin{array}{ll}
\sinh(\omega (t_\eps(s) - t_\mu)) &\mbox{if $s \in [a,b]$}\\
0 &\mbox{if $s \in I \backslash [a,b]$.}
\end{array}
\right.
\]
Clearly, it is $\tau \in H^1_0(I,\R)$; thus, it is $(\tau,0)\in
T_{z_{\eps}}\Omega^1(p,q;\M)$ and $f'_\eps(z_\eps)[(\tau,0)] = 0$.
By means of careful estimates, the main properties of function $y
= \sinh\sigma$, (\ref{alfa})--(\ref{derivate}) and
(\ref{maggioroenergia}) with (\ref{po}), suitable choices of $\mu$
and $\omega$ allow one to prove that there exists $\theta > 0$
(depending on $\mu$ but independent of $\eps$) such that $\omega
\sqrt\theta < 2$ and
\[
\int_a^b (t'_\eps)^2\ \cosh(\omega (t_\eps - t_\mu))\ ds\ \le\
\theta\ \int_a^b \cosh(\omega (t_\eps - t_\mu))\ ds\ .
\]
Hence, by applying \cite[Lemma 3.4]{GiaMa95}, a constant $L^*>0$
exists, independent of $\eps$, such that it has to be
$\displaystyle{\rm ess}\sup_{s\in I}t_\eps(s) \le L^*$.
Similarly, one can argue if it is $\displaystyle{\rm
ess}\sup_{s\in I}(-t_\eps(s)) > t_\mu$.
Hence, in any case it has to be $|t_\eps|_\infty \le L^*$.\\
Now, we have to prove that also $\| t'_\eps\|$ is bounded
independently of $\eps$ if (\ref{po}) holds. For this proof, let
us start from $f'_\eps(z_\eps)[(\tau,0)] = 0$ with
\[
\tau(s) = \sinh(\omega (t_\eps(s) - j^*(s)), \qquad s \in I,
\]
with arbitrary $\omega >0$ (later on, some restrictions on it will be added).
By making use of ``good'' tricks
as Young inequality
\[
\int_0^1 | t'| \cosh(\omega(t_\eps - j^*)) ds \le
\frac{\delta}{2}
 \int_0^1  (t')^2 \cosh(\omega(t_\eps - j^*)) ds
 + \frac1{2 \delta} \int_0^1 \cosh(\omega(t_\eps - j^*)) ds,
\]
for any $\delta > 0$, and of estimates due to the boundedness of
$|t_\eps|_\infty$ such as
\[
1\ \le \ \cosh(\omega(t_\eps(s) - j^*(s)))\ \le \gamma
\quad\mbox{for all $s \in I$,}
\]
($\gamma$ independent of $\eps$), and by taking into
account the given hypotheses (\ref{alfa}) -- (\ref{derivate}) and
by doing straightforward calculations, a suitable choice of a
large enough $\omega$ and a small enough $\delta$ (so to have
strictly positive all the coefficients in the inequality below),
allow one to prove that
\[
\begin{split}
&( 2 \omega \nu\ -\ \frac{K N}{\lambda}\ - K - \omega |\Delta_t| N
\delta)\ \| t'\|^2 + (2 \omega - \frac{K}{\lambda})
\psi'_\eps(\| t'\|^2) \|t'\|^2\\
&\le \ \gamma\ (\frac{K L}{\lambda} + 2 \omega \Delta_t^2
\psi'_\eps(\| t'\|^2) + N |\Delta_t| \frac{\omega}{\delta}) .
\end{split}
\]
Hence, the proof follows from (\ref{primaprop}).
\end{proof}

Thus, from now on in this section, let us assume that the hypotheses of Theorem
\ref{maintheoremsplitting} are satisfied and, as a further step,
let us introduce a finite dimensional decomposition of
$W^1(t_p,t_q)$.

So, let us consider the orthonormal basis $\{\sin(l \pi s)\}_{l \in \N}$
of $H^1_0(I,\R)$ and for any  $m \in \N$ define
\[
W_m = H_m + j^*\ , \quad\hbox{with}\;
H_m = span \left\{ \sin(l \pi s) : \; l \in\{1, 2,\dots, m\} \right\}\ ,
\]
and
\[
Z_m = W_m \times \Omega^1(x_p,x_q;\mo),\quad
f_{\eps,m} = f_\eps|_{Z_m} \; \hbox{(for any $\eps > 0$).}
\]

Clearly, by arguing as in Proposition \ref{Palais0}, the following lemma
can be stated.

\begin{lemma}
\label{PSlemma}
For any $\eps > 0$ and $m \in \N$, the
functional $f_{\eps,m}$ satisfies condition $(PS)$.
\end{lemma}

Obviously, fixed $\eps > 0$, once we are able to find critical points of
$f_{\eps,m}$ on $Z_m$, it is necessary to come back to the
critical points of $f_{\eps}$ on the whole manifold $\Omega^1(p,q;\M)$.
Thus, to this aim, we need the following proposition.

\begin{lemma}
\label{galerkin}
Assume that for all $m \in \N$ there exists $z_m \in Z_m$ critical point of
$f_{\eps,m}$ on $Z_m$ such that
\[
c_1 \le f_{\eps,m} (z_m) \le c_2
\]
for two given constants $c_1$ and $c_2$ (independent of $m$).
Then, sequence $(z_m)_{m \in\N}$ converges, up to a subsequence, to a
critical point $z \in \Omega^1(p,q;\M)$ of $f_\eps$ such that
\begin{equation}\label{ppp}
c_1 \le f_\eps (z) \le c_2\ .
\end{equation}
\end{lemma}

\begin{proof}
By arguing as in the proof of Proposition \ref{Palais0} it follows that
$(z_m)_{m\in\N}$ is a bounded sequence of $\Omega^1(p,q;\M)$; hence, up to a subsequence,
$z \in \Omega^1(p,q;\M)$ exists so that $z_m \wk z$ weakly in $\Omega^1(p,q;\M)$.\\
Set $z_m = (t_m,x_m)$, $z = (t,x)$, from Lemma \ref{splitting} it follows that
two sequences $\xi_m \in T_{x_m}\Omega^1(x_p,x_q;\mo)$
and $\nu_m \in H^1(I,\R^N)$ exist such that (\ref{sequences}) holds.
Thus, denoted by $P_m: H^1_0(I,\R) \to H_m$ the orthogonal projection between such two
spaces and $\tau_m = P_m(t - j^*) - (t_m - j^*)\in H_m \equiv T_{t_m}W_m$,
we have $f'_{\eps,m} (z_m)[(\tau_m,\xi_m)] = 0$ for all $m \in \N$.
Whence, the same arguments of the second part of the proof
of Proposition \ref{Palais0} show that
$z_m \to z$ strongly in $\Omega^1(p,q;\M)$ while
from \cite[Lemma 3.4]{BF94} it follows that
$f'_\eps(z) = 0$ and (\ref{ppp}) holds.
\end{proof}

Now, we want to apply Theorem \ref{saddlepoint} to functional $f_\eps$ in $\Omega^1(p,q;\M)$.\\
Fixing $y \in \Omega^1(x_p,x_q;\mo)$ 
and $C^1$, let us define the following sets:
\begin{eqnarray}
S &=& \{j^*\} \times \Omega^1(x_p,x_q;\mo) ,\label{esse}\\
Q(R) &=& \big\{(t,y) \in \Omega^1(p,q;\M) :\ \|t' - \Delta_t\| \le R\big\} \quad
\hbox{(for any $R > 0$),}\nonumber
\end{eqnarray}
with $\displaystyle\frac{dj^*}{ds}(s) \equiv \Delta_t$. Obviously, taken $\eps_0 \in\ ]0,1]$ so that
$\Delta_t^2 < \frac{1}{\eps_0}$, from (\ref{penalita}) it follows
\[
\psi_\eps(\|\Delta_t\|^2) = 0\quad \hbox{for all $\eps \in\ ]0,\eps_0]$;}
\]
hence, (\ref{alfa}) and (\ref{beta}) imply
\begin{equation}
\label{provacinque}
\inf f_\eps(S)\ \ge\ - N \Delta_t^2
\quad \mbox{for all $\eps \in\ ]0,\eps_0]$.}
\end{equation}
On the other hand, by means of suitable computations (for more details, see
the proof of \cite[Theorem 1.4]{CMS00}), two positive constants $\bar k_1$ and $\bar k_2$
exist, independent of $\eps$ and $\Delta_t$, such that for all $z \in Q(R)$ it is
\[
f_\eps(z)\ \le\ f(z)\ \le\
\bar k_1 +\bar k_2\ \| t'\|\ -\ \nu\ \| t'\|^2
\le \bar k_1 + \bar k_2\ \sqrt{R^2 + \Delta_t^2}
\ -\ \nu\ (R^2 + \Delta_t^2).
\]
So, not only it is
\[
\sup f_\eps(Q(R)) < +\infty\quad\hbox{for all $R > 0$}
\]
but also there exists $R$ large enough so that
\[
\sup f_\eps(\partial Q(R)) < \inf f_\eps(S).
\]
Thus, the existence of a critical point for the action functional
$f$ in $\Omega^1(p,q;\M)$, i.e., a geodesic joining $p$ to $q$ in $\M$,
follows from Theorem \ref{saddlepoint} and Lemmas
\ref{noeps}, \ref{PSlemma}, \ref{galerkin}.
\smallskip

Now, in order to apply Theorem \ref{teomolteplicita}, it
is necessary exploiting the topological properties of the space of
curves $\Omega^1(p,q;\M)$. To this aim, the following result is
basic and extends Proposition \ref{FadellHusseini} to the relative
category of a product set involving $\Omega^1(x_p,x_q;\mo)$ (for
more details, see \cite[Corollary 3.2]{FH94}).

\begin{proposition}
\label{FadellHusseini2} Let $\mo$ be a simply connected and
non--contractible smooth manifold, $x_p$ and $x_q$ two points of
$\mo$ and $D^m$ the unit disk in $\R^m$ with boundary $S^m$. Then,
for any $k \in \N$, there exists a compact subset $V_{m,k}$ of
$D^m \times\Omega^1(x_p,x_q;\mo)$ such that
\[
\cat_{D^m \times\Omega^1(x_p,x_q;\mo), S^m \times\Omega^1(x_p,x_q;\mo)} (V_{m,k}) \ge k .
\]
\end{proposition}

\begin{remark}
\label{indippartespaziale} Let $(V_{m,k})_{m,k}$ be the family of
compact subsets of the product manifold
$D^m\times\Omega^1(x_p,x_q;\mo)$ obtained in Proposition
\ref{FadellHusseini2}. Fixing $k \in \N$, by means of
\cite[Proposition 2.12, Theorems 2.14 and 3.1]{FH94} it follows
that, for each $m \in \N$ the set $V_{m,k}$ has the same
projection on $\Omega^1(x_p,x_q;\mo)$, i.e., there exists $V^x_k$
compact subset of $\Omega^1(x_p,x_q;\mo)$, independent on $m$, and
$V^t_{m,k} \subset D^m$ such that $V_{m,k} = V^t_{m,k} \times
V^x_k$.
\end{remark}

\begin{lemma}
\label{lemmauno}
There exists a continuous map
$\varrho : \R_+ \to \R_+$ such that
\[
z = (t,x) \in \Omega^1(p,q;\M)\ ,\; \| t' - \Delta_t\| = \varrho(\| x'\|) \quad
\then \quad f(z) \le \left\{
\begin{array}{ll}
- 2 N \Delta_t^2 &\hbox{if $\Delta_t \ne 0$,}\\
- 1&\hbox{if $\Delta_t = 0$,}\end{array}\right.
\]
where $N$ is defined as in (\ref{beta}).
\end{lemma}

\begin{proof}
Taken $z = (t,x) \in \Omega^1(p,q;\M)$, then (\ref{beta}), (\ref{derivate})
and assumption (\ref{alfa2}) allow one to prove that
\[
f(z)\ \le\ (\lambda_1 + K |\Delta_t|)\ \| x'\|^2 +
K \| x'\|^2 \| t' - \Delta_t\| - \nu\ \| t' - \Delta_t\|^2\ .
\]
Hence, it is enough to define
\[
\varrho (r)\ =\ \frac{K r^2}{2 \nu}\ +\
\sqrt{\frac{K^2 r^4}{4 \nu^2}\ +\
\frac{(\lambda_1 + K |\Delta_t|) r^2 + 2 N \Delta_t^2}{\nu}}.
\]
\end{proof}

Now, fixing $m \in \N$,
let us define the ``cylinder''
\[
C_m = \{z = (t,x) \in Z_m : \| t' - \Delta_t\| = \varrho(\| x'\|)\}\ .
\]
For each $\eps > 0$, it is not difficult to show that (\ref{penalita}), (\ref{provacinque})
and Lemma \ref{lemmauno} imply
\[
\sup f_{\eps,m}(C_m)\  <\ \inf f_{\eps,m}(S)
\]
being $S$ as in (\ref{esse}).
Furthermore, the subset $C_m$ is a strong deformation retract of $Z_m \backslash S$;
hence, from Proposition \ref{deformazione} it follows
$\cat_{Z_m,C_m} (Z_m \backslash S) = 0$.

\begin{lemma}
\label{lemmaquattro}
Let $\mo$ be 1--connected.
For any $m$, $k \in \N$ there exists a compact subset $K_{m,k}$ of
$Z_m$, whose projection on $\Omega^1(x_p,x_q;\mo)$ is independent of $m$, such that
\[
\cat_{Z_m,C_m} (K_{m,k}) \ge k .
\]
\end{lemma}

\begin{proof} Let us consider the following sets:
\[
D_m = \{t \in W_m : \| t' - \Delta_t\| \le 1\} ,\quad
\tilde D_m = D_m \times \Omega^1(x_p,x_q;\mo) ,
\]
\[
\Sigma_m = \partial D_m = \{t \in W_m : \| t' - \Delta_t\| = 1\} ,
\quad \tilde\Sigma_m = \Sigma_m \times\Omega^1(x_p,x_q;\mo).
\]
By Proposition \ref{FadellHusseini2},
there exists a compact set $V_{m,k}$ in $\tilde B_m$ such that
\[
\cat_{\tilde B_m,\tilde\Sigma_m}(V_{m,k}) \ge k .
\]
Furthermore, we can construct a retraction from $Z_m$ onto $\tilde B_m$,
and an homeomorphism from $Z_m$ onto itself so that Proposition \ref{retratto}
applies and the proof follows from Remark \ref{indippartespaziale}
(for more details, see \cite[Lemma 7.4]{CMS00}).
\end{proof}

\begin{proof}[Proof of Theorem \ref{maintheoremsplitting}.]
Assumed that $\mo$ is a 1--connected manifold (it can be done without loss of generality,
 due to the nature of the bounds \eqref{alfa}--\eqref{alfa2})
and $\eps$ small enough, from Lemmas \ref{PSlemma} and
\ref{lemmaquattro} it follows that Theorem \ref{teomolteplicita}
can be applied to each functional $f_{\eps,m}$ in $Z_m$.
Hence, there exists a sequence of critical points
$(z^k_{\eps,m})_{k} \subset Z_m$ of $f_{\eps,m}$ such that
\[
f_{\eps,m}(z^k_{\eps,m}) \ge \inf f_{\eps,m}(S)\ ,\quad
\lim_{k \to +\infty} f_{\eps,m}(z^k_{\eps,m}) = \sup f_{\eps,m}(Z_m)
= +\infty .
\]
Moreover, carefull estimates allow one to prove that
for each $k \ge 1$ there exists a constant $\gamma_k > 0$,
independent of $\eps$ and $m$, such that
\[
f_{\eps,m}(z^k_{\eps,m}) \le \gamma_m ,
\]
while for each $c > 0$ there exists $k_c \in \N$, independent
of $\eps$ and $m$, such that
\[
f_{\eps,m}(z^k_{\eps,m}) \ge \lambda c - \Delta_t^2 N\quad \hbox{for all $k \ge k_c$}
\]
(being $\lambda$, $N$ as in the given hypotheses).\\
Whence, Lemmas \ref{noeps} and \ref{galerkin} imply that for any $k \ge k_c$
there exists a critical point $z^k \in \Omega^1(p,q;\M)$ of $f$ such that
\[
\lambda c - \Delta_t^2 N \le f(z^k) \le \gamma_k
\]
and the end of the proof of Theorem \ref{maintheoremsplitting}
follows from the arbitrariness of $c > 0$.
\end{proof}

\begin{remark}
To our knowledge, up to now  no better result of geodesic
connectedness has been obtained in the orthogonal splitting case.
We emphasize that the assumptions in Theorem
\ref{maintheoremsplitting} imply global hyperbolicity \cite{Sa97}.
Even though every globally hyperbolic spacetime is orthogonal
splitting (\cite{BeSa05}, see also \cite{MS}),  it may not satisfy
assumptions in Theorem \ref{maintheoremsplitting} (as there are
non-geodesically connected counterexamples, including those in the
static case). So,
 following the
careful approaches introduced in the stationary case, it should be
interesting to obtain a result similar to that one of Theorem
\ref{maintheoremsplitting} also under weaker  assumptions on the
metric or under intrisinsic hypotheses more related to the
geometry of the manifold.
\end{remark}


{\footnotesize \printindex }

\end{document}